\documentclass[11pt, draft]{amsart}
\usepackage{amsmath,amsthm, amscd, amssymb, amsfonts}
\usepackage[dvips]{graphics}
\usepackage{hhline}
\usepackage[OT2,OT1]{fontenc}
\usepackage[all]{xy}
\usepackage{eucal}

\usepackage[dvips, dvipsnames, usenames]{color}

\renewcommand{\_}[1]{_{\left( #1 \right)}}
\renewcommand{\^}[1]{^{\left( #1 \right)}}
\newcommand{\ot}{{\otimes}}

\newcommand{\sE}{s_{i,M}}
\newcommand{\trid}{\triangleright}

\newcommand{\abofet}{{\ac(\st, \Oc_2^3, \lambda)}}

\newcommand{\VP}{\widetilde{V}}

\newcommand{\pl}[1]{\partial^L_{#1}}
\newcommand{\pr}[1]{\partial^R_{#1}}

\newcommand{\plb}[1]{\overline\partial^L_{#1}}

\newcommand{\Ind}{\operatorname{Ind}}

\newcommand{\ev}{\operatorname{ev}}

\newcommand{\supp}{\operatorname{supp}}

\newcommand{\co}{\operatorname{co}}
\newcommand{\coev}{\operatorname{coev}}

\newcommand{\pro}{\operatorname{pr}}
\newcommand{\trasp}{\tau}

\newcommand{\grd}{\operatorname{gr-dual} }

\newcommand{\vd}{{}^*V}

\newcommand\cyr{%
\renewcommand\rmdefault{wncyr}%
\renewcommand\sfdefault{wncyss}%
\renewcommand\encodingdefault{OT2}%
\normalfont \selectfont} \DeclareTextFontCommand{\textcyr}{\cyr}

\newcommand{\wo}{\Wg_0(M)}

\newcommand{\ac}{{\mathcal A}}
\newcommand{\toba}{{\mathcal B}}
\newcommand{\C}{{\mathcal C}}
\newcommand{\D}{{\mathcal D}}

\newcommand{\K}{{\mathcal K}}

\newcommand{\m}{\mathcal{M}}

\newcommand{\Oc}{{\mathcal O}}
\newcommand{\Pc}{{\mathcal P}}

\newcommand{\rfl}{{\mathcal R}}
\newcommand{\Ss}{{\mathcal S}}

\newcommand{\Wg}{{\mathcal W}}

\newcommand\Bg{\mathfrak B}

\newcommand{\J}{{\mathfrak J}}
\newcommand{\Kb}{\mathfrak K}

\newcommand{\W}{{\mathfrak W}}

\newcommand\g{\mathfrak{g}}
\newcommand\h{\mathfrak{h}}

\newcommand{\ku}{\Bbbk}
\newcommand{\Cc}{{\mathbb C}}
\newcommand{\Ib}{{\mathbb I}}
\newcommand{\Mb}{{\mathbb M}}
\newcommand{\N}{{\mathbb N}}

\newcommand\st{\mathbb S_3}
\newcommand\sk{\mathbb S_4}

\newcommand\sn{\mathbb S_n}
\newcommand\dn{\mathbb D_n}
\newcommand\dt{\mathbb D_3}

\newcommand{\Z}{{\mathbb Z}}
\def\zt{\Z^{\theta}}

\newcommand{\ydh}{{}^H_H\mathcal{YD}}

\newcommand{\ydvh}{{}^{\ac(V)}_{\ac(V)}\mathcal{YD}}
\newcommand{\ydrh}{{}^{R\# H}_{R\# H}\mathcal{YD}}

\newcommand{\ydg}{{}^{\Cc G}_{\Cc G}\mathcal{YD}}

\newcommand{\Alg}{\operatorname{Alg}}

\newcommand{\End}{\operatorname{End}}
\newcommand{\Aut}{\operatorname{Aut}}

\newcommand\sgn{\operatorname{sgn}}
\newcommand\ad{\operatorname{ad}}
\newcommand\Hom{\operatorname{Hom}}

\newcommand{\lmax}{L^{\rm max}}
\newcommand{\lt}{{\langle Z\rangle}}

\newcommand{\rsys}{\boldsymbol{\Delta }^{re}}
\newcommand{\fd}{finite-dimensional}
\newcommand{\infd}{infinite-dimensional}

\newcommand{\cou}{\varepsilon }

\numberwithin{equation}{section}\theoremstyle{plain}

\newtheorem{thmintro}{Theorem}

\newtheorem{theorem}{Theorem}[section]
\newtheorem{lema}[theorem]{Lemma}
\newtheorem{mlema}[theorem]{Main Lemma}
\newtheorem{cor}[theorem]{Corollary}

\newtheorem{prop}[theorem]{Proposition}

\theoremstyle{definition}
\newtheorem{definition}[theorem]{Definition}

\newtheorem{exa}[theorem]{Examples}

\theoremstyle{remark}
\newtheorem{obs}[theorem]{Remark}

\newcommand\id{\operatorname{id}}

\newcommand\op{\operatorname{op}}
\newcommand\cop{\operatorname{cop}}

\newcommand\sw{\operatorname{\textmd{\rm o}}}

\def\pf{\begin{proof}}
\def\epf{\end{proof}}

\theoremstyle{remark}

\begin{document}

\renewcommand{\baselinestretch}{1.2}

\thispagestyle{empty}
\title[The Nichols algebra of a semisimple Yetter-Drinfeld module]
{The Nichols algebra of a semisimple Yetter-Drinfeld module}
\author[andruskiewitsch, heckenberger and schneider]{Nicol\'{a}s Andruskiewitsch}
\address{Facultad de Matem\'{a}tica, Astronom\'{i}a y
F\'{i}sica\\\newline\indent Universidad Nacional de C\'{o}rdoba \\
\newline\indent CIEM - CONICET
\\ \newline\indent (5000) Ciudad Universitaria \\ C\'{o}rdoba
\\Argentina} \email{andrus@famaf.unc.edu.ar}
\author[]{Istv\'an Heckenberger}
\address{Mathematisches Institut\\ \newline\indent  Universit\"at M\"unchen
\\ \newline\indent  Theresienstr. 39,
D-80333 Munich, Germany}
\email{i.heckenberger@googlemail.com}
\author[]{Hans-J\"urgen Schneider}
\address{Mathematisches Institut\\ \newline\indent  Universit\"at M\"unchen
\\ \newline\indent  Theresienstr. 39,
D-80333 Munich, Germany}
\email{Hans-Juergen.Schneider@mathematik.uni-muenchen.de}
\thanks{The work of N. A.
was partially supported by Ag. Cba. Ciencia, CONICET, Foncyt and
Secyt (UNC).
The work of I.H. was partially supported by DFG within a Heisenberg fellowship
at the University of Munich}

\begin{abstract} We study the Nichols algebra
of a semisimple Yetter-Drin\-feld module and introduce new invariants such as
real roots. The crucial ingredient is
a ``reflection'' in the class of such Nichols algebras.
We conclude
the classifications of \fd{} pointed Hopf algebras
over $\st$, and of \fd{} Nichols algebras over $\sk$.
\end{abstract}

\maketitle

\setcounter{tocdepth}{2} \tableofcontents

\section*{Introduction}

{\bf 1.} Although semisimple complex Lie algebras $\mathfrak{g}$ cannot be deformed there are highly interesting $q$-deformations $U_q(\mathfrak{g})$ of their enveloping algebras {\em as Hopf algebras} with generic $q$ introduced by Drinfeld and Jimbo around 1985. 

As an algebra, $U_q(\mathfrak{g})$ is generated by elements $E_i,F_i,K_i^{\pm}, 1 \leq i \leq n$. The Hopf algebra $U_q(\mathfrak{g})$ is determined by the +-part $U_q^+(\mathfrak{g})= k\langle E_1,\dots,E_n \rangle$ since $U_q(\mathfrak{g})$ is essentially  a Drinfeld double of $U_q^+(\mathfrak{g})$. The algebra $U_q^+(\mathfrak{g})$ has a very easy and beautiful description as the Nichols algebra (or quantum symmetric algebra) $\mathcal{B}(W)$ of a finite-dimensional vector space 
\begin{equation}\label{example}
W=\mathbb{C} E_1 \oplus \cdots \oplus \mathbb{C} E_n
\end{equation}
together with a grading and an action of a free abelian group $G$ with basis $K_1,\dots,K_n$. Each $E_i$ has degree $K_i$ and the action of $G$ is given by 
\begin{equation}\label{action}
K_i E_jK_i^{-1}=q^{d_ia_{ij}} \text{ for all }i,j.
\end{equation}
Here $(d_ia_{ij})$ is the symmetrized Cartan matrix. Thus $W$ has the structure of a Yetter-Drinfeld module over the group algebra $k[G]$, and the vector spaces $\mathbb{C}E_i$ are one-dimensional Yetter-Drinfeld submodules. See Section 1 for the definition of Yetter-Drinfeld modules.

In the same way Nichols algebras also determine the small quantum groups $u_q(\mathfrak{g})$, $q$ a root of unity, introduced by Lusztig, and the generalizations of the quantum groups $U_q(\mathfrak{g})$ to Kac-Moody Lie algebras, see \cite{L, R, R2, Sbg}.

\medbreak{\bf 2.} Nichols algebras appeared in the work of Nichols \cite{N}. They are defined for Yetter-Drinfeld modules $W$ over any Hopf algebra $H$ (with bijective antipode) instead of the group algebra $H=k[G]$. The category of Yetter-Drinfeld modules over $H$ is braided, and the Nichols algebra can be defined by the following universal property: The tensor algebra $T(W)$ is a braided Hopf algebra where the elements of $W$ are primitive. Then
\begin{equation}\label{N}
\mathcal{B}(W) = T(W)/I_W,
\end{equation}
where $I_W$ is the largest coideal of $T(W)$ spanned by elements of $\mathbb{N}$-degree $\geq 2$.

The smash product $\mathcal{B}(W) \# H$ (called bosonization) is a Hopf algebra in the usual sense, and to  understand Nichols algebras of Yetter-Drinfeld modules  in general is of fundamental importance for the general theory of Hopf algebras. Nichols algebras form a crucial part of the $\mathbb{N}$-graded Hopf algebra associated to the coradical filtration of a Hopf algebra whose  coradical is a Hopf subalgebra \cite{AS-jalg}. An important class of such Hopf algebras are pointed Hopf algebras, that is Hopf algebras where the coradical is a group algebra (or equivalently, where all the simple comodules are one-dimensional). The quantum groups $U_q(\mathfrak{g})$ and all their variants are pointed.

\medbreak The definition of the Nichols algebra is easy. The inherent conceptual difficulty of understanding Nichols algebras is their very indirect definition by a universal property. In general there is no method to actually determine the Nichols algebra $\mathcal{B}(W)$ for a given $W$, for example to calculate the dimensions of the $\mathbb{N}$-homogeneous components of $\mathcal{B}(W)$ or to compute the defining relations, that is to compute generators of the unknown ideal $I_W$.

The relations of the Nichols algebra of the Yetter-Drinfeld module \eqref{example} are the quantized Serre relations, see \cite[33.1.5]{L} for a proof of this deep result.

\medbreak{\bf 3.} During the last few years several classification results for Hopf algebras were obtained based on the theory of Nichols algebras and following the procedure proposed in \cite{AS-jalg}. This program has been particulary successful for finite-dimensional pointed Hopf algebras with abelian
group of group-like elements \cite{AS-05}.

Let $H$ be a finite-dimensional cosemisimple complex Hopf algebra
and let us consider a finite-dimensional Hopf algebra $A$ such that its coradical
$A_0$ is a Hopf algebra isomorphic to $H$. To solve the problem of
classifying all such Hopf algebras $A$, we have to
address two fundamental questions: Given a Yetter-Drinfeld module
$W$ over $H$,
\begin{enumerate}
     \item[(a)]\label{questions}
     decide when $\dim \toba(W) < \infty$, and

    \item[(b)] describe a suitable set of defining relations of
    $\toba(W)$.
\end{enumerate}

Now, since $H$ is cosemisimple, it is also semisimple \cite{LR};
then the category $\ydh$  of Yetter-Drinfeld modules over $H$ is
semisimple \cite{Ramq}. Therefore we just need to consider the
questions (a) and (b) in the following cases:

\begin{enumerate}
     \item[(i)] when $W$ is an irreducible Yetter-Drinfeld module, and

    \item[(ii)]  when $W = V_1 \oplus \dots \oplus V_\theta$ is a direct sum
of irreducible Yetter-Drinfeld modules, under the assumption that the answers
to questions (a) and (b) are known for $V_1, \dots, V_\theta$.
\end{enumerate}

As applications of the main results of the present paper we obtain new information in case (ii) when not all the $V_i$ are one-dimensional.

\medbreak  If $V = \Cc v$ is a one-dimensional
Yetter-Drinfeld submodule over $H$, then it determines a
group-like element $g\in G(H)$ and a character $\chi\in \Alg(H,
\Cc)$ defining the coaction and the action of $H$. Let $q =
\chi(g)$. The Nichols algebra of $V$ is easy to determine: it is
either the polynomial algebra $\Cc[v]$, when $q = 1$ or is not a
root of 1, or else it is the truncated polynomial algebra
$\Cc[v]/(v^N)$, when $q$ is a root of 1 of order $N >1$. In other
words, questions (a) and (b) have completely satisfactory answers
in this case.

\medbreak Assume next that $W = \Cc v_1 \oplus \dots \oplus \Cc
v_\theta $ is a direct sum of \emph{one-dimensional}
Yetter-Drinfeld submodules over $H$. Let $g_i\in G(H)$ and
$\chi_i\in \Alg(H, \mathbb{C})$ determined by the submodule $\ku v_i$ as
above. Let $q_{ij} = \chi_j(g_i)$, $1\le i,j\le \theta$. Note that in the classical situation of the Yetter-Drinfeld module \eqref{example} for $q$ a root of unity, $q_{ij} = q^{d_ia_{ij}}$. The
Nichols algebra of $W$ can be viewed as a ``gluing'' of the
various Nichols subalgebras $\toba(\mathbb{C} v_i)$ along the generalized
Dynkin diagram with vertices $1, \dots, \theta$; there is a line
joining the vertices $i$ and $j$ if $q_{ij}q_{ji} \neq 1$, and
then the line is labelled by the scalar $q_{ij}q_{ji}$, resembling
the classical Killing-Cartan classification of semisimple Lie
algebras.

\medbreak Assume moreover that $W$ is of Cartan type, that is,  there exist
$a_{ij}\in \mathbb Z$ such that $q_{ij}q_{ji} = q_{ii}^{a_{ij}}$
for any $i\neq j$; the classical example of Cartan type is $q_{ij}=q^{d_ia_{ij}}$ for all $i,j$ where $(d_ia_{ij})$ is the symmetrized Cartan matrix. Then $A =
(a_{ij})$ is a generalized Cartan matrix, and we have complete
answers to questions (a) and (b) above:

\begin{enumerate}

   \item[(a)] \emph{$\dim\toba(W) < \infty$ if and only if $A$ is of finite
    type.}

    \item[(b)] \emph{If $\dim\toba(W) < \infty$, then the ideal of relations of $\toba(W)$ is generated by the quantum
    Serre relations and appropriate powers of the root vectors.}
    \end{enumerate}

These results were proved in \cite{AS-adv} under some restrictions
on the orders of the $q_{ij}$'s by reduction to the theory of
quantum groups. Part (a) was shown without any restriction in
\cite{H3}. 

\medbreak The classification of the matrices $(q_{ij})$ whose
corresponding Nichols algebras are \fd{} is given in \cite{H2}. In
general, one Cartan matrix is not sufficient to describe the Nichols algebra, a family of generalized Cartan matrices is needed. The main instrument to control them is the Weyl groupoid --
introduced already in \cite{H3}. As for the defining relations,
these are not yet known, except in the standard case \cite{Ag} when all the Cartan matrices are the same, and in the general case of rank two, that is $\theta=2$ \cite{H4}. Their precise description  is more delicate than (b) above.

\medbreak {\bf 4.} Let now $W = V_1 \oplus \dots \oplus V_\theta$
be an arbitrary direct sum of irreducible Yetter-Drinfeld modules.
In analogy with the situation of Cartan type, it was proposed to
consider the $V_i$'s as ``fat points'' of a generalized Dynkin
diagram (or some kind of generalized Cartan matrix) and, assuming
the knowledge of the Nichols algebras $\toba(V_i)$, to describe
the Nichols algebra $\toba(W)$ as a ``gluing'' of the various
Nichols subalgebras $\toba(V_i)$ along it \cite[p. 41]{bariloche}.
Because of \cite{H2}, it is clear that just one generalized Cartan
matrix would not be enough, and that we would need to attach to
our $W$ a collection of generalized Cartan matrices. This is what
we do in the present paper.

\medbreak {\bf 5.} Let us now proceed with a detailed description
of our results which in fact hold in a much more general context. Let
$\ku$ be an arbitrary field and let $H$ be any Hopf algebra with
bijective antipode. Let 
\begin{equation}
W = V_1 \oplus \dots \oplus V_\theta
\end{equation}
be a direct sum of \fd{} irreducible Yetter-Drinfeld modules over
$H$. Assume for simplicity that the adjoint action of $\toba(W)$
on itself is locally finite. We fix an index $i$, $1\le i\le
\theta $. Let
\begin{align}\label{cartanmatrices}
  a_{ij} &= 1 - \text{top degree of } \ad\, \toba(V_i)(V_j), \quad i\neq
  j,\text{ and } a_{ii} = 2.
\end{align}
Then $(a_{ij})_{1\le i,j\le \theta}$ is a generalized Cartan
matrix attached to $W$; note that a version of the quantum Serre
relations holds by definition. We define $V'_i=V_i^*$, $V'_j$ as the top homogeneous
component of $\ad\, \toba(V_i)(V_j)$ if $i\neq j$, and
\begin{align}\label{V}
  W' &= V'_1 \oplus \dots \oplus V'_\theta.
\end{align}

Let $\K$ be the algebra of coinvariant elements of $\toba(W)$ with
respect to the right coaction of $\toba(V_i)$, and $\#$ denotes
the smash product introduced in Definition \ref{def:ksmash}. Our
first main result is the key step for the construction of the
family of Cartan matrices generalizing  \cite[Prop. 1]{H3} where all $V_i$ are one-dimensional.

\begin{thmintro}\label{theo:first} There is an isomorphism
\begin{align}\label{intro:iso}
  \toba(W')\simeq \K\#\toba(V_i^*).
\end{align}

In particular,  if $\dim \toba(W) < \infty$, then $\dim
\toba(W)=\dim \toba(W')$.
\end{thmintro}

The assignment $W\mapsto W'$ is a generalized $i$-th reflection. Theorem \ref{theo:first} allows to find by iterated reflections a class of new Nichols algebras $\mathcal{B}(W')$ of the same dimension as $\mathcal{B}(W)$. This defines an equivalence relation between non-isomorphic Nichols algebras. The collection of generalized Cartan matrices
we are looking for consists of the generalized Cartan matrices of the Nichols algebras in the equivalence class of $\mathcal{B}(W)$. We are now in a position to define real roots
of $\mathcal{B}(W)$, see Section \ref{ss:Weylgroupoid}.

In order to prove Theorem \ref{theo:first}, we have to overcome several difficulties. The proof of \cite[Prop. 1]{H3} depends on the existence of PBW-bases shown  by Kharchenko \cite{Kh}. In our case where not all the $V_i$ are one-dimensional such bases do not exist in general. Another difficulty in the general case is to prove irreducibility of the Yetter-Drinfeld modules $V_i'$ in \eqref{V}. Our proof of Theorem 1 can
not rely on the usual characterization of a Nichols algebra as a
braided Hopf algebra with special properties, because it does not
seem possible to describe the comultiplication of
$\K\#\toba(V_i^*)$ explicitly. Instead, we use a new
characterization of Nichols algebras in terms of braided
derivations, see Theorem \ref{theo:quotients-with-derivations}.
This new characterization is a powerful tool to deal with Nichols
algebras; we expect many applications of it.


\medbreak {\bf 6.} Having defined the collection of
generalized Cartan matrices and reflections attached to our $W$, 
the following questions arise:

\begin{enumerate}

   \item[(A)] To develop a theory of generalized root systems that
   correspond to our collections of generalized Cartan matrices,
   including classifications of suitable classes of them.

    \item[(B)] To obtain answers to questions (a) and
    (b) in page \pageref{questions} on the Nichols algebra $\toba(W)$ from the structure of
    its generalized root systems.
    \end{enumerate}

These matters are out of the scope of the present paper. In \cite{HS} the generalized root system of $\mathcal{B}(W)$ is defined (under the restriction that $H$ is semisimple or more generally that all finite tensor powers of $W$ are semisimple). These root systems satisfy the axioms introduced in \cite{HY} and studied in \cite{CH08a, CH08b}.

We present however a partial answer to question (B). Let us say that
$W$ is \emph{standard} if the generalized Cartan matrix
$(a'_{ij})_{1\le i,j\le \theta }$ corresponding to $W'$ coincides
with the generalized Cartan matrix $(a_{ij})_{1\le i,j\le \theta
}$ corresponding to $W$, for all $W'$ obtained from $W$ by
finitely many reflections.

\begin{thmintro}\label{theo:second} If $W$ is standard and $\dim \toba(W)< \infty$, then the
generalized Cartan matrix is of finite type.
\end{thmintro}

See Thm.~\ref{theorem:gpd-nichols-finite}. By \cite[Corollary 7.4]{HS} the
converse of Theorem 2 is true, that is, if $W$ is standard, $\dim
\toba(V_i')<\infty $ for all $i$ and all $W'$ obtained from $W$ by iterated reflections, and if $(a_{ij})_{1\le i,j\le \theta
}$ is of finite type, then $\dim \toba(W)<\infty $. Using the results of the present paper a necessary and sufficient criterion for $\dim \toba(W)< \infty$ is given in the general non-standard case in \cite[Theorem 7.3]{HS}.

\medbreak {\bf 7.}  There is at the present moment no general
method to deal with questions (a) and (b) for \emph{irreducible}
Yetter-Drinfeld modules over a finite non-abelian group. In fact,
we know very few examples with finite dimension. The first
examples, calculated in 1995, correspond to the transpositions in
$\sn$, $n=3,4,5$ \cite{MS}. As an application of Theorem
\ref{theo:second} for $\st$ and $\sk$, we prove that $\mathcal{B}(W)$ is infinite-dimensional if $W$ is not irreducible. (In \cite{HS}  this result is generalized to all finite simple groups and to all symmetric groups $\sn$, $n \geq 3$.) This allows to conclude the classifications of \fd{}
pointed Hopf algebras over $\st$, Theorem \ref{theo:s3}, and of
\fd{} Nichols algebras over $\sk$, Theorem \ref{thm:s4}. The group
$\st$ is the first non-abelian group $G$ where the classification
of \fd{} pointed Hopf algebras with coradical $\ku G$ is known,
and where a Hopf algebra other than the group algebra exists.
Recently, some groups that admit no \fd{} pointed Hopf algebra
except the group algebra were found: $\mathbb{A}_n$, $n\geq 5$,
$n\neq 6$ \cite{AF2,AFGV} and
$SL(2,q)$ with $q$ even \cite{FGV}.
 Theorems \ref{theo:s3} and \ref{thm:s4} can be rephrased
in terms of racks, giving rise to new techniques to establish that
some Nichols algebras have infinite dimension \cite{AF3}. These
techniques have been applied in \cite{AFZ, AFGV, AFGVe}.

\medbreak {\bf 8.}  The paper is organized in four sections,
besides this introduction. In Sect.~\ref{section:preliminaries} we
collect several well-known results that will be used later on. In
Sect.~\ref{section:qdo} we use quantum differential operators to
give a new characterization of Nichols algebras.
Sect.~\ref{section:weyl} is the bulk of the paper: We construct
the reflection of a semisimple Yetter-Drinfeld module satisfying
some hypothesis (for instance, having \fd{} Nichols algebra),
discuss the notion of ``standard'' modules, and prove our main
theorems. In Sect.~\ref{section:appl} we state a few general
consequences of the theory in the previous sections, and then
prove the classification results for $\st$ and $\sk$ alluded
above. We also include a result on Nichols algebras over the
dihedral group $\dn$ with $n$ odd.

\medbreak In the paper $H$ denotes a Hopf algebra with bijective
antipode $\Ss $.

\section{Preliminaries}\label{section:preliminaries}

\subsection{Notation}\label{subsection:notation}

Let $\ku$ be a field. All vector spaces, algebras, coalgebras,
Hopf algebras, unadorned tensor products and unadorned Hom spaces
are over $\ku$. If $V$ is a vector space and $n\in \N$, then
$V^{\ot n}$ or $T^n(V)$ denote the $n$-fold tensor product of $V$
with itself. We use the notation $\langle\,,\,\rangle: \Hom(V,
\ku) \times V \to \ku$ for the standard evaluation. We identify
$\Hom(V, \ku)\ot \Hom(V, \ku)$ with a subspace of $\Hom(V\ot V,
\ku)$ by the recipe
$$
\langle f\ot g, v\ot w\rangle = \langle f, w\rangle\langle g,
v\rangle
$$
for $f, g\in \Hom(V, \ku)$, $v, w\in V$. Consequently, we identify
$\Hom(V, \ku)^{\ot n}$ with a subspace of $\Hom(V^{\ot n}, \ku)$,
$n\in \N$, via
\begin{equation}\label{eqn:vndual}
\langle f_n\ot \dots \ot f_1, v_1\ot \dots \ot v_n\rangle =
\prod_{1\le i \le n} \langle f_i, v_i\rangle,
\end{equation}
for $f_1, \dots, f_n\in \Hom(V, \ku)$, $v_1, \dots, v_n\in V$.

\medbreak Let $\theta\in \N$ and let $\Ib = \{1, \dots, \theta\}$.
Let $V = \oplus_{\alpha\in \zt} V_\alpha$ be a
$\zt$-graded vector space. If $\alpha = (n_1, \dots,
n_{\theta})\in \zt$, then let $\pro_\alpha=
\pro_{n_1, \dots, n_{\theta}}: V\to V_\alpha$ denote the projection
associated to this direct sum. We identify $\Hom(V_\alpha, \ku)$
with a subspace of $\Hom(V, \ku)$ via the transpose of
$\pro_\alpha$. The graded dual of $V$ is
\begin{align}
  V^{\grd} = \oplus_{\alpha\in \zt} \Hom(V_\alpha, \ku)
  \subset \Hom(V, \ku).
  \label{eq:grdual}
\end{align}
If $V = \oplus_{\alpha\in \zt} V_{\alpha}$ is a
$\zt$-graded vector space, then the support of $V$ is $\supp
V := \{\alpha\in \zt \,|\, V_{\alpha} \neq 0\}$.

\medbreak  Let $C$ be a coassociative coalgebra. Let
$\Delta^n: C \to C^{\ot (n+1)}$ denote the $n$-th
iterated comultiplication of $C$. Let $G(C)$ denote the set of
group-like elements of $C$. If $g, h\in G(C)$, then let
${\Pc}_{g,h}(C)$ denote the space
$\{x \in C\,|\, \Delta(x) = g \otimes x + x \otimes h\}$
of $g$, $h$ skew-primitive elements of $C$. If $C$
is a braided bialgebra, then ${\Pc}(C) := {\Pc}_{1,1}(C)$. The
category of left (resp. right) $C$-comodules is denoted ${}^C\m$,
resp. $\m^C$. We use Sweedler's notation for the comultiplication
of $C$: If $x\in C$, then $\Delta(x) = x\_{1}\otimes x\_{2}$.
Similarly, the coaction of a left $C$-comodule $M$ is denoted
$\delta(m) = m\_{-1} \ot m\_{0}\in C\ot M$, $m\in M$.

\begin{obs}\label{rem:C*alg}
  The dual vector space $C^* = \Hom(C, \ku)$ is an algebra
with the convolution product: $\langle fg, c\rangle = \langle g,
c\_{1}\rangle\, \langle f, c\_{2}\rangle$, cf. \eqref{eqn:vndual},
for $f,g\in C^*$, $c\in C$. The reader should be warned that
usually one writes $C^*{}^{\op }$ for this algebra, see
\cite[Sect.\,1.4.1]{Mo}. With our convention -- forced by
\eqref{eqn:vndual} -- any left $C$-comodule becomes a \emph{left}
$C^*$-module by
\begin{equation}\label{transpose-action}
f\cdot m = \langle f, m\_{-1}\rangle m\_{0},
\end{equation}
$f\in C^*$, $m\in M$. Indeed, if also $g\in C^*$, then
\begin{align*}
f\cdot (g\cdot m) &= \langle g, m\_{-1}\rangle f\cdot  m\_{0} =
\langle f, m\_{-1}\rangle \langle g, m\_{-2}\rangle m\_{0} \\
&= \langle fg, m\_{-1}\rangle m\_{0}  = (fg)\cdot m.
\end{align*}
\end{obs}

Recall that a graded coalgebra is a coalgebra $C$
provided with a grading $C = \oplus_{m\in \N_0} C^m$ such that
$\Delta(C^m) \subset \oplus_{i+j =m} C^i \ot C^j$. Then the graded
dual $C^{\grd}$ is a subalgebra of $C^*$.

Let $\Delta_{i,j}: C^{m} \to C^i \ot C^j$ denote the composition
$\pro_{i,j}\Delta$, where $m = i+ j$. More generally, if $i_1,
\dots, i_n\in \N_0$ and $i_1 + \dots + i_n = m$, then
$\Delta_{i_1, \dots, i_n}$ is the composition  $\pro_{i_1, \dots,
i_n}\Delta^{n-1}$:
\begin{equation}\label{deltaij}
\xymatrix{C^{m} \ar[d]_{\Delta_{i_1, \dots, i_n}}
\ar[0,1]^-{\Delta^{n-1}} &   \oplus_{j_1 + \dots + j_n = m}C^{j_1}
\ot \dots \ot C^{j_n} \ar@{>>}[1,-1]^{\qquad\pro_{i_1, \dots,
i_n}}
\\ C^{i_1} \ot \dots \ot C^{i_n}. &}
\end{equation}

\begin{obs}\label{subcom-gen}
Let $C$ be a coalgebra, let $M\in {}^{C}\m$ and let $Z\subset M$
be a vector subspace. Then the subcomodule generated by $Z$ is
\begin{equation}\label{subcom-gen-u}
C^*\cdot Z = \ku\text{-span of }\{\langle f, z\_{-1}\rangle\,
z\_{0}\,|\, z\in Z, \, f\in C^*\}.
\end{equation}
If $C = \oplus_{m\in \N_0} C^n$ is a graded coalgebra, then
\begin{equation}\label{subcom-gen-gr}
C^*\cdot Z = \ku\text{-span of }\{\langle f, z\_{-1}\rangle\,
z\_{0}\,|\, z\in Z, \, f\in C^{\grd}\}.
\end{equation}
\end{obs}

\pf Clearly, \eqref{subcom-gen-u} is the subcomodule generated by
$Z$. Assume that $\dim Z< \infty$. Then there exists $m\in \N$
such that $\delta(Z) \subset \oplus_{0\le n\le m} C^n\ot M$.
Therefore, in \eqref{subcom-gen-u} it suffices to take
$$f\in\left(\oplus_{ n> m} C^n\right)^\perp \simeq \left(\oplus_{0\le n\le m} C^n\right)^*
\subset \oplus_{n\ge 0} \left(C^n\right)^*.$$ If $\dim Z$ is
arbitrary, then $$C^*\cdot Z = C^*\cdot\Big(\sum_{Z'\subset Z:
\dim Z' < \infty} Z'\Big) = \sum_{Z'\subset Z\,|\,\dim Z' < \infty}
\big(C^*\cdot Z'\big),$$ proving the assertion. \epf

\subsection{Yetter-Drinfeld modules}\label{subsection:yd}
Our reference for the theory of Hopf algebras is \cite{Mo}. Recall that
$H$ is a Hopf algebra with bijective antipode $\Ss$. The adjoint
representation of $H$ on itself is the algebra map $\ad: H\to \End
H$, $\ad x (y) = x\_{1} y \Ss(x\_{2})$, $x,y\in H$. Then
\begin{equation}\label{ad}
\ad x (yy') = \ad (x\_{1}) (y) \ad(x\_{2})(y'),
\end{equation}
$x,y, y'\in H$. That is, $H$ is a left $H$-module algebra via the
adjoint.

\medbreak Let $\ydh$ be the category of Yetter-Drinfeld modules
over $H$; $V\in \ydh$ is a left $H$-module and a left $H$-comodule
such that
\begin{equation}\label{yd}
\delta(h\cdot x) = h\_{1}x\_{-1}\Ss(h\_{3})\ot h\_{2}\cdot x\_{0},
\end{equation}
$h\in H$, $x\in V$. It is well-known that $\ydh$ is a braided
tensor category, with braiding $c_{V, W}: V\ot W \to W\ot V$,
$c_{V, W}(v\ot w) = v\_{-1}\cdot w\ot v\_{0}$, $V, W\in \ydh$,
$v\in V$, $w\in W$. We record that the inverse braiding is given
by
\begin{equation}\label{inversebraiding}
  c^{-1}_{V, W}(v\ot w) = w\_{0}\ot \Ss^{-1}(w\_{-1})\cdot v,
\end{equation}
$V, W\in \ydh$, $v\in V$, $w\in W$.

\begin{obs}\label{ydsubm}
Let $V\in \ydh$.

(i) If $U\subset V$ is an $H$-submodule, then the subcomodule
$H^* \cdot U$ generated by $U$ is a Yetter-Drinfeld submodule of
$V$.

(ii) If $T\subset V$ is an $H$-subcomodule, then the submodule
$H\cdot T$ generated by $T$ is a Yetter-Drinfeld submodule of $V$.
\end{obs}

\pf (i).  If $u\in U$, $f\in H^*$ and $h\in H$, then
$$
h\cdot(\langle f, u\_{-1}\rangle\, u\_{0}) = \left\langle f, \Ss(
h\_{1})(h\_{2}\cdot u)\_{-1} h\_{3}\right\rangle\, (h\_{2}\cdot
u)\_{0} \in H^*\cdot U
$$
by \eqref{yd}.  (ii) is also a direct consequence of \eqref{yd}.
\epf

Let $V\in \ydh$ be \fd. The left and right duals of
$V$ are respectively denoted $^*V$ and $V^*$. As vector spaces,
$^*V = V^* = \Hom(V, \ku)$. Their structures of Yetter-Drinfeld
modules are determined by requiring that the following natural
maps are morphisms in $\ydh$:
\begin{align*}
&\ev: V^* \ot V\to \ku, & &\coev: \ku\to V \ot V^*, \\
&\ev: V \ot {}^*V\to \ku, & &\coev: \ku\to {}^*V \ot V,
\end{align*}
cf. \cite[Def. 2.1.1]{BK}. Thus $V^*$ has action and coaction
given by
\begin{align}\label{vd1}
\langle h\cdot f, v\rangle &= \langle f, \Ss(h)\cdot v\rangle,
\\\label{vd2} f\_{-1}\langle f\_{0}, v\rangle &= \Ss^{-1}(v\_{-1})\langle f, v\_{0}\rangle,
\end{align}
$f \in V^*$, $v\in V$. Albeit evident, we record that \eqref{vd2}
is equivalent to
\begin{align}\label{vd2bis}
\Ss(f\_{-1})\langle f\_{0}, v\rangle &= v\_{-1}\langle f,
v\_{0}\rangle,
\end{align}
$f \in V^*$, $v\in V$. Notice that \eqref{vd1}
provides $V^* = \Hom(V, \ku)$ with
an $H$-module
structure, regardless of whether $\dim V$ is finite or not.

It is easy to see that $T^n(V^*)$ is a Yetter-Drinfeld submodule of
$(T^n(V))^*$ via the identification \eqref{eqn:vndual}. Also, the
evaluation $\langle\,,\,\rangle: V^* \times V \to \ku$ satisfies
\begin{equation}\label{braided-duality}
\langle c_{V^*}(f\ot g), v\ot w\rangle = \langle f\ot g,
c_{V}(v\ot w)\rangle,
\end{equation}
$f, g\in V^*$, $v,w\in V$.

\pf We compute
\begin{align*} \langle c_{V^*}(f\ot g), v\ot
w\rangle &= \langle f\_{-1}\cdot g\ot f\_{0}, v\ot w\rangle =
\langle f\_{0}, v\rangle\langle f\_{-1}\cdot g, w\rangle
\\&= \langle f\_{0}, v\rangle\langle g, \Ss(f\_{-1})\cdot w\rangle
= \langle f, v\_{0}\rangle\langle g, v\_{-1}\cdot w\rangle
\\&= \langle f\ot g, v\_{-1}\cdot w\ot v\_{0}\rangle
= \langle f\ot g, c_{V}(v\ot w)\rangle.
\end{align*}
\epf

Analogously, $^*V$ has action and coaction given by $\langle v,
h\cdot f\rangle = \langle  \Ss^{-1}(h)\cdot v, f\rangle$, $f\_{-1} \langle
v, f\_{0}\rangle = \Ss(v\_{-1})\langle v\_{0}, f\rangle$, $f \in
\vd$, $v\in V$.

\begin{obs}\label{lema:double-dual}
  One has $V\simeq V^{**}$ for any \fd{} $V\in \ydh$
  \cite[(2.2.6)]{BK}. Explicitly, if we identify
  $V$ and $V^{**}$ as vector spaces via the map $v\mapsto \varphi _v$,
  where $\langle \varphi _v,f\rangle :=\langle f,v\rangle $ for all $f\in V^*$ and $v\in V$,
  then the isomorphism $\psi_V: V^{**}\to V$ in $\ydh$ and its inverse
  $\phi_V := \psi^{-1}_V$ are given by
\begin{align}\label{eqn:double-dual}
  \psi_V(\varphi _v) &=  \Ss^{-2}(v\_{-1})\cdot v\_0,
\\
  \label{eqn:double-dualbis}
  \phi_V(v) &=  \Ss((\varphi _v)\_{-1})\cdot (\varphi _v)\_0,
  \qquad v\in V.
\end{align}
Further, \eqref{vd1} and \eqref{vd2} imply that
\begin{align}
  \delta (\varphi _v)=&\Ss^{-2}(v\_{-1})\ot \varphi _{v\_0}, \notag \\
  \langle \phi _V(v), f\rangle =& \langle v\_{-1}\cdot f,v\_0\rangle .
  \label{eq:phiVv,f}
\end{align}
\end{obs}

\subsection{Smash coproduct}\label{subsection:smashcoproduct}
We shall need later the following
well-known facts.
Let $C\in {}^H\m$ be a left comodule coalgebra-- that
is, the comultiplication of $C$ is a comodule map. Let us denote
the comultiplication of $C$ by the following variation of
Sweedler's notation: If $c\in C$, then $\Delta(c) = c\^{1}\otimes
c\^{2}$. Let $C\# H$ be the corresponding smash coproduct: This is
the vector space $C\ot H$ (with generic element $c\#h$)  with
comultiplication
\begin{align} \Delta(c\# h) = c^{(1)} \# (c^{(2)})_{(-1)} h_{(1)} \otimes
  (c^{(2)})_{(0)} \#  h_{(2)},
  \label{eq:smashcopr}
\end{align}
$c\in C$, $h\in H$.
Let $p_C = \id\ot \cou: C\# H \to C$ and $p_H =
\cou\ot \id: C\# H \to H$ be the canonical coalgebra
projections. Let $\trasp: H\ot C \to C\ot H$ be given by
$$\trasp(h\ot c) = c\_{0}\ot \Ss^{-1}(c\_{-1})h,\qquad h\in H,
\, c\in C.$$

\begin{lema}\label{lema:smashcopr}
Let $M\in {}^{C\# H}\m$ with coaction $\delta_{C\# H}$.
Hence also $M\in {}^{C}\m$ with coaction
$\delta_{C} = (p_C\ot \id)\delta_{C\# H}$ and $M\in {}^{H}\m$ with
coaction $\delta_{H} = (p_H\ot \id)\delta_{C\# H}$. Then the
following hold.
\begin{enumerate}
  \item[(i)]
$\delta_{C\#H} = (\id \ot\delta_{H})\delta_{C} = (\trasp\ot \id)(\id
\ot\delta_{C})\delta_{H}$.
  \item[(ii)]
If $N\subset M$ is both a $C$-subcomodule and an
$H$-subcomodule, then it is a $C\# H$-subcomodule.
\item[(iii)]
If $Z\subset M$ is an $H$-subcomodule, then the
$C$-subcomodule generated by $Z$ is a $C\# H$-subcomodule.
\end{enumerate}
\end{lema}

\pf Let $m\in M$ and write $\delta_{C\# H}(m) = m\_{C, -1}\#
m\_{H, -1}\ot m\_{0}$. We spell out the coassociativity in this
notation:
\begin{multline}\label{coassociativity}
m\_{C, -1}\# m\_{H, -1}\ot m\_{0,C, -1}\# m\_{0,H, -1}\ot m\_{0,0}
\\= (m\_{C, -1})^{(1)} \# ((m\_{C, -1})^{(2)})_{(-1)} (m\_{H,
-1})_{(1)} \\ \otimes ((m\_{C, -1})^{(2)})_{(0)} \# (m\_{H,
-1})_{(2)} \ot m\_{0}.
\end{multline}
Applying $p_C\ot p_H \ot \id$ to \eqref{coassociativity}, we get
\begin{align*}
(\id \ot\delta_{H})\delta_{C}(m)&= m\_{C, -1}\cou( m\_{H,
-1})\# \cou(m\_{0,C, -1}) m\_{0,H, -1}\ot m\_{0,0} \\ &=
m\_{C, -1}\# m\_{H, -1} \ot m\_{0} = \delta_{C\# H}(m).
\end{align*}
Applying $(\trasp\ot \id)(p_H\ot p_C \ot \id)$ to
\eqref{coassociativity}, we get
\begin{align*}
(\trasp\ot \id)&(\id \ot\delta_{C})\delta_{H}(m) =   \trasp
\left((m\_{C, -1})_{(-1)} m\_{H, -1} \otimes (m\_{C,
-1})_{(0)}\right) \ot m\_{0} \\ &= (m\_{C, -1})_{(0)} \ot
\Ss^{-1}\left((m\_{C, -1})_{(-1)}\right)(m\_{C, -1})_{(-2)} m\_{H,
-1}  \ot m\_{0}\\ &= \delta_{C\#H}(m).
\end{align*}
Now (ii) follows from the first equality in
Lemma~\ref{lema:smashcopr} (i). Finally, the equality of the first
and third expressions in Lemma~\eqref{lema:smashcopr} (i) gives
that the $C\#H$-subcomodule generated by $Z$ is contained in (and
hence it coincides with) the $C$-subcomodule generated by $Z$.
This gives (iii). \epf

\subsection{Braided Hopf algebras and
bosonization}\label{subsection:bosonization}
We briefly summarize
results from \cite{Ra}, see also \cite{Mj}.
Let $A$ be a Hopf
algebra provided with Hopf algebra maps $\pi: A\to H$, $\imath:
H\to A$, such that $\pi\iota = \id_H$. In other words, we have a
commutative diagram in the category of Hopf algebras:
$$
\xymatrix{& &
H\ar@{=}[1,0]\ar@{_{(}->}[1,-2]_{\iota}\\A\ar[0,2]^{\pi}& & H. }
$$
Let  $R = A^{\co H} = \{a\in A\,|\,(\id\otimes \pi_{H})\Delta (a) =
a\otimes 1\}$. Then $R$ is a braided Hopf algebra in $\ydh$.
Following the notation in Subsection
\ref{subsection:smashcoproduct}, let $\Delta(r) = r\^{1}\otimes
r\^{2}$ denote the coproduct of $r\in R$ (or any other braided
Hopf algebra). Explicitly, $R$ is a subalgebra of $A$, and
\begin{equation}\label{smash2}
\begin{aligned}h\cdot r &= h_{(1)} r\Ss(h_{(2)}), \\ r\_{-1}\ot r\_{0} &=
\pi(r\_{1})\ot r\_{2}, \\
r\^{1}\ot r\^{2} &= \vartheta_R(r\_{1})\ot r\_{2},
\end{aligned}\end{equation}
$r\in R$, $h\in H$. Here $\vartheta_R: A \to R$ is the map defined
by
\begin{equation}\label{vartheta}\vartheta_R(a) = a\_{1}\iota\pi(\Ss
(a\_{2})),\end{equation} $a\in A$. It can be easily shown that
\begin{equation}\label{proptheta}
\vartheta_R(rh) = r\cou(h),\qquad \vartheta_R(hr) = h\cdot
r\end{equation} for  $r\in R$, $h\in H$. Reciprocally, let $R$ be
a braided Hopf algebra in $\ydh$.  A construction discovered by
Radford, and interpreted in terms of braided categories by Majid,
produces a Hopf algebra $R\# H$ from $R$. We call $R\# H$ the
bosonization of $R$. As a vector space, $R\# H = R\ot H$; if $r\#
h := r\ot h$, $r\in R$, $h\in H$, then the multiplication and
comultiplication of $R\# H$ are given by
\begin{equation}\label{smash1}
\begin{aligned}(r\# h)(s\# f) &= r (h_{(1)} \cdot s)\# h_{(2)}f, \\ \Delta(r\# h) &=
r^{(1)} \# (r^{(2)})_{(-1)} h_{(1)} \otimes (r^{(2)})_{(0)}
 \#  h_{(2)}.
 \end{aligned}
 \end{equation}
The maps $$\pi_{H}: R\# H \to H \text{ and }\imath: H \to R\# H,
\quad \pi_{H}(r\# h) = \cou(r)h, \quad \imath(h) = 1\# h,$$
$r\in R$, $h\in H$, are Hopf algebra homomorphisms; we identify
$H$ with the image of $\imath$. Hence
\begin{equation}\label{smash3}
r\_{1}\ot r\_{2} = r\^{1}(r\^{2})\_{-1}\ot (r\^{2})\_{0},
\end{equation}
$r\in R$. The map $p_{R}: R\# H \to R$, $p_{R}(r\# h) =
r\cou(h)$, $r\in R$, $h\in H$, is a coalgebra
homomorphism -- see page \pageref{subsection:smashcoproduct}. We
shall write $rh$ instead of $r\# h$, $r\in R$, $h\in H$. The
antipodes $\Ss_R$ of $R$ and $\Ss = \Ss_{R\# H}$ of $R\# H$ are
related by
\begin{align}\label{antipodas}
\begin{aligned}\Ss_R(r) &= r\_{-1} \Ss(r\_{0}), \\
 \Ss(r) &= \Ss(r\_{-1}) \Ss_R(r\_{0}),
\end{aligned}\end{align} $r\in R$.
The antipode $\Ss_R$ is a morphism of Yetter-Drinfeld modules. Let
$\mu$ be the multiplication of $R$ and $c\in\End(R\otimes R)$ be
the braiding.
Then $\Ss _R$
is anti-multiplicative and anti-comultiplicative
in the following sense:
\begin{align}\label{eqn:antipodatrenzada-anti}
\begin{aligned}\Ss_R\mu &= \mu (\Ss_R\ot \Ss_R)c = \mu c(\Ss_R\ot \Ss_R), \\
\Delta\Ss_R &=  (\Ss_R\ot \Ss_R)c\Delta = c(\Ss_R\ot \Ss_R)\Delta ,
\end{aligned}\end{align}
see for instance \cite[1.2.2]{AG}. The adjoint representation of
$R$ on itself is the algebra map $\ad_c: R\to \End R$, $\ad_c x(y)
= \mu(\mu\otimes\Ss)(\id\otimes c)(\Delta\otimes\id)(x\otimes y)$,
$x,y\in R$. That is,
\begin{equation}\label{braided-adj}
\ad_c x(y) = x\^{1}[(x\^2)\_{-1}\cdot y]\Ss((x\^2)\_{0}) = \ad x(y)
\end{equation}
for all $x,y\in R$, where the second equality follows immediately
from \eqref{smash2} and  \eqref{antipodas}. If $x\in {\Pc} (R)$,
then
\begin{align}
  \ad_c x(y) =xy - (x\_{-1}\cdot y)x\_{0}
  \label{eq:adcxy}
\end{align}
for all $y\in R$.
Similarly, define
$$\ad _{c^{-1}} x(y)=xy-y\_0(\Ss^{-1}(y\_{-1})\cdot x) $$
for $x\in {\Pc}(R)$, $y\in R$. We record the next well-known
remark for further reference.
\begin{obs}\label{prim-ydh}
The space of primitive elements $\Pc(R)$ is a Yetter-Drinfeld
submodule of $R$. \qed
\end{obs}

The next consequences of \eqref{eqn:antipodatrenzada-anti} will be
used later.

\begin{lema}\label{S}
  \begin{enumerate}
    \item[(i)]
Let $x\in \Pc(R)$, $y\in R$. Then
\begin{equation}
\label{S1}\ad_c x(\Ss_{R}(y)) = \Ss_{R}(\ad_{c^{-1}} x(y)).
\end{equation}
\item[(ii)] Let  $X$ be a Yetter-Drinfeld submodule of $R$ and let $K$
be the subalgebra generated by $X$. Then $\Ss_{R}(K)$ is the
subalgebra generated by $\Ss_{R}(X)$.
  \end{enumerate}
\end{lema}

\pf Since $\Ss_{R}(x) = -x$, (i) follows directly from
\eqref{eqn:antipodatrenzada-anti}: $\ad_c x(\Ss_{R}(y)) =
-\mu(\Ss_{R}\ot \Ss_{R})(\id - c) (x\ot y)
\overset{\eqref{eqn:antipodatrenzada-anti}}= -\Ss_{R}\mu(c^{-1} -
\id) (x\ot y) = \Ss_{R}(\ad_{c^{-1}} x(y))$.

(ii). If $X$, $Y$ are Yetter-Drinfeld submodules of $R$, then $XY$
is also a Yetter-Drinfeld submodule and $\Ss_R(XY) =
\Ss_R(Y)\Ss_R(X)$ by \eqref{eqn:antipodatrenzada-anti}. This
implies immediately (ii). \epf

\begin{obs}\label{smash-bos}
Let $K$ be a left $A$-module algebra, that is,
$K$ is a left $H$-module algebra and a left $R$-module such that
the action $\cdot $ of $R$ on $K$ satisfies equation
$r\cdot (k\widetilde{k}) = \big(r\^{1}\cdot ((r\^{2})\_{-1}\cdot k)\big)
\big((r\^{2})\_{0}\cdot \widetilde{k}\big)$
for all $r\in R$, $k, \widetilde k \in K$.

  (i) The smash product $K\# A$ is a right $H$-comodule
  algebra via the coaction $(\id \# \id \ot \pi)(\id \# \Delta)$,
  with subalgebra of coinvariants $K\# R$.
  According to \eqref{eq:smashcopr}, the product in the last is given by
$$(k\#r)(k'\#r') = k(r\^{1}(r\^{2})\_{-1})\cdot k' \#
(r\^{2})\_{0}r',$$ $k,k'\in K$, $r,r'\in R$.

(ii) The multiplication induces a linear isomorphism $R\ot K \to
K\# R$. The inverse map is given by $k\#r \mapsto r\_{2}\ot
\Ss ^{-1}(r\_{1})\cdot k$.
\end{obs}

\begin{obs}\label{obs:cop} Let $B$ be a braided bialgebra.
Let $B^{\cop}$ denote the algebra $B$ together with the
comultiplication $c^{-1}\Delta$; this is a braided Hopf algebra
but with the inverse braiding, see \cite[Prop. 2.2.4]{AG}.
Clearly, $\Pc(B) = \Pc(B^{\cop})$.
\end{obs}

\subsection{Nichols algebras}

Let $V\in \ydh$.   The tensor algebra $T(V)$ is a braided Hopf
algebra in $\ydh$. A very important example of braided Hopf
algebra in $\ydh$ is the Nichols algebra  $\toba(V)$ of $V$; this
is the quotient of $T(V)$ by a homogeneous ideal $\J = \J(V)$,
generated by (some) homogeneous elements of degree $\ge 2$. See
\cite{AS-cambr} for the precise definition and main properties of
Nichols algebras, and the relation with pointed Hopf algebras.

\medbreak Another description of the ideal $\J(V)$ is as the
kernel of the quantum symmetrizer introduced by Woronowicz
\cite{Wo}, see \cite{Sbg}. Let $\mathbb B_{n}$ be the braid group
in $n$ letters and let $\pi: \mathbb B_{n} \to \mathbb S_{n}$ be a
natural projection; it admits a set-theoretical section $s:\mathbb
S_{n} \to \mathbb B_{n}$ called the Matsumoto section. Let
$\mathfrak S_{n} := \sum_{\sigma \in \mathbb S_{n}} s(\sigma)$.
The braid group $\mathbb B_{n}$ acts on $T^n(V)$ via $c$ and the
homogeneous component $\J^n(V)$ of $\J (V)$ equals $\ker \mathfrak
S_{n}$. Thus $\toba(V)$ depends (as algebra and coalgebra) only on
the braiding $c$. We write $\toba(V) = \toba(V,c)$, $\J(V) =
\J(V,c)$.

The Nichols algebra has a unique grading $\toba(V) = \oplus_{n\in
\N_0}\toba^n(V)$ such that $\toba^1(V)=V$, the multiplication and the
comultiplication are graded, and the action and the coaction of $H$ are
homogeneous.

If $\dim V<\infty $, then
there exists a bilinear form $\langle\, , \,\rangle: T({V}^*)
\times T(V) \to \ku$ such that
\begin{align}\label{eqn:duality-homogeneous-tensor}
\langle T^n({V}^*) , T^m(V)\rangle &= 0, \qquad n\neq m,
\\ \label{eqn:duality-tensor} \langle f_n \dots  f_1, x\rangle &= \langle
f_n\ot \dots \ot f_1, \Delta_{1, \dots, 1}(x)\rangle
\end{align}
for $f_1, \dots, f_n\in {V}^*$, $x\in T^n(V)$, $n\in \N_0$. It
satisfies the following properties:

\begin{align}\label{properties:duality1-tensor}
\langle fg, x\rangle &= \langle f, x\^{2}\rangle\langle g,
x\^{1}\rangle,
\\\label{properties:duality2-tensor}
\langle f, xy\rangle &= \langle f\^{2}, x\rangle\langle f\^{1},
y\rangle,
\\\label{properties:duality3-tensor}
\langle h\cdot f, x\rangle &= \langle f, \Ss(h)\cdot x\rangle,
\\\label{properties:duality4-tensor} f\_{-1}\langle f\_{0}, x\rangle &= \Ss^{-1}(x\_{-1})\langle f, x\_{0}\rangle
\end{align}
for all $f, g\in T({V}^*)$, $x, y\in T(V)$, $h\in H$. This was
first observed in \cite{Mj93}, see also \cite[10.4.13]{Mj-book}. A
combination of the explicit formulas in \cite[10.4.13]{Mj-book}
and \cite[Eqs.\,(3.25), (3.26)]{Wo} shows that
$$\Delta_{1,\ldots ,1}= \mathfrak S_{n}$$
for all $n\in \N $,
that is, $\J (V,c)$ is the radical of the form in the second argument.
More precisely, the following holds.

\begin{prop}\label{prop:duality}
  \cite[Thm.\,3.2.29]{AG}
Assume that $V\in \ydh $ such that $\dim V<\infty $.
  Then there exists a non-degenerate bilinear form $$\langle\, ,
\,\rangle: \toba({V}^*) \times \toba(V) \to \ku$$ such that
\begin{align}\label{eqn:duality-homogeneous}
\langle\toba^n({V}^*) , \toba^m(V)\rangle &= 0, \qquad n\neq m,
\\ \label{eqn:duality} \langle f_n \dots  f_1, x\rangle &= \langle
f_n\ot \dots \ot f_1, \Delta_{1, \dots, 1}(x)\rangle,
\end{align}
for $f_1, \dots, f_n\in {V}^*$, $x\in \toba^n(V)$, $n\in \N_0$. It
satisfies \eqref{properties:duality1-tensor},
\eqref{properties:duality2-tensor},
\eqref{properties:duality3-tensor}, and
\eqref{properties:duality4-tensor} for all $f, g\in \toba({V}^*)$,
$x, y\in \toba(V)$, $h\in H$.
\end{prop}

This proposition tells that
\begin{align}
  \toba(V)^{\grd}\simeq\toba(V^*),
  \label{eq:BVgrdual}
\end{align}
where $\toba(V)^{\grd}$ is the graded dual of $\toba(V)$, see
\eqref{eq:grdual}.

\begin{lema}\label{lema:toba-c-menos-uno} $\J(V, c) = \J(V,c^{-1})$
and $\toba(V,c) \simeq \toba(V,c^{-1})$ as algebras.
\end{lema}
\pf Let $\toba(V)^{\cop}$ be the opposite coalgebra, see Remark
\ref{obs:cop}. Clearly, the algebra $\toba(V)^{\cop}$ is generated
in degree one, and $\Pc\left(\toba(V)^{\cop}\right) =
\Pc\left(\toba(V)\right) = V$. Hence $\toba(V)^{\cop} \simeq
\toba(V,c^{-1})$, and $\J(V, c) = \J(V,c^{-1})$. \epf

\begin{lema}\label{rmk:gradedcondition}
Let $x = \sum_{n\ge 1} x(n)\in \toba(V)$, with $x(n)\in
\toba^n(V)$. Assume that $x\^{1} \ot \pro_1(x\^{2}) = 0$. Then
$x=0$.
\end{lema}

\pf  From  $0 =x\^{1} \ot \pro_1(x\^{2}) =\sum_{n\ge 1} x(n)\^{1}
\ot \pro_1(x(n)\^{2})$ we conclude that $\Delta_{n-1, 1}(x(n)) =
x(n)\^{1} \ot \pro_1(x(n)\^{2}) = 0$, since $x(n)\^{1} \ot
\pro_1(x(n)\^{2})$ $\in \toba^{n-1}(V)\ot \toba^{1}(V)$. But
$\Delta_{n-1, 1}$ is injective in a Nichols algebra, hence $x(n) =
0$ for all $n$ and \emph{a fortiori} $x=0$. \epf

For simplicity, we write $\ac(V) = \toba(V)\# H$ for the
bosonization of $\toba(V)$. Then $\ac(V) = \oplus_{n\in
\N_0}\ac^n(V)$, where $\ac^n(V) = \toba^n(V)\# H$, is a graded
Hopf algebra.

\section{The algebra of quantum differential
operators}\label{section:qdo}

We now discuss two algebras of quantum differential operators that
appeared frequently in the literature. For quantum groups, it
seems that they were first defined in \cite{Ka}, see also
\cite[Chapter 15]{L}. For Yetter-Drinfeld modules over finite
group algebras, see \cite{G-jalg}.

\subsection{The algebra of quantum differential operators}

Let $B$ be a brai\-ded bialgebra in $\ydh$. Then the space of
linear endomorphisms $\End B$ is an associative algebra with
respect to the convolution product: $T*S\,(b) = T(b\^2)S(b\^1)$,
$T, S\in \End B$, $b\in B$, a convention coherent with
\eqref{eqn:vndual}. Since $B$ is a left and right comodule over
itself via the comultiplication, it becomes a left and right
module over $B^*$. If $\xi\in B^*$, then we define the quantum
differential operators $\pl{}, \pr{}: B^* \to \End B$ as the
representations associated to those actions. That is,
\begin{equation}\label{eqn:defiqdo-gral}
\pl{\xi}(b) = \langle \xi, b\^1\rangle b\^2, \qquad \pr{\xi}(b)  =
\langle \xi, b\^2\rangle b\^1, \qquad b\in B,\, \xi\in B^*.
\end{equation}
Let also $L, R:B\to \End B$ be the left and right regular
representations.

\medbreak If $\xi,\zeta\in B^*$, then clearly $\pl{\zeta}\pr{\xi}
= \pr{\xi}\pl{\zeta}$. Other basic properties of the quantum
differential operators are stated in the next lemma.

Recall that $A^{\sw}$ denotes the Sweedler dual of an algebra $A$.
Explicitly,
\begin{align*}A^{\sw} &= \{f\in \Hom(A, \ku)\,|\,\ker f \text{
contains a left ideal $I$ of finite codimension}\}.
\end{align*}

\begin{lema}\label{lema:qdo-gral}

{\rm (i)} The maps $\pl{}: B^* \to \End B$ and $\pr{}: B^{*\op}
\to \End B$ are injective algebra homomorphisms.

{\rm (ii)} If $B$ is a braided Hopf algebra with bijective
antipode,
    then the maps $\Psi^L, \Psi^R: B\ot B^* \to \End B$,
    $\Psi^L(b\ot \xi) = L_b\circ \pl{\xi}$,
    $\Psi^R(b\ot \xi) = R_b\circ \pr{\xi}$,
    are injective.

{\rm (iii)} If $\xi\in B^{\sw}$ and $b,c\in B$, then
    \begin{align}\label{eqn:leibniz-pl-gral}
\pl{\xi} (bc) &= \langle \xi\^2, b\^1 \rangle (b\^2)\_{0} \,
\pl{\Ss^{-1}((b\^2)\_{-1})\cdot\xi\^1}  (c),
\\\label{eqn:leibniz-pr-gral} \pr{\xi} (bc) &=
\pr{(\xi\^2)\_{0}}(b)\, \Ss\left((\xi\^2)\_{-1}\right) \cdot \pr{
\xi\^1}(c).
\end{align}

{\rm (iv)} If $\xi\in \Pc(B^{\sw})$ and $b,c\in B$, then
    \begin{align}\label{eqn:leibniz-pl-prim}
\pl{\xi} (bc) &= b\_{0} \pl{\Ss^{-1}(b\_{-1})\cdot\xi} (c) +
\pl{\xi}(b) c,
\\\label{eqn:leibniz-pr-prim} \pr{\xi} (bc) &= b \pr{\xi}(c)
+ \pr{\xi\_{0}}(b)\, \Ss\left(\xi\_{-1}\right) \cdot c.
\end{align}

{\rm (v)} Let $U$ be a Yetter-Drinfeld submodule of
      $\Pc(B^{\sw})$. Let $S$ be the subalgebra of $B^{\sw}$
    generated by $U$.
    Then $\D^L(B, U) :=
    L(B)\circ \pl{}(S)$ and $\D^R(B, U) := R(B)\circ \pr{}(S)$ are
    subalgebras of $\End B$.
\end{lema}

\pf (i). If $\pl{\xi}(b) = 0$, then $\langle \xi, b\rangle =
\cou\pl{\xi}(b) = 0$; thus $\pl{}$ is injective -- and
similarly for $\pr{}$. Now, if $b\in B$, $\xi, \zeta\in B^*$, then
\begin{align*}
\pl{\zeta}\pl{\xi}(b) &= \langle \xi, b\^1\rangle \pl{\zeta}
(b\^2) = \langle \xi, b\^1\rangle \langle \zeta, b\^2\rangle b\^3
= \langle \zeta *\xi, b\^1\rangle  b\^2 = \pl{\zeta *\xi}(b),
\\ \pr{\zeta}\pr{\xi}(b) &= \langle \xi, b\^2\rangle
\pr{\zeta} (b\^1) = \langle \xi, b\^3\rangle \langle \zeta, b\^2\rangle b\^1
= \langle \xi *\zeta,
b\^2\rangle  b\^1 = \pr{\xi *\zeta}(b).
\end{align*}

(ii). Let $\sum_i b_i\ot\xi_i\in \ker \Psi^L$, and
assume that the $b_i$'s are linearly independent. Thus $\sum_i b_i
\langle \xi_i, b\^1\rangle b\^2 = 0$ for any $b\in B$. Therefore
\begin{align*}
\sum_i b_i \langle \xi_i, b\rangle  = \sum_i b_i \langle \xi_i,
b\^1\rangle b\^2 \Ss_B(b\^3)= 0 \implies \langle \xi_i, b\rangle  = 0
\end{align*}
for all $i$ and $b\in B$; hence $\xi_i = 0$ for all $i$. The
argument for $\Psi^R$ is similar.

(iii). We compute
{\allowdisplaybreaks
\begin{align*}
  \pl{\xi} (bc) &\overset{\phantom{\eqref{vd1}}}=
  \langle \xi, (bc)\^1\rangle (bc)\^2
  = \langle \xi\^1\ot \xi\^2, b\^1\ot (b\^2)\_{-1}\cdot c\^1\rangle
  (b\^2)\_{0} c\^2 \\
  &\overset{\phantom{\eqref{vd1}}}= \langle \xi\^2, b\^1 \rangle (b\^2)\_{0}\,
  \langle \xi\^1,(b\^2)\_{-1}\cdot c\^1\rangle  c\^2 \\
  &\overset{\eqref{vd1}}= \langle \xi\^2, b\^1 \rangle (b\^2)\_{0} \,
  \langle \Ss^{-1}\big((b\^2)\_{-1}\big)\cdot\xi\^1, c\^1\rangle  c\^2; \\
  \pr{\xi} (bc) &\overset{\phantom{\eqref{vd1}}}=
  \langle \xi, (bc)\^2\rangle (bc)\^1
  = \langle \xi\^1\ot \xi\^2, (b\^2)\_{0} \ot c\^2\rangle
  b\^1 (b\^2)\_{-1}\cdot c\^1 \\
  &\overset{\phantom{\eqref{vd1}}}= \langle \xi\^2, (b\^2)\_{0} \rangle
  b\^1 \,\langle \xi\^1,c\^2\rangle (b\^2)\_{-1}\cdot c\^1 \\
  &\overset{\eqref{vd2bis}}= \langle (\xi\^2)\_{0}, b\^2 \rangle
  b\^1 \,\langle \xi\^1,c\^2\rangle \Ss\left((\xi\^2)\_{-1}\right) \cdot c\^1.
\end{align*}
}

Now (iv) follows at once from (iii). Next,
\eqref{eqn:leibniz-pl-prim} and \eqref{eqn:leibniz-pr-prim} say
that
\begin{align}\label{eqn:leibniz-pl-appl} \pl{\xi} \circ L_b &=
L_{b\_{0}} \circ\pl{\Ss^{-1}(b\_{-1})\cdot\xi}  + L_{\pl{\xi}(b)},
\\\label{eqn:leibniz-pr-appl} \pr{\xi}\circ R_c &=  R_{\pr{\xi}(c)}
+ R_{\Ss\left(\xi\_{-1}\right) \cdot c} \circ\pr{\xi\_{0}}
\end{align}
for $\xi\in \Pc(B^{\sw})$, $b,c\in B$. These equalities imply (v).
\epf

\begin{exa}\label{exa:qdo} (i). If $B$ is a usual bialgebra, then
the generalized Leibniz rules \eqref{eqn:leibniz-pl-gral} and
\eqref{eqn:leibniz-pr-gral} simply say that $$\pl{\xi} (bc) =
\pl{\xi\^2} (b) \pl{\xi\^1}(c), \quad \pr{\xi} (bc) = \pr{\xi\^2}
(b) \pr{\xi\^1}(c).$$

\medbreak (ii). Let $W\in \ydh$ be \fd{} and let $B =
\toba(W)$. By Prop.~\ref{prop:duality}, there exists an
embedding $\toba(W^*) \to \toba(W)^*$ and we can consider the
algebras of \emph{quantum differential operators} \begin{align*}
\D^L(W) &:= \D^L(\toba(W), W^*) = L(\toba(W))\circ
\pl{}(\toba(W^*)), \\ \D^R(W) &:= \D^R(\toba(W), W^*) =
R(\toba(W))\circ \pr{}(\toba(W^*));
\end{align*} these are subalgebras of $\End B$ by
Lemma~\ref{lema:qdo-gral}(v).

\medbreak  Let $x\in\toba(W)$ be homogeneous of degree $p$. Let us
write in this case
$$
\Delta(x) = x\ot 1 + 1\ot x + \sum_{0<r<p} x'_r\ot x''_{p-r}.
$$
Here we use a symbolic notation with $x'_r\in\toba^r(W)$,
$x''_{p-r}\in\toba^{p-r}(W)$. If $f\in W^*$ and $p>1$, then
$\pl{f}(x) = \langle f, x'_1\rangle x''_{p-1}$, $\pr{f}(x) =
x'_{p-1}\langle f, x''_1\rangle$. Also, $\pl{f}(w) = \langle f,
w\rangle  = \pr{f} (w)$ for $w\in W$.

The following fact is well-known and goes back essentially to
\cite{N}: If $x\in \toba(W)$ and $\pl{f}(x) = 0$ for all $f\in
W^*$, then $x\in \ku$.

\medbreak (iii). Let $W$ and $\toba(W)$ as in (ii).
Assume that $W$ admits a basis $v_1, \dots,
v_\theta$ such that $\delta(v_i) = g_i \ot v_i$, for some $g_i\in
G(H)$, $1\le i \le \theta$. Let $f_1, \dots, f_\theta$ be the dual
basis; then $\delta(f_i) = g_i^{-1} \ot f_i$, $1\le i \le \theta$.
Set $\partial_i = \pr{f_i}$. Then
$$ \partial_i(bc) = b\partial_i(c) + \partial_i(b) \,g_i\cdot c,
\qquad b,c\in \toba(W).
$$
Similarly, let $\Alg(H,\ku )$ be the group of algebra homomorphisms
from $H$ to $\ku$;
it acts on $B$ by $\chi\cdot b = \langle\chi,
b\_{-1}\rangle b\_0$. Suppose that $W$ admits a basis $v_1, \dots,
v_\theta$ such that $h\cdot v_i = \chi_i(h) v_i$, for some
$\chi_i\in \Alg(H,\ku )$, $1\le i \le \theta$. Let $f_1, \dots,
f_\theta$ be the dual basis; then $h\cdot f_i = \chi_i^{-1}(h)
f_i$, $1\le i \le \theta$. Set $\underline{\partial}_i =
\pl{f_i}$. Then
$$ \underline{\partial}_i(bc) = (\chi_i\cdot b)\underline{\partial}_i(c) + \underline{\partial}_i(b) c,
\qquad b,c\in \toba(W).
$$
\end{exa}

\begin{prop}\label{prop:qdo-toba} Let $W\in \ydh$ be \fd.

(1). The map
$$\Psi ^L:\toba (W)\otimes \toba(W^*)\to \D^L(W)$$ is a linear isomorphism.

(2).  The map $\Theta:T(W \oplus W^*) \to \D^L(W)$, $(v, f)
\mapsto L_v\circ \pl{f}$, $v\in W$, $f\in W^*$, induces an algebra
isomorphism $\vartheta: T(W \oplus W^*)/I \to \D^L(W)$, where $I$
is the two-sided ideal generated by

\begin{enumerate}
    \item[(i)] the relations of $\toba(W)$,
    \item[(ii)]  the relations of $\toba(W^*)$,
    \item[(iii)] the relations
    \begin{align}\label{eqn:rel-qdo} f v =
v\_{0}\, \Ss^{-1}(v\_{-1})\cdot f  + \pl{f}(v), \qquad v\in W, \,
f\in W^*.
\end{align}
\end{enumerate}

\end{prop}

\medbreak  If $v\in W$, $f\in W^*$, then \eqref{eqn:rel-qdo}
implies that
\begin{align}\label{eqn:rel-qdo-bis} vf =
(v\_{-1}\cdot f) \, v\_{0}  - \pl{v\_{-1}\cdot f}(v\_0).
\end{align}

\pf By what was already said, $\Theta$ induces $\vartheta$ and
this is surjective. Indeed, \eqref{eqn:leibniz-pl-appl} says more
generally that
 \begin{align}\label{eqn:rel-qdogral} f x =
x\_{0}\, \Ss^{-1}(x\_{-1})\cdot f  + \pl{f}(x), \qquad x\in
\toba(W), \, f\in W^*.
\end{align}
Clearly, the inclusions of $W$ and $W^*$ induce algebra maps $j_W:
\toba(W) \to T(W \oplus W^*)/I$,
$j_{W^*}: \toba(W^*) \to T(W \oplus W^*)/I$.
Let $\mu$ be the multiplication of
$T(W \oplus W^*)/I$.
Then \eqref{eqn:rel-qdo} guarantees that $\mu
\circ(j_W\ot j_{W^*})$ is surjective. But the following diagram
commutes:
\begin{equation*}
\xymatrix{\toba(W)\ot\toba(W^*) \ar[1,1]_{\Psi^L} \ar[0,2]^-{\mu
\circ(j_W\ot j_{W^*})} & &  T(W \oplus W^*)/I
\ar@{>>}[1,-1]^{\vartheta}
\\ &\D^L(W), &}
\end{equation*}
and $\Psi^L$ is a linear isomorphism by Lemma~\ref{lema:qdo-gral}
(ii). Thus $\mu \circ(j_W\ot j_{W^*})$ and $\vartheta$ are
isomorphisms. \epf

\begin{cor}\label{cor:qdo-ydh} Let $W\in \ydh$ be \fd.
Then $\D^L(W)$ is an algebra in $\ydh$.
\end{cor}

\pf
Straightforward. \epf

\begin{definition}
\label{def:ksmash} Let $\toba(W)\#\toba(W^*)$ denote the vector
space $\toba(W)\otimes \toba(W^*)$ with the multiplication
transported along the isomorphism $\Psi ^L:\toba(W)\otimes
\toba(W^*)\to \D^L(W)$. Thus
\begin{align}\label{eqn:BWBW*} \xi \, b &=
\langle \xi\^2, b\^1 \rangle (b\^2)\_{0} \#
\Ss^{-1}((b\^2)\_{-1})\cdot \xi\^1
\end{align}
for $b\in \toba(W)$ and $\xi\in \toba(W^*)$ by
\eqref{eqn:leibniz-pl-gral}. The multiplication map
  $\mu :\toba(W^*)\otimes \toba(W)\to \toba (W)\#\toba(W^*)$
  is an isomorphism, with inverse map $\mu ^-$ given by
  \begin{align}
  \mu^-:b\xi \mapsto
  \langle \Ss_{\toba(W^*)}( (b\_{-1}\cdot \xi )\^2),b\_0{}\^1\rangle
  (b\_{-1}\cdot \xi )\^1\otimes b\_0{}\^2
    \label{eq:multinv}
  \end{align}
  for $b\in \toba(W)$, $\xi \in \toba(W^*)$. Note that $\Psi ^L:\toba(W)\#\toba(W^*)\to \D^L(W)$ is an isomorphism
  in $\ydh $, where $H$ acts and coacts diagonally on
  $\toba(W)\#\toba(W^*)$, see Corollary~\ref{cor:qdo-ydh}.
\end{definition}

\begin{obs}\label{rem:Heisenberg}
  Alternatively to the above construction, the algebra $\D^L(W)$ can be
  obtained as the subalgebra $\toba (W)\#\toba(W^*)$
  of the Heisenberg double $\ac (W)\# \ac (W)^{\circ }$.
  Here for any Hopf algebra $A$ the Heisenberg double $A\# A^{\circ }$
  is the smash product algebra corresponding to the left action of
  the Hopf dual $A^\circ $, see our convention in Remark \ref{rem:C*alg},
  given by the left $A$-coaction on $A$ via $\Delta $. The embedding of
  $\toba(W)\otimes \toba(W^*)$ into $\ac (W)\# \ac (W)^{\circ }$ is given by
  the inclusion of $\toba(W)$ and the map
  $$\toba(W^*)\ni f\mapsto \langle f, \cdot \rangle \ot \cou
  \in (\toba (W)\# H)^{\circ }=\ac (W)^{\circ }.$$
  One can check that $\toba(W)\otimes \toba(W^*)\subset
  \ac(W)\# \ac(W)^{\circ }$
  is a subalgebra and that this algebra structure on
  $\toba (W) \otimes \toba(W^*)$ coincides with $\toba (W) \# \toba(W^*)$
  as in Definition \ref{def:ksmash}. Further, the restriction of the map in
  Remark \ref{smash-bos} (ii) coincides with the map in \eqref{eq:multinv}.
  These facts will not be used in the sequel.
\end{obs}

\begin{obs}\label{re:multinv}
Let $K$ be a subalgebra of $\toba(W)$ and $\Kb$ be a
  braided Hopf subalgebra of $\toba(W^*)$ such that
\begin{itemize}
    \item $K$ is an $H$-subcomodule,
    \item $\Kb$ is an $H$-submodule,
    \item $\pl{\xi}(b) =\langle \xi ,b\^1 \rangle b\^2\in K$
      for all $b\in K$, $\xi\in \Kb$.
\end{itemize}
Then $K\ot \Kb $ is a subalgebra of $\toba(W)\#\toba(W^*)$,
denoted by $K\# \Kb$. Again, the multiplication map
  $\mu :\Kb \otimes K\to K\# \Kb $
  is an isomorphism. If $K\subset \toba(W)$ and $\Kb \subset \toba(W^*)$ are subobjects in $\ydh
  $ then $K\#\Kb $ is a subalgebra of $\toba(W)\#\toba(W^*)$ in $\ydh $.
\end{obs}


\begin{obs}\label{obs:Dlgraded} Let $\Gamma$ be an abelian group. Assume that
  $W = \oplus_{\gamma\in\Gamma} V_{\gamma}$ is a \fd{}
$\Gamma$-graded Yetter-Drinfeld module; $W^* \simeq
\oplus_{\gamma\in\Gamma} V^*_{\gamma}$ becomes a $\Gamma$-graded
Yetter-Drinfeld module with $\deg V^*_{\gamma} = - \gamma$. Then
$\toba(W)$,
$\toba(W)\#\toba(W^*)$,
and $\D^L(W)$ are $\Gamma $-graded
algebras.
\end{obs}

\pf The tensor algebras $T(W)$ and $T(W \oplus W^*)$ inherit
the $\Gamma$-grading of Yetter-Drinfeld modules
in the usual way: $\deg (V_{\gamma_1} \ot \dots V_{\gamma_s})
=\gamma_1 + \dots + \gamma_s$. By definition, the braiding $c$
preserves homogeneous components; thus $\toba(W)$ inherits the
grading. Now the relations \eqref{eqn:rel-qdo} are also
homogeneous, hence
$\toba(W)\#\toba(W^*)$ and $\D^L(W)$ are $\Gamma$-graded
algebras. \epf

\subsection{Braided derivations}\label{subsection:skew-derivations}
We next give a characterization of Nichols algebras in terms of
quantum differential operators suitable for our later purposes. Recall
that the kernel of the counit of a bialgebra $B$ is denoted by
$B^+$.

\medbreak   First, let $B$ be a braided bialgebra and consider
$B^{\cop}$ as in Remark \ref{obs:cop}. We write $\Delta(x) =
x^{[1]}\ot x^{[2]}$ to distinguish from the previous coproduct.
Thus $\Delta(xy) = x^{[1]}{y^{[1]}}\_0\ot
\big(\Ss^{-1}({y^{[1]}}\_{-1}) \cdot x^{[2]}\big)y^{[2]}$, for
$x,y\in B$, \emph{cf.} \eqref{inversebraiding}. Let $\xi\in B^*$
and let $\plb{\xi}\in \End B$ be $\pl{\xi}$ for this bialgebra,
that is
\begin{align}
  \plb{\xi}(x) = \langle \xi, x^{[1]}\rangle x^{[2]}.
  \label{eq:plbxi}
\end{align}
Then
\begin{align}\label{plinverso}
\plb{\xi} (xy) &= \big({\xi^{[1]}}\_{-1} \cdot \plb{\xi^{[2]}}(x)\big)
\plb{{\xi^{[1]}}\_0}(y).
\end{align}
Indeed,
\begin{align*}
\plb{\xi} (xy) &= \big\langle \xi, x^{[1]}{y^{[1]}}\_0\big\rangle
\,\big( \Ss^{-1}({y^{[1]}}\_{-1}) \cdot x^{[2]}\big)y^{[2]} \\
&= \big\langle \xi^{[2]}, x^{[1]} \big\rangle\,
\Ss^{-1}({y^{[1]}}\_{-1}) \cdot x^{[2]} \,\big\langle \xi^{[1]},
{y^{[1]}}\_0\big\rangle  y^{[2]}
\\ &\overset{\makebox[0pt]{\tiny \eqref{vd2}}}= \big\langle
\xi^{[2]}, x^{[1]} \big\rangle\, {\xi^{[1]}}\_{-1} \cdot x^{[2]}
\,\big\langle {\xi^{[1]}}\_0,y^{[1]}\big\rangle  y^{[2]} \\
&= \big({\xi^{[1]}}\_{-1} \cdot \plb{\xi^{[2]}}(x)\big)
\plb{{\xi^{[1]}}\_0}(y).
\end{align*}
If $\xi\in\Pc(B) = \Pc(B^{\cop})$, then
\begin{align}\label{plinverso-prim} \plb{\xi} (xy) &= (\xi\_{-1}
\cdot x )\,\plb{\xi\_0}(y) + \plb{\xi}(x) \,y.
\end{align}

\medbreak Part (i) of the following theorem is well-known, but
part (ii) seems to be new.

\begin{theorem}\label{theo:quotients-with-derivations}
  Let $W\in \ydh$ be \fd.
Let $I\subset T(W)^+$ be a 2-sided ideal, stable under the action
of $H$. Let $R = T(W)/I$ and let $\pi: T(W) \to R$ be the
canonical projection.
\begin{enumerate}
  \item[(i)] Assume that $I$ is a homogeneous Hopf ideal, so that $R$ is a
graded braided Hopf algebra quotient of $T(W)$, and that $I\cap W
= 0$. Then for any $f\in W^*$ there exists a map $d_f\in \End R$
such that for all $x,y\in R$, $v\in W$,
\begin{align}
\label{derivadauno-bis} d_f(xy) &= (f\_{-1} \cdot x)\,
d_{f\_{0}}(y)
+ d_{f}(x)\, y,\\
\label{derivadados} d_f(\pi(v)) &= \langle f, v\rangle.
\end{align}
  \item[(ii)] Conversely, assume that for any $f\in W^*$ there exists a
map $d_f\in \End R$ such that
\eqref{derivadauno-bis} and \eqref{derivadados} hold.
Then $I\subseteq \J(W)$, that is,
there exists a unique surjective algebra map
$\Omega: R \to \toba(W)$ such that
$\Omega (\pi (w))=w$ for all $w\in W$.
Moreover
\begin{align}
\label{factorizacion-omega} \Omega d_f &= \plb{f}\Omega.
\end{align}
\end{enumerate}
\end{theorem}

\pf (i). We have $R = \oplus_{n\ge 0} R^n$ with $R^1 \simeq W$; we
identify $W^*$ with a subspace of $R^*$, see Subsection
\ref{subsection:notation}. Hence, if $f\in W^*$, then $d_f :=
\plb{f}$ satisfies \eqref{derivadauno-bis} by
\eqref{plinverso-prim}.

(ii). We apply (i) to $I = 0$; set $D_f\in \End T(W)$,
$D_f = \plb{f}$ for $f\in W^*$. Note that \eqref{derivadados}
implies that $\pi$ restricted to $W$ is injective. We claim that
$d_{f}\pi = \pi D_{f}$, that is, the following diagram commutes:
\begin{equation}\label{factorizacion-d}
\xymatrix{T(W)  \ar[d]_{\pi} \ar[0,2]^-{D_{f}}& &T(W) \ar[d]^{\pi}
\\ R  \ar[0,2]^-{d_f}&& R. }
\end{equation}
For, let $\delta_f = d_{f}\pi$, $\widetilde{\delta}_f = \pi D_{f}:
T(W) \to R$, and let $x,y\in T(W)$. Then
\begin{align*}
d_{f}(\pi(xy)) &= (f\_{-1} \cdot \pi(x))\, d_{f\_{0}}(\pi(y)) +
d_{f}(\pi(x))\, \pi(y);
\\ \pi D_{f}(xy) &= \pi\big((f\_{-1} \cdot x)\, D_{f\_{0}}(y)
+ D_{f}(x)\, y\big)\\
&= (f\_{-1} \cdot \pi(x))\, \pi D_{f\_{0}}(y) + \pi D_{f}(x)\, \pi
(y),
\end{align*}
by the hypothesis on $I$. Also $d_f(\pi(v)) = \langle f, v\rangle
= \pi D_{f}(v)$ for $v\in W$. Thus the set of all $x\in T(W)$ such
that $\delta_f(x) = \widetilde{\delta}_f (x)$ is a subalgebra that
contains $W$; hence $d_{f}\pi = \pi D_{f}$. (This shows that such
a map $d_f$ is unique when it exists; hence $f\mapsto d_f$ is
linear in $f$). In other words, $ D_{f}(I) \subset I$. Let
$\langle\, , \,\rangle: T({W}^*) \times T(W) \to \ku$ be the
bilinear form defined by \eqref{eqn:duality-homogeneous} and
\eqref{eqn:duality}, but with respect to $c^{-1}$. We know that
$\J(W, c^{-1})$ is the (right) radical of this form, and $\J(W,
c^{-1})= \J(W, c)$ by Lemma~\ref{lema:toba-c-menos-uno}; so we
need to show that $\langle T(W^*) , I\rangle = 0$, or equivalently
that $\langle T^n(W^*) , I\rangle = 0$ for all $n\ge 0$. If $n=0$,
then this is clear as $\cou(I) = 0$. If $n=1$, $f\in W^*$ and
$x\in I$, then $\langle f, x\rangle = \cou(\langle f,
x^{[1]}\rangle x^{[2]}) = \cou( D_{f}(x)) \in \cou(I) = 0$. If $n
>1$, $g\in T^{n-1}(W^*)$, $f\in W^*$ and $x\in I$, then
$$\langle gf, x\rangle =
\langle f, x^{[1]}\rangle\langle g, x^{[2]}\rangle = \langle g,
 D_{f}(x)\rangle \in \langle g, D_{f}(I)\rangle\subset \langle
g,I\rangle = 0.$$
In the following diagram, the big and upper squares commute by (i)
and \eqref{factorizacion-d}, respectively:
\begin{equation*}
\xymatrix{T(W)  \ar@/_2pc/[2,0]_{p} \ar[d]_{\pi}
\ar[0,2]^-{D_{f}}& &T(W) \ar[d]^{\pi} \ar@/^2pc/[2,0]^{p}
\\ R  \ar[d]_{\Omega} \ar[0,2]^-{d_{f}}& &R \ar[d]^{\Omega}
\\ \toba(W)  \ar[0,2]^-{\plb{f}}&& \toba(W). }
\end{equation*}
Hence $\plb{f}\Omega\pi = \plb{f}p = pD_f = \Omega\pi D_f = \Omega
d_f \pi$, and since $\pi$ is surjective, $\plb{f}\Omega = \Omega
d_f $. \epf

There are other versions of this theorem. Taking
\eqref{eqn:leibniz-pl-prim} or \eqref{eqn:leibniz-pr-prim} into
consideration, we have similar results replacing the requirement
\eqref{derivadauno-bis} by either of the following:
\begin{align}
\label{derivadauno-izq} d_f(xy)
&= x\_{0} d_{\Ss^{-1}(x\_{-1})\cdot f}(y) + d_{f}(x) y,\\
\label{derivadauno} d_f(xy) &= x d_f(y) +
d_{f\_0}(x)\Ss(f\_{-1})\cdot y,
\end{align}
where $x,y\in R$, $f\in W^*$. The proof goes exactly as for
Theorem \ref{theo:quotients-with-derivations}.

\medbreak  The results in Theorem
\ref{theo:quotients-with-derivations} motivate the following
definition.

\begin{definition}\label{defi:skew-der} Let $M\in {}^H\m$,
$R$ an algebra, $T$ an $H$-module algebra, $\wp: R\to T$ an
algebra map, and let $d: M\to \Hom(R, T)$ be a linear map, denoted
by $f\mapsto d_f$. Following \cite{Mj93} we say that $d$ is a
\emph{family of braided derivations}
if for all $x,y\in R$, $f\in M$,
\begin{equation}
\label{derivadagral} d_f(xy) = (f\_{-1} \cdot \wp(x))\,
d_{f\_{0}}(y) + d_{f}(x)\, \wp(y). \end{equation}
\end{definition}

We are mostly concerned with the case when
$R = T$ and $\wp = \id$.
In this case we say that
$d$ is a \emph{family of braided derivations of $R$.}

\begin{definition}\label{de:canfam}
  Let $W\in \ydh $.
  The family $d^W:W^*\to \End \toba(W)$ of braided derivations
  of $\toba (W)$
  with $d^W_f(w)=\langle f,w\rangle $ for all $f\in W^*$ and
  $w\in W$, see Theorem \ref{theo:quotients-with-derivations} (i),
  is called the \textit{canonical family of braided derivations}
  of $\toba (W)$.
\end{definition}

Our next goal is to develop basic properties of families of
braided derivations which will be useful in the sequel.

\begin{lema}\label{le:fambrderTV}
  Let $M\in {}^H\m$, $V$ a vector space, $T$ an $H$-module algebra,
  and $\wp : V\to T$ a linear map. Then any family of braided derivations
  $d:M\to \Hom (T(V),T)$ determines a linear map $d^1:M\to
  \Hom (V,T)$ by letting $d^1_f=d_f|_V$, $f\in M$.
  Conversely, any linear map $d^1:M\to \Hom
  (V,T)$ gives rise to a unique family of braided derivations $d:M\to \Hom
  (T(V),T)$, where $d_f|_V=d^1_f$, $f\in M$.
\end{lema}

\pf If $d$ is a family of braided derivations, then linearity of
$d$ gives that $d^1:M\to \Hom (V,T)$ is a linear map. On the other
hand, if $V$, $\wp :V\to T$, and $d^1:M\to \Hom (V,T)$ are given,
then $\wp $ extends uniquely to an algebra map $\wp :T(V)\to T$,
and the formula
\begin{align*}
  d_f(v_1v_2\cdots v_n)=\sum _{i=1}^n &(f\_{1-i}\cdot \wp (v_1))
  (f\_{2-i}\cdot \wp (v_2)) \cdots (f\_{1}\cdot \wp (v_{i-1}))\\
  &\times d_{f\_{0}}(v_i)\wp (v_{i+1})\cdots \wp (v_n),
\end{align*}
where $v_j\in V$ for all $j=1,\ldots ,n$, defines a family of
braided derivations $d:M\to \Hom (T(V),T)$ for $M$, $\wp $,
$T(V)$, and $T$. The uniqueness of $d$ as a family of braided
derivations follows from \eqref{derivadagral} and the fact that
$V$ generates the algebra $T(V)$. \epf

\begin{lema}\label{le:fambrderTVI}
  Let $M\in {}^H\m$, $V$ a vector space,
  $T$ an $H$-module algebra, $\wp : T(V)\to T$ an algebra map,
  and $d:M\to \Hom (T(V),T)$ a family of braided derivations. Let $I\subset T(V)$
  be an ideal with $\wp (I)=0$. Assume that $I$ is generated by
  a subset $J\subset I$, and define $R=T(V)/I$. The following
  are equivalent.
  \begin{enumerate}
    \item[(i)]
      $d$ induces a family of braided derivations $d^R:M\to \Hom (R,T)$ by
      letting $d^R_f(x+I):=d_f(x)$ for $x\in T(V)$, $f\in M$,
    \item[(ii)] $d_f(I)=0$ for all $f\in M$,
    \item[(iii)]
      $d_f(x)=0$ for all $f\in M$ and all generators $x\in J$.
  \end{enumerate}
\end{lema}

\pf The implications (i)$\Rightarrow $(ii) and (ii)$\Rightarrow
$(iii) are trivial. By \eqref{derivadagral}, the linearity of
$d_f$, and since $\wp (I)=0$, one obtains (iii)$\Rightarrow $(ii).
Finally, since $d^R:M\to \Hom (R,T)$ is a well-defined linear map,
for the implication (ii)$\Rightarrow $(i) it is sufficient to
check \eqref{derivadagral}. The latter holds since $I$ is an ideal
and $\wp (I)=0$. \epf

For the next theorem we need a compatibility relation between the maps
$\pl{g}$ and $\ad _c$.

\begin{lema}\label{le:plrel}
  Let $W\in \ydh $, $w\in \toba(W)$, $x\in W$, and $g\in W^*$. Then
\begin{align*}
  \pl{w\_{-1}\cdot g}( \ad _c (w\_{0})(x))
  = {}& \pl{w\_{-1}\cdot g}( w\_{0} )x \\
  &-(w\_{-2}\cdot x\_0)\pl{w\_{-1}\Ss^{-1}(x\_{-1})\cdot g}(w\_{0}).
\end{align*}
\end{lema}

\pf The definition of $\ad _c$ and \eqref{eqn:leibniz-pl-prim}
imply that {\allowdisplaybreaks
\begin{align*}
  &\pl{w\_{-1}\cdot g}( \ad _c (w\_{0})(x))
  = \pl{w\_{-1}\cdot g}( w\_{0}x)
  -\pl{w\_{-2}\cdot g}( (w\_{-1}\cdot x)w\_{0})\\
  &\quad = \pl{w\_{-1}\cdot g}( w\_{0})x
  -\pl{w\_{-2}\cdot g}( w\_{-1}\cdot x)w\_{0}
  +w\_{0}\pl{\Ss ^{-1}(w\_{-1})w\_{-2}\cdot g}(x)\\
  &\qquad
  -(w\_{-2}\cdot x\_{0})\pl{\Ss^{-1}(w\_{-3}x\_{-1}\Ss (w\_{-1}))w\_{-4}\cdot
  g}(w\_{0})\\
  &\quad = \pl{w\_{-1}\cdot g}( w\_{0} )x
  -\langle w\_{-2}\cdot g,w\_{-1}\cdot x\rangle w\_{0}\\
  &\qquad +w\langle g,x\rangle
  -(w\_{-2}\cdot x\_0)\pl{w\_{-1}\Ss^{-1}(x\_{-1})\cdot g}(w\_{0}).
\end{align*}
} The claim of the lemma now follows from \eqref{vd1}. \epf

We now show a very general way of constructing a family of braided
derivations of $\toba(W)\#\toba(W^*)$. This will be crucial in the
proof of Theorem \ref{theo:main}, but it may be of independent
interest. Recall the notion of canonical family of braided
derivations $d^W$, see Definition \ref{de:canfam}.

\begin{theorem}\label{theo:skew-der-qdo}
  Let $W\in \ydh$ be \fd.
  For all $w\in W$ and $f=\phi _W(w) \in W^{**}$,
  see \eqref{eqn:double-dualbis},
  define $d_f\in \End (\toba (W)\#\toba (W^*))$ by

\begin{align}\label{eqn:def-dv}
  d_f(x\# g) &=  -\ad_c w(x) \# g + (f\_{-1}\cdot x) \# \, d^{W^*}_{f\_0}(g)
\end{align}

\noindent for all $x\in \toba(W)$ and $g\in \toba(W^*)$. Then
$d:W^{**}\to \End (\toba(W)\#\toba(W^*))$ is a family of braided
derivations of $\toba(W)\#\toba(W^*)$.
\end{theorem}

\pf By Proposition \ref{prop:qdo-toba} and Definition
\ref{def:ksmash} there exists a unique algebra map $\wp :T(W\oplus
W^*)\to \toba (W)\#\toba (W^*)$ with $\wp (w)=w$, $\wp (g)=g$ for
all $w\in W$, $g\in W^*$. Let $d':W^{**}\to \Hom (T(W\oplus
W^*),\toba(W)\#\toba(W^*))$ be the unique family of braided
derivations with this $\wp $ and with
\begin{gather*}
  d'_f(x)=-\ad _c (\psi _W(f))(x),\qquad
  d'_f(g)=\langle f,g\rangle ,
\end{gather*}
where $f\in W^{**}$, $x\in W$, and $g\in W^*$, see Lemma~\ref{le:fambrderTV}.
We are going to use
the implication (iii)$\Rightarrow $(i)
in Lemma~\ref{le:fambrderTVI} to show that $d'$ induces the family $d$ of
braided derivations of $\toba(W)\#\toba(W^*)$. Indeed, one has
$d_f(z)=d'_f(z)$ for $z\in W\oplus W^*$. Further,
for $w=\psi _W(f)$ the map $-\ad _c w\in \End \toba(W)$
satisfies the relation
\begin{align*}
  &-\ad _c w (xy) = -wxy+(w\_{-1}\cdot (xy))w\_0\\
  &\quad = -wxy + (w\_{-1}\cdot x)w\_0 y
  - (w\_{-1}\cdot x)w\_0 y +(w\_{-2}\cdot x)(w\_{-1}\cdot y)w\_0\\
  &\quad =-\ad _c w(x)\, y-(w\_{-1}\cdot x)\,\ad _c (w\_0)(y).
\end{align*}
Thus, since $\psi _W$ is an $H$-comodule map, the restriction of
$d'_f$ to $T(W)$ coincides with $-\ad _c w\circ \pi _W$,
where $\pi _W:T(W)\to \toba(W)$ is the canonical map. Moreover,
the restriction of $d'_f$ to $T(W^*)$ is precisely $d^{W^*}_f\circ
\pi _{W^*}$. It remains to
show that $d'$ induces a family of braided derivations of
$\toba(W)\#\toba(W^*)$. By the
previous claims the latter family then has to coincide with the family $d$ of
linear maps $d_f$, where $f\in W^{**}$, and hence $d$
is a family of braided derivations.

To see that $d'$ induces a family of braided derivations of
$\toba(W)\#\toba(W^*)$, we have to check
that $d'_f$ vanishes on the generators (i)--(iii) in
Prop.~\ref{prop:qdo-toba}.
Since the restriction of $d'_f$ to $T(W)$ coincides with
$-\ad _c (\psi (f))\circ \pi _W$, one gets $d'_f(z)=0$ for all $z\in \J(W)$.
Similarly one has $d'_f(z)=d_f(\pi _{W^*}(z))$ for all $z\in \J(W^*)$.
Thus it suffices to check that
\begin{align}\label{eq:d'frel}
  d'_f(g\ot x-x\_{0}\ot (\Ss ^{-1}(x\_{-1})\cdot g)-\pl{g}(x))=0
\end{align}
for all $x\in W$, $g\in W^*$, and $f\in W^{**}$. Note that $d'_f(\pl{g}(x))=0$
since $\pl{g}(x)\in \ku $.
Let now $x,w\in W$, $g\in W^*$, and $f=\phi _W(w)\in W^{**}$. By definition of
$d'_f$ one gets
{\allowdisplaybreaks
\begin{align*}
  &d'_f(g\ot x-x_{(0)}\ot (\Ss^{-1}(x\_{-1})\cdot g))
  =d'_f(g)x+(f\_{-1}\cdot g)d'_{f\_{0}}(x)\\
  &\quad -d'_f(x\_{0})(\Ss^{-1}(x\_{-1})\cdot g)
  -(f\_{-1}\cdot x\_{0})d'_{f\_{0}}(\Ss^{-1}(x\_{-1})\cdot g)\\
  &\quad =\langle f,g\rangle x-(w\_{-1}\cdot g)\,\ad _c (w\_{0})(x)
  +\ad _c w(x\_0)\,(\Ss^{-1}(x\_{-1})\cdot g)\\
  &\qquad -(w\_{-1}\cdot x\_0)\langle \phi _W(w\_0),\Ss^{-1}(x\_{-1})\cdot g
  \rangle .\\
  \intertext{Now \eqref{eqn:rel-qdogral} and Lemma~\ref{le:plrel}
  allow to simplify this expression further:}
  &\quad =\langle f,g\rangle x-\ad _c (w\_{0})(x\_{0})\,
  (\Ss^{-1}(w\_{-1}x\_{-1})w\_{-2}\cdot g)\\
  &\qquad -\pl{w\_{-1}\cdot g}( \ad _c (w\_0)(x))
  +\ad _c w(x\_0)\,(\Ss^{-1}(x\_{-1})\cdot g)\\
  &\qquad -(w\_{-1}\cdot x\_0)\langle \phi _W(w\_0),\Ss^{-1}(x\_{-1})\cdot g
  \rangle \\
  &\quad =\langle f,g\rangle x
  -\langle w\_{-1}\cdot g,w\_0\rangle x
  +(w\_{-2}\cdot x\_0)\langle w\_{-1}\Ss^{-1}(x\_{-1})\cdot g,w\_{0}\rangle \\
  &\qquad -(w\_{-1}\cdot x\_0)\langle \phi _W(w\_0),\Ss^{-1}(x\_{-1})\cdot g
  \rangle .
\end{align*}
Using the relation $f=\phi _W(w)$ and \eqref{eq:phiVv,f} twice,
the latter expression becomes zero. This proves \eqref{eq:d'frel}.
} \epf

\section{Reflections of Nichols algebras}\label{section:weyl}

This section is devoted to the construction of ``reflections'',
see \eqref{eq:rfliM}. Based on them we introduce and study new
invariants of Nichols algebras in $\ydh $, see Definition
\ref{def:rsys}. Then we discuss the particular class of standard
semisimple Yetter-Drinfeld modules.

\subsection{Braided Hopf algebras with projection}
\label{subsection:braidedproj}
We begin by considering a commutative diagram of braided Hopf
algebras in $\ydh$:
$$
\xymatrix{& &
R\ar@{=}[1,0]\ar@{_{(}->}[1,-2]_{\iota}\\S\ar[0,2]^{\pi_R}& & R. }
$$
Here and below we use subscripts to distinguish between the
various projections, coactions, etc. By bosonization, we get a
commutative diagram of Hopf algebras:
$$
\xymatrix{& & R\# H\ar@{=}[1,0]\ar@{_{(}->}[1,-2]_{\iota}
 \\S \# H\ar[0,2]_{\pi_{R \# H}}\ar[1,2]_{\pi_{H,S}}& &R \# H
\ar[1,0]^{\pi_{H,R}}\\ && H.}
$$
Clearly, the projections $\pi_{H,R}: R \# H \to H$ and $\pi_{H,S}:
S \# H \to H$ satisfy
\begin{equation}\label{pis0}
\pi_{H,R}\pi_{R \# H} = \pi_{H,S}.
\end{equation}
We propose to study this situation through the subalgebra of
coinvariants
\begin{equation}\label{K0}
K := (S\# H)^{\co R\# H}.
\end{equation}
We collect some basic properties of $K$.

\begin{lema}\label{lema:basicK0}
\begin{enumerate}
    \item[(i)] $K$ is a braided Hopf algebra in $\ydrh$ and the
    multiplication induces an isomorphism
    $$
K\# (R\# H) \simeq S\# H.
    $$
    \item[(ii)]$K = S^{\co R} = \{x\in S\,|\,x\^{1}\otimes
\pi_{R}(x\^{2}) = x\otimes 1\}$ is a subalgebra of $S$ and the
multiplication induces an algebra isomorphism
    $$
K\# R \simeq S, \qquad \text{\emph{cf.} Remark \ref{smash-bos}.}
    $$
    \item[(iii)] $K$ is a Yetter-Drinfeld submodule over $H$ of $S$  and
\begin{equation}\label{deltah=deltarh}
\delta_H(x) = (\pi_{H,R}\ot \id)\delta_{R\# H}(x), \qquad x\in K.
\end{equation}
    \item[(iv)] $\Ss_S(K)$ is a subalgebra and Yetter-Drinfeld submodule over $H$ of
    $S$.
\end{enumerate}
\end{lema}

\pf (i). By the general theory of biproducts.

(ii). Let $x\in K$. By \eqref{pis0}, $x\_{1} \ot \pi_{H,
S}(x\_{2}) = x\_{1} \ot \pi_{H, R}\pi_{R\# H}(x\_{2}) = x\ot 1$;
hence $x\in S$. Now, $$x\ot 1 = x\_{1} \ot \pi_{R\# H}(x\_{2}) =
x^{(1)} (x^{(2)})_{(-1)} \otimes \pi_{R}((x^{(2)})_{(0)})$$ by
\eqref{smash3}. Applying the $H$-coaction to the second tensorand
and then $(\mu_S\ot \id)(\id\ot \Ss\ot \id)$, we get
\begin{align*}
x\ot 1 &= x\^{1}(x\^{2})\_{-2}
\Ss((x\^{2})\_{-1})\ot\pi_{R}((x\^{2})\_{0})
=x\^{1}\otimes\pi_{R}(x\^{2}), \end{align*} since $\pi_{R}$ is
$H$-colinear. Thus $x\in S^{\co R}$.

Conversely, let $x\in S^{\co R}$. Applying the $H$-coaction to the
second tensorand of the equality $x\^{1}\otimes \pi_{R}(x\^{2}) =
x\otimes 1$, and since $\pi_{R}$ is $H$-colinear, we get $x\ot 1 =
x^{(1)} (x^{(2)})_{(-1)} \otimes \pi_{R}((x^{(2)})_{(0)}) = x\_{1}
\ot \pi_{R\# H}(x\_{2})$. Hence $x\in K$. The multiplication gives
rise to an isomorphism because of the analogous fact in (i).

(iii). Clearly, $K$ is an $H$-submodule of $S$. From
\eqref{smash2} and \eqref{pis0} we get \eqref{deltah=deltarh}.
Thus $K$ is also an $H$-subcomodule, and \emph{a fortiori} a
Yetter-Drinfeld submodule, of $S$.

(iv) follows from (iii) and the properties of the antipode,
\emph{cf.} \eqref{eqn:antipodatrenzada-anti}. \epf

\subsection{The algebra $\K$}\label{subsection:basicK}
We next work in the following
general setting. Let $V$, $W$ be Yetter-Drinfeld modules over $H$
such that $V$ is a direct summand of $W$ in $\ydh$. Or, in other
words, we have a commutative diagram in $\ydh$:
$$
\xymatrix{& &
V\ar@{=}[1,0]\ar@{_{(}->}[1,-2]_{\iota}\\W\ar[0,2]^{\pi}& & V. }
$$
Set $\VP = \ker \pi$, so that $W = V\oplus \VP$ in $\ydh$.  By
functoriality of the Nichols algebra, we have a commutative
diagram of graded Hopf algebras in $\ydh$:
$$
\xymatrix{& &
\toba(V)\ar@{=}[1,0]\ar@{_{(}->}[1,-2]_{\iota}\\\toba(W)\ar[0,2]_{\pi_{\toba(V)}}&
& \toba(V) .}
$$
By bosonization, we get a commutative diagram of graded Hopf
algebras:

\begin{equation}\label{eqn:commdiag-avaw}
\xymatrix{& & \qquad \ac(V) = \toba(V) \#
H\ar@{=}[1,0]\ar@{_{(}->}[1,-2]_{\iota}
 \\\ac(W) = \toba(W) \# H\ar[0,2]_{\pi_{\ac(V)}}\ar[1,2]_{\pi_{H,W}}& &\ac(V)
 = \toba(V)  \# H\ar[1,0]^{\pi_{H,V}}
 \\ && H.}
\end{equation}
As before, the projections $\pi_{H,V}: \ac(V) \to H$ and
$\pi_{H,W}: \ac(W) \to H$ satisfy
\begin{equation}\label{pis}
\pi_{H,V}\pi_{\ac(V)} = \pi_{H,W}.
\end{equation}

The main actor of this section is the subalgebra of coinvariants
\begin{equation}\label{K}
\K := \ac(W)^{\co \ac(V)}.
\end{equation}

\begin{lema}\label{basicK}
\begin{enumerate}
    \item[(i)] $\K$ is a graded braided Hopf algebra in $\ydvh$ and the
    multiplication induces an isomorphism
    $$
    \K\# \ac(V) \simeq \ac(W).
    $$
    \item[(ii)] $\K = \toba(W)^{\co \toba(V)} = \{x\in \toba(W):
      x\^{1}\otimes \pi_{\toba(V)}(x\^{2}) = x\otimes 1\}$ is a graded
      subalgebra of $\toba(W)$ and the multiplication induces a homogeneous
      isomorphism
    $$
    \K\# \toba(V) \simeq \toba(W).
    $$
    \item[(iii)] $\K$ is a Yetter-Drinfeld submodule over $H$ of $\toba(W)$
      and
\begin{equation}\label{deltah=deltaav}
  \delta_H(x) = (\pi_{H,V}\ot \id)\delta_{\ac(V)}(x), \qquad x\in \K.
\end{equation}
    \item[(iv)] $\K\cap W = \VP \subset \Pc(\K)$.
\end{enumerate}
\end{lema}

\pf (i) to (iii) are consequences of Lemma~\ref{lema:basicK0}
except for statements ``$\K$ is graded'', that follow since
$\pi_{\ac(V)}$ is homogeneous.

(iv). If $x\in W$, then $x\^{1}\otimes \pi_{\toba(V)}(x\^{2}) =
x\ot 1 + 1\ot \pi_{\toba(V)}(x)$. Hence $x\in W\cap \K$ if and
only if $x\in \ker  \pi _{\toba(V)}\cap W= \VP$. Moreover, if
$x\in \VP$, then $\vartheta_{\K}(x) = x$, thus $\Delta_{\K}(x) =
x\ot 1 + 1\ot x$, cf. \eqref{smash2}. \epf

\subsection{The module $L$}\label{subsection:basicL}
We keep the notation of Subsection \ref{basicK}. Let $U$ be a
Yetter-Drinfeld submodule over $H$ of $\VP$. We define
\begin{equation}\label{defL}
L := \ad \toba(V) (U).
\end{equation}
In other words, $L$ is the vector subspace of $\ac(W)$ spanned by
the elements
\begin{equation}\label{defLm}
\ad_c (x_1)(\dots (\ad_c (x_m)(y))), \qquad x_h\in V,\, 1\le h\le m, \,
y \in U,
\end{equation}
for $m \ge 0$. We collect some basic properties of $L$.

\begin{lema}\label{basicL}
\begin{enumerate}
    \item[(i)] $L = \ad \ac(V) (U)$.
    \item[(ii)] $L = \oplus_{m\in \N} L^m$, where
      $L^m = L\cap \toba^m(W)$; $L^1 = U$.
    \item[(iii)] $L$ is a graded Yetter-Drinfeld submodule over $\ac(V)$ of $\Pc(\K)$.
    \item[(iv)] $L$ is a graded Yetter-Drinfeld submodule over $H$ of $\Pc(\K)$.
\item[(v)] For any $x\in L$, we have
\begin{align}\label{crucial}
\Delta_{\ac(W)} (x) &\in x\ot 1 + \ac(V) \ot L,\\
\label{crucial2} \Delta_{\toba(W)} (x) &\in x\ot 1 + \toba(V) \ot
L.
\end{align}
\item[(vi)] If $x\in L$ and $\pi_{\ac(V)}(x\_{1}) \ot \pro_1(x\_{2}) = 0$, then $x=0$.
\item[(vii)]
  If $0\neq L'$ is an $\ac(V)$-subcomodule of $L$, then $L'\cap U
\neq 0$.
\end{enumerate}
\end{lema}

\pf (i) follows from $\ad \ac(V) (U) = \ad \toba(V) \ad H (U)
\subset \ad \toba(V) (U)$.

(ii). It is clear that $L$ is a graded subspace of $\toba(W)$
since $\toba(V)$ is graded and $U$ is homogeneous. Indeed,
for all $m\in \N _0$ the space $L^{m+1}$
is the span of the elements \eqref{defLm}.

(iii). We know that $U \subset \Pc(\K)$ by Lemma~\ref{basicK}
(iv). Hence $L \subset \Pc(\K)$ by Remark \ref{prim-ydh}. We show
that $U$ is also an $\ac(V)$-subcomodule. If $y\in U$, then
\begin{align*}
\delta_{\ac(V)}(y) &= (\pi_{\ac(V)}\ot \id)(y\ot 1 + y\_{-1}\ot
y\_{0}) = y\_{-1}\ot y\_{0},
\end{align*}
because $\pi_{\ac(V)}(y) = 0$ (since $y\in \VP = \ker \pi$) and
$\pi_{\ac(V)}(y\_{-1}) = y\_{-1}$ (since $y\_{-1}\in H$). By (i)
and Remark \ref{ydsubm} (ii), $L$ is a Yetter-Drinfeld submodule
over $\ac(V)$ of $\Pc(\K)$. Finally, $L^m = L\cap \K^m$, being the
intersection of two Yetter-Drinfeld submodules, is a
Yetter-Drinfeld submodule itself.

(iv) follows from (iii) and \eqref{deltah=deltaav}. We prove
\eqref{crucial}: If $x = \ad z(y)$, where $z\in \toba(V)$ and
$y\in U$, then
\begin{align*}
\Delta_{\ac(W)}(x) &= z\_{1}y\_{1}\Ss(z\_{4}) \ot
z\_{2}y\_{2}\Ss(z\_{3}) \\
&= z\_{1}y\Ss(z\_{4}) \ot z\_{2}\Ss(z\_{3}) +
z\_{1}y\_{-1}\Ss(z\_{3}) \ot \ad (z\_{2})(y\_{0})
\\ &\in x\ot 1 +  \ac(V) \ot L,
\end{align*}
since $z\_{1}y\_{-1}\Ss(z\_{3}) \in \toba(V)\#H$ and
$ \ad (z\_{2})(y\_{0}) \in L$. Here again we used that $y\_{1} \ot y\_{2}
= y\ot 1 + y\_{-1} \ot y\_{0}$. Now
\begin{align*}
\Delta_{\toba(W)}(x) &= (\vartheta_{\toba(W)} \ot
\id)\Delta_{\ac(W)}(x) \\ &\in \vartheta_{\toba(W)}(x)\ot 1 +
\vartheta_{\toba(W)}(\ac(V)) \ot L \\ &= x\ot 1 + \toba(V) \ot L.
\end{align*}
by \eqref{proptheta}, showing \eqref{crucial2}.

(vi).  By \eqref{crucial}, for some
$y_i\in \ac(V)$,
$\ell_i \in
L$, we have
\begin{align*}
0 &= \pi_{\ac(V)}(x\_{1}) \ot \pro_1(x\_{2}) = \pi_{\ac(V)}(x)\ot
\pro_1(1) + \sum_i \pi_{\ac(V)}(y_i) \ot \pro_1(\ell_i)  \\&=
\sum_i y_i \ot \pro_1(\ell_i) = x\_{1} \ot \pro_1(x\_{2}) =
x\^{1}(x\^{2})\_{-1}\ot \pro_1\left((x\^{2})\_{0}\right).
\end{align*}
As the projection $\pro_1$ is $H$-colinear, we infer that
\begin{align*}
0 &= x\^{1}(x\^{2})\_{-2}\ot (x\^{2})\_{-1}\ot
\pro_1\left((x\^{2})\_{0}\right)\\
&\overset{\text{applying }(\mu\ot\id)(\id\ot\Ss\ot\id)}{\implies}
\quad x\^{1} \ot \pro_1(x\^{2}) = 0.
\end{align*}
Since $x\in L
\subset \sum_{n\ge 1} \toba^n(W)$, we conclude that $x=0$ by
Lemma~\ref{rmk:gradedcondition}.

(vii). Let $0\neq x\in L'$ and write $x = \sum_{1\le m \le p}
x(m)$ with $x(m)\in L^m$ and $y:= x(p) \neq 0$. By (vi),
$$
0\neq \pi_{\ac(V)}(y\_{1}) \ot \pro_1(y\_{2}) \in \ac^{p-1}(V)\ot
\toba^1(W).
$$
Let now $F\in \Hom(\ac(V), \ku )$ such that the restriction of $F$
to $\ac^m(V)$ is 0 for all $m\neq p-1$. We claim that
$$
F\pi_{\ac(V)}(x\_{1})x\_{2} = F\pi_{\ac(V)}(y\_{1})
\pro_1(y\_{2}).
$$
Indeed,
\begin{align*}
F\pi_{\ac(V)}(x\_{1})x\_{2} &=\sum_{1\le m \le p}
F\pi_{\ac(V)}(x(m)\_{1})x(m)\_{2}
\\&=F\pi_{\ac(V)}(y\_{1})y\_{2}
\\ &=F\pi_{\ac(V)}(y\_{1}) \pro_1(y\_{2}).
\end{align*}
Here the second and third equalities are clear from the assumption
on $F$; if $m\le p$ then $\pi_{\ac(V)}(x(m)\_{1})\ot x(m)\_{2} \in
\oplus_{0\le h \le m}\ac^{m-h}(V) \ot \ac^{h}(W)$. Applying $F$ we
get 0 except $m = p$, $h=1$. Choosing $F$ appropriately, we have
$$
0\neq F\pi_{\ac(V)}(x\_{1})x\_{2} = F\pi_{\ac(V)}(y\_{1})
\pro_1(y\_{2})\in L'\cap U.
$$
\epf

Part (vii) of Lemma~\ref{basicL} implies some strong restrictions
on the Yetter-Drinfeld submodules of $L$.

\begin{prop}\label{ldirectsum} Assume that $U = U_1 \oplus\cdots \oplus U_\theta$ in
$\ydh$. Let $L_i = \ad \ac(V) (U_i)$. Then $L= L_1 \oplus\cdots
\oplus L_\theta$  in $\ydvh$. \end{prop}

\pf We have to show that the sum $L_1 +\cdots + L_\theta$ is
direct. Suppose that $L_i \cap (\sum_{j\neq i} L_j) \neq 0$; then
$L_i \cap (\sum_{j\neq i} L_j) \cap U \neq 0$ by Lemma~\ref{basicL}(vii). Note
that $(\sum_{j\neq i} L_j) \cap U =
(\sum_{j\neq i} L_j)^1 \cap U = \sum_{j\neq i} U_j$. Thus $L_i
\cap (\sum_{j\neq i} L_j) \cap U = U_i \cap (\sum_{j\neq i} U_j)
\neq 0$, a contradiction.  \epf

Clearly, if $U' \subsetneq U$ in $\ydh$, then $L' := \ad\ac(V) (U')
\subsetneq L$ in $\ydvh$. Hence, if $L$ is irreducible in $\ydvh$,
then   $U$ is irreducible in $\ydh$. The converse holds because of
Lemma~\ref{basicL}(vii).

\begin{prop}\label{lirreducible} If  $U$ is irreducible in $\ydh$,
then $L$ is irreducible in $\ydvh$.
\end{prop}

\pf Let $0\neq L'$ be a subobject of $L$ in $\ydvh$. Then $L'\cap
U \neq 0$ by Lemma~\ref{basicL}(vii). Since both $L'$ and $U$ are
$H$-stable, $L'\cap U$ is an $H$-submodule of $U$. It is an
$H$-subcomodule of $U$ by \eqref{deltah=deltaav}; thus $L'\cap U
\hookrightarrow U$ in $\ydh$. By the irreducibility assumption,
$L'\cap U = U$, hence $L = \ad \ac(V) (U) \subseteq L'$. \epf

If $U = \VP$, then we have the following property,
important for our later considerations.

\begin{prop}\label{Kgrnlj} The algebra $\K$ is generated
  by $\ad \toba(V) (\VP)$.
\end{prop}

\pf Let $\K'$ be the subalgebra of $\K$ generated by $\ad \toba(V) (\VP)$
and let
$X$ be the image of $\K'\# \toba(V)$ under the isomorphism $\K\#
\toba(V) \simeq \toba(W)$ given by multiplication. It suffices to
prove that $X = \toba(W)$. Since $V \subset X$ and $\VP \subset
X$, one gets $W\subset X$;
it remains then to show that $X$ is a subalgebra
of $\toba(W)$. For this, observe that $\K'$ is stable under the
adjoint action of $\ac(V)$. Indeed, $\ad x(yy') = \ad (x\_{1}) (y)
\ad (x\_{2})(y')$, for all $x\in \ac(V)$, $y, y' \in
\ad \toba(V) (\VP )$. Hence, if
$x\in V$ and $y\in \K'$, then $xy =  \ad_c x(y) + (x\_{-1}\cdot
y)x\_0 \in \K' + \K'\# V \subset X$. As both $\K'$ and $\toba(V)$
are subalgebras, we conclude that $X$ is a subalgebra and the
proposition follows. \epf

We now introduce the following finiteness condition on $U$. Recall that
$L=\ad \toba(V)(U)$.
\begin{equation}\label{Lmax}
  \text{$L^M\neq 0$ and $L^p = 0$ for some $M\in \mathbb{N}$
  and all $p>M$. }
\tag{F}
\end{equation}
Clearly, a sufficient condition for \eqref{Lmax} is that $L =
\oplus_{m\in \N} L^m$ has finite dimension. In this case, $\dim U<
\infty$ too.

If $M$ is determined by \eqref{Lmax}, then we write
\begin{align}
  \lmax := L^M.
  \label{eq:lmax}
\end{align}

\begin{lema}\label{basicLt}
Assume that $U$ satisfies Condition (F). Let $Z$ be a
Yetter-Drinfeld submodule over $H$ of $\lmax$ and
$\lt$ the $\toba(V)$-subcomodule of $L$ generated
by $Z$.
\begin{enumerate}
  \item[(i)] $\lt = \oplus_{m=1}^M \lt^m$, where
      $\lt^m = \lt \cap \toba^m(W)$ for all $m$, and $\lt^M = Z$.
    \item[(ii)] $\lt$ is the $\ac(V)$-subcomodule of $L$ generated by $Z$.
    \item[(iii)] $\lt$ is a graded Yetter-Drinfeld submodule over $\ac(V)$ of $L$.
    \item[(iv)] $\lt$ is a graded Yetter-Drinfeld submodule over $H$ of $L$.
  \end{enumerate}
\end{lema}

\pf
By \eqref{subcom-gen-gr}, $\lt$ is the vector subspace of
$L$ spanned by the elements
\begin{equation*}
  \langle f, z\_{-1}\rangle \, z\_{0}, \quad \text{where } z\in Z, \,
f\in \toba^n(V)^*, \qquad n\ge 0,
\end{equation*}
where $z\_{-1}\ot z\_{0}=\delta _{\toba(V)}(z)$.
Let $z\in Z$ and
$f\in \toba^n(V)^*$. We obtain that $\langle f, z\_{-1}\rangle\, z\_{0} \in
\lt^{M-n}$, since
$$z\_{-1}\ot z\_{0} = \pi_{\toba(V)} (z\^{1})\ot
z\^{2}\in \oplus_{m\in \N_0} \toba^{m}(V) \ot \toba^{M-m}(W).$$
This proves (i). Now (ii) follows from Lemma~\ref{lema:smashcopr}
(iii); then (iii)
follows from (ii), Assumption (F), and Remark \ref{ydsubm} (i),
while (iv) follows from (iii) and \eqref{deltah=deltaav}. \epf

We can now present the first ingredient of our construction in
\eqref{eq:rfliM}.

\begin{theorem}\label{lmaxirreducible} Suppose that $U$ is irreducible
  in $\ydh$ and satisfies Condition~(F).
  Then $\lmax$ is irreducible in $\ydh$
  and $L$ is generated by $\lmax$ as a $\toba(V)$-comodule.
\end{theorem}

\pf By Proposition \ref{lirreducible}, $L$ is irreducible in
$\ydvh$. If $0\neq Z\hookrightarrow\lmax$ in $\ydh$, then $0\neq
\lt = \toba(V)^*\cdot Z\hookrightarrow L$ in $\ydvh$ by
Lemma~\ref{basicLt} (iii). Thus $\lt = L$, and $Z = \lt^M = L^M =
\lmax$ by Lemma~\ref{basicLt} (i). \epf

\subsection{Reflections}
\label{subsection:reflected}

\medbreak Let us fix $\theta \in \mathbb{N}$ and let
$\Ib =\{1,\ldots ,\theta \}$.
Let $\C_\theta $ denote the class of all families
$$M=(M_1,\dots ,M_\theta )$$
of \fd{} irreducible Yetter-Drinfeld modules $M_j\in \ydh $, where
$j\in \Ib $.
Two families $M,M'\in \C_\theta $ are called
\textit{isomorphic} if $M_j$ is
isomorphic to $M'_j$ in $\ydh $ for all $j\in \Ib $.
In this case we write $M\simeq M'$.

\medbreak Let $(\alpha _1,\ldots
,\alpha _\theta )$ be the standard basis of $\zt $.
Let $M =(M_1,\ldots ,M_\theta )\in \C_\theta $ and
\begin{align}
  W=\oplus _{j=1}^\theta M_j.
  \label{eq:Wdec}
\end{align}
Define a $\zt$-grading on $W$ by
$\deg M_j=\alpha _j$ for all $j\in \Ib $.
We fix $i\in\Ib $ and set
$$ V = M_i,\qquad
\VP = \bigoplus_{j\in \Ib,\, j\neq i}M_j.$$ Thus, we are
in the situation of Subsections \ref{subsection:basicK} and
\ref{subsection:basicL}. Let
\begin{equation}\label{defLj}
  L_j := \ad \toba(V) (M_j) \qquad \text{for $j\in \Ib\setminus
  \{i\}$.}
\end{equation}
Thus, $L_j$ is the vector subspace of $\toba(W)$ spanned by the
elements $$\ad_c (x_1)(\dots(\ad_c (x_m)(y))),\quad x_h\in
M_i,\, 1\le h\le m,\, y \in M_j,\, m \ge 0.$$
Recall that $\K =\ac(W)^{\co \ac(V)}=\toba (W)^{\co \toba (V)}$,
see \eqref{K} and Lemma~\ref{basicK} (ii). Consider the
$\zt$-grading on the algebras $\toba(W)$ and $\toba(V)$ discussed
in Remark \ref{obs:Dlgraded}, page~\pageref{obs:Dlgraded}. Then
the algebras $\ac(W)$ and $\ac(V)$ are also $\zt$-graded, by
setting $\deg H = 0$. Since the map $\pi_{\ac(V)}$ in
\eqref{eqn:commdiag-avaw} is homogeneous, the algebra $\K$
inherits this grading. Then $L_j$ is a $\zt$-graded subspace of
$\K$ and $\supp L_j \subset \alpha_j + \N_0 \alpha_i$. Let
\begin{equation}\label{eqn:defi-cartan-matrix}
  -a^M_{ij} := \sup \{h\in \N_0\,|\,\alpha_j + h \alpha_i\in \supp L_j\}.
\end{equation}
Then either $a^M_{ij}\in \Z_{\leq 0}$ (when $\supp L_j$ is finite),
or $a^M_{ij} = -\infty$. Let also $a^M_{ii} = 2$.

We introduce the following finiteness conditions for $M$.
\begin{enumerate}
    \item[$(F_i)$]$\dim L_j$ is finite for all $j\in \Ib$, $j\neq i$,
\end{enumerate}
or, equivalently,
\begin{enumerate}
    \item[$(F'_i)$]$\supp L_j$ is finite for all $j\in \Ib$, $j\neq i$.
\end{enumerate}
Note that $(F_i)$ means that $a^M_{i j}>-\infty $ for all $j\in
\Ib \setminus \{i\}$.

\begin{obs}\label{rmk:gkdim}
  It would be interesting to find an \emph{a priori} condition
guaranteeing that $(F_i)$ holds.
Obviously, if $\dim \toba(W)<\infty$, then $\dim L_j < \infty$ for
all $j$. Because of \cite{R2}, we believe that $(F_i)$ holds
whenever the Gelfand-Kirillov dimension of $\toba(W)$ is finite.
\end{obs}

Assume that $M$ satisfies Condition~$(F_i)$. Let $\sE \in
GL(\theta, \Z)$ and
\begin{align}
  \rfl _i(M):= (M'_1,\ldots ,M'_\theta )\in \C_\theta
  \label{eq:rfliM}
\end{align}
be given by
\begin{align}\label{eqn:betaj}
\sE(\alpha_j) &= \alpha_j - a^M_{ij}\alpha_i 
, \qquad j\in \Ib,
\\
\label{eqn:defRi}
M'_j &= \begin{cases}
  L_j^{\max }                & \text{if $j\neq i$,} \\
  {M_i}^{\,*} = V^* & \text{if $j=i$.}\end{cases}
\end{align}
Notice that $\rfl_i(M)$ is an object of $\C_{\theta }$ by Theorem
\ref{lmaxirreducible}. We say that $\rfl _i$ is \textit{the $i$-th
reflection}. The linear map $\sE $ is a reflection in the sense of
\cite[Ch.\,V, \S2.2]{B68},
that is, $\sE ^2=\id $ and the rank of $\id -\sE $ is $1$.

We embed $V^*$ into $W^*$ via the decomposition of $W$ in
\eqref{eq:Wdec}. Then
$$\K\# \toba({}V^*)\subset \toba(W)\#\toba(W^*)$$
is a subalgebra, see Definition \ref{def:ksmash}. Further, $\K \#
\toba(V^*)$ is a $\zt $-graded algebra in $\ydh $ with
$\deg x=\sE (\alpha _i)=-\alpha _i$ for
all $x\in V^*$, see Remark \ref{obs:Dlgraded}.

\begin{lema}
  The map $T(\K\oplus V^*)\to
  \K \# \toba(V^*)$, $\K\oplus V^*\ni (x,f)\mapsto x\#1+1\#f$, induces
  an algebra isomorphism $T(\K\oplus V^*)/I\to \K\#\toba(V^*)$, where $I$ is
  the two-sided ideal generated by
  \begin{itemize}
    \item[(i)] the elements $x\otimes y-xy$, where $x,y\in \K$, and $1_\K -
      1_{T(\K \oplus V^*)}$,
    \item[(ii)] the relations of $\toba(V^*)$,
    \item[(iii)] the elements
\begin{align}\label{eq:V*Krel} g\otimes x -
x\_{0}\otimes \Ss^{-1}(x\_{-1})\cdot g  -\pl{g}(x), \qquad
x\in \K , \, g\in V^*.
\end{align}
  \end{itemize}
  \label{le:KBgenrel}
\end{lema}

\pf See the proof of Proposition \ref{prop:qdo-toba}. \epf

Let $W'=\oplus _{j=1}^\theta M'_j$.
Then $W'$ is contained in $\K \# \toba (V^*)$ via the embeddings
\begin{align*}
  &M'_j\subset \K \simeq \K \# \ku \subset \K \# \toba (V^*) \qquad
  \text{for all $j\not=i$,}\\
  &M'_i=V^*\simeq \ku \# V^*\subset \K \# \toba (V^*)^+.
\end{align*}
Moreover, $W'$ inherits a $\zt $-grading from $\toba(W)\#\toba(W^*)$:
One has $\deg M'_j=\sE (\alpha _j)$ for all $j\in \Ib$.

\begin{lema}\label{lema:gen-by-Wprime}
The algebra $\K\# \toba(V^*)$ is generated by $W'$.
\end{lema}

\pf
Let $\Bg =\ku \langle W'\rangle $ be the subalgebra of $\K \# \toba (V^*)$
generated by $W'$.
Since $W'\in \ydh $, $\Bg$ is a subobject of $\toba(W)\#\toba(W^*)$
in $\ydh $ by Corollary~\ref{cor:qdo-ydh}.
Fix $j\neq
i$ and pick $x\in L_j\cap \Bg$, $f \in V^*$. Then $f x = x\_{0}\,
\Ss^{-1}(x\_{-1})\cdot f + \pl{f}(x)$ by \eqref{eqn:rel-qdogral}.
Now, $L_j\cap \Bg$ being a Yetter-Drinfeld submodule over $H$,
this says that $\pl{f}(x)\in \Bg$. But
$$
\pl{f}(x) = \langle f, x\^1\rangle x\^2 \in  \langle f, x\rangle 1
+ \langle f, \toba(V) \rangle L_j \subset L_j,
$$
by \eqref{crucial2}. This shows that $L_j\cap \Bg$ is a
$\toba(V)$-subcomodule of $L_j$; indeed, $\langle f,x\_{-1}\rangle
x\_{0} = \langle f,\pi_{\toba(V)} (x\^{1})\rangle x\^{2}=\langle
f, x\^1\rangle x\^2$. We conclude that $L_j\cap \Bg = L_j$ by
Lemma~\ref{basicLt} (iii) and Prop.~\ref{lirreducible}, since
$0\neq L_j\cap \Bg \supset L_j^{\max}$.
Hence $L_j\subset \Bg $ for $j\in \Ib \setminus \{i\}$,
and Prop.~\ref{Kgrnlj} implies that $\K\subset \Bg $.
This proves the lemma.
\epf

\medbreak Here is our first main result.

\begin{theorem}\label{theo:main}
  Let $M=(M_1,\ldots ,M_\theta )\in \C_\theta $
  and $i\in \Ib $ such that $M$ satisfies Condition~$(F_i)$.
  Let $V=M_i$, $W=\oplus _{j\in \Ib }M_j$,
  $\K =\toba(W)^{\co \toba(V)}$,
  $M'=\rfl_i(M)$ and $W'=\oplus _{j\in \Ib }M'_j$.
  Define a $\zt $-grading on $W'$ by $\deg x=\sE(\alpha _j)$
  for all $x\in M'_j$, $j\in \Ib $.
  \begin{enumerate}
    \item
  The
inclusion $W'\hookrightarrow \K\#\toba(V^*)$
  induces
a $\zt$-homogeneous isomorphism
$\toba(W') \simeq \K\#\toba(V^*)$
of algebras
and of Yetter-Drinfeld modules over $H$.
\item The family $\rfl_i(M)$ satisfies Condition $(F_i)$, and
  $s_{i,\rfl _i(M)}=\sE $, $\rfl _i^2(M)\simeq M$.
  \end{enumerate}

\end{theorem}

We prove the theorem in several steps.
The strategy of the proof is the following.
First we define a surjective algebra map $\Omega :\K\#\toba(V^*)\to \toba
(W')$. Then we conclude that the same construction can be performed for $M'$
instead of $M$, and that (2) holds.
Finally we prove that $\Omega $ is bijective.
The restriction of the
inverse map of $\Omega $ to $W'$ is the given embedding of $W'$ in $\K \#
\toba (V^*)$.

For the definition of $\Omega $, see Prop.~\ref{pr:defOm}, we use
the characterization of Nichols algebras in Theorem
\ref{theo:quotients-with-derivations} (ii). In the next lemma we
prove the existence of the required family of braided derivations.

\begin{mlema}\label{le:exderiv}
There is a unique family $d:W'{}^{*}\to \End (\K \# \toba (V^*))$
of braided derivations of $\K \# \toba (V^*)$
such that
\begin{align}\label{eq:dfw'}
  d_f(w')=\langle f,w'\rangle
\end{align}
for all $f\in W'{}^* \simeq V^{**} \oplus
\oplus_{j\in \Ib \setminus \{i\}}(L_j^{\max})^*$, $w'\in W'$.
Moreover, for all $v\in V$, $f=\phi_V(v)\in V^{**}$,
and $x\in \K $ equation
$d_f(x)=-\ad_c v(x)$ holds.
\end{mlema}


\pf The family $d$ is unique since $\K\#\toba(V^*)$ is generated
by $W'$, see Lemma~\ref{lema:gen-by-Wprime}. By Definition
\ref{defi:skew-der} it is sufficient to show that
\begin{enumerate}
  \item there exists a family
$d:V^{**}\to \End (\K \# \toba (V^*))$
of braided derivations of $\K \# \toba (V^*)$ such that
$d_f(w')=\langle f,w'\rangle $ for all $f\in V^{**}$ and $w'\in W'$,
  \item for all $j\in \Ib \setminus \{i\}$ there exists a family
$d:(L_j^{\max})^*\to \End (\K \# \toba (V^*))$
of braided derivations of $\K \# \toba (V^*)$ such that
$d_f(w')=\langle f,w'\rangle $ for all $f\in (L_j^{\max})^*$ and $w'\in W'$.
\end{enumerate}

First we prove (1). Let $d:V^{**}\to \End (\toba(W)\#\toba(W^*))$
be the restriction to $V^{**}$ of the family of braided
derivations in Theorem \ref{theo:skew-der-qdo}. By
\eqref{eqn:def-dv} one gets
\begin{align}
  d_f(x)=-\ad _c v(x)\quad \text{for all $v\in V$, $f=\phi _V(v)\in V^{**}$,
  $x\in \K $.}
  \label{eq:dfx}
\end{align}
Thus $d_f(\K)\subset\K$ for all $f\in V^{**}$, and
$d_f(\toba(V^*))\subset \toba (V^*)$ since $d_f(w')=\langle
f,w'\rangle $ for all $w'\in V^*$ by \eqref{eqn:def-dv}. Hence $d$
induces a family of braided derivations of $\K \# \toba(V^*)$ by
restriction. The relation $d_f(w')=\langle f,w'\rangle =0$ for
$w'\in L_j^{\max}$, $j\not=i$, follows from the definition of
$L_j^{\max}$, and the second claim of the lemma holds by
\eqref{eq:dfx}.

\medbreak
To prove (2), let $j\in \Ib \setminus \{i\}$. We first define
a family
$d:(L_j^{\max})^*\to \End (\K )$
of braided derivations of $\K $.
Then we extend $d$ to a family of braided
derivations of $\K \# \toba(V^*)$.

Recall from \eqref{eqn:betaj} that $\sE (\alpha _j)=\alpha _j-a^M_{i j}\alpha
_i$. Define $d_F: \toba(W) \to \toba(W)$ for any $F\in
\toba(W^*)_{-\sE (\alpha _j)}$ by
\begin{equation}\label{eqn:dfLj}
d_F(x) := \langle F, (x\^2)\_0\rangle \Ss^{-1}((x\^2)\_{-1})\cdot
x\^1, \qquad x\in \toba(W),
\end{equation}
see \eqref{eq:plbxi}. Then
\begin{equation}\label{eqn:dfLj-seanula}
d_F(x) = 0 \qquad \text{if } x\in L_h, h\neq j,i, \text{ or }
x\in L_j^m,\ m < 1 - a^M_{ij}.
\end{equation}
Indeed, if $x\in L_h^m$, where $h\in \Ib \setminus \{i\}$, $m\in
\N $, then by Lemma~\ref{basicL} (iii) and \eqref{crucial2} one
gets
%
\begin{align*} \Delta(x)  \in  x\ot 1 + 1\ot x +
\sum_{0< r< m}
\toba(W)_{r\alpha_i}\ot\toba(W)_{\alpha_h + (m-1-r)\alpha_i}.
\end{align*}
Hence $\langle F, (x\^2)\_0\rangle = 0$ whenever $h\not=j$ or $h=j$, $m<1-a_{i
j}^M$. Further, if $x\in L_j^{1-a^M_{i j}}$ then
\begin{equation}\label{eqn:dfLj-max}
  d_F(x) = \langle F, x\rangle \qquad \text{for all } x\in
  L_j^{1-a^M_{ij}}.
\end{equation}
We next claim that
\begin{equation}\label{eqn:dfLj-K}
  d_F(xy) := d_F(x)y + (F\_{-1}\cdot x) d_{F\_0}(y) \qquad
  \text{for all }x, y\in \K.
\end{equation}
Let $x,y\in \K $. Then
\begin{equation}\label{eq:dFxy}
\begin{aligned}
d_F(xy) = &\langle F, (x\^2)\_0(y\^2)\_0\rangle \\
&\quad \times \Ss^{-1}((x\^2)\_{-1}(y\^2)\_{-1}) \cdot
[x\^1((x\^2)\_{-2}\cdot y\^1)].
\end{aligned}
\end{equation}
Now $\langle F, (x\^2)\_0(y\^2)\_0\rangle = \langle F\^1,
(y\^2)\_0\rangle\langle F\^2, (x\^2)\_0\rangle$.
Further,
\begin{align*}
  \Delta(F) - F\ot 1 - 1\ot F \in
  \sum _{0<r<1-a^M_{i j}}(&\toba(W^*)_{-r\alpha _i}\ot
  \toba(W^*)_{-\sE (\alpha _j)+r\alpha _i}\\
  &+\toba(W^*)_{-\sE (\alpha _j)+r\alpha _i}\ot \toba(W^*)_{-r\alpha _i}).
\end{align*}
Since $\K \subset \toba(W)$ is a left coideal and $\langle F',\K \rangle=0$
for all $F'\in \toba(W^*)_{-r\alpha _i}$ and $r>0$, one gets
$$\langle F, (x\^2)\_0(y\^2)\_0\rangle = \langle F,(y\^2)\_0\rangle
\cou ((x\^2)\_0)+\cou((y\^2)\_0)\langle F,(x\^2)\_0\rangle.$$ This
means that $d_F$ behaves in the same way as $d_{F'}$ for primitive
$F'$, and hence \eqref{eqn:dfLj-K} follows from
\eqref{plinverso-prim}.
%
%

We point out two consequences of the claim \eqref{eqn:dfLj-K}. First,
this shows that $d_F(\K) \subset \K$; indeed, $\K$ is generated as
an algebra by $L$ and we know already that $d_F(L) \subset \K$ by
\eqref{eqn:dfLj-seanula} and \eqref{eqn:dfLj-max}. Second, the
inclusion $L_j^{1-a^M_{ij}}\subset \toba(W)_{\sE (\alpha _j)}$ induces a
projection $\pi: \toba(W^*)_{-\sE (\alpha _j)} \to \big(L_j^{1-a^M_{ij}}\big)^*$;
then $d_F\in \End \K$ depends only on $f = \pi(F)$. For, if
$\pi(F) = 0$, then $d_F = 0$ on $L$ by \eqref{eqn:dfLj-seanula}
and \eqref{eqn:dfLj-max}. Hence $d_F = 0$ on $\K$ by
\eqref{eqn:dfLj-K}.
Thus we have constructed the desired family $d:(L_j^{\max})^*\to
\End \K $ of braided derivations of $\K $.

Now we extend $d$ to a family of braided derivations of $\K\#\toba(V^*)$
by letting
\begin{equation}\label{eqn:df-final}
d_f(x g) = d_f(x) g, \qquad x\in \K,\, g\in \toba(V^*).
\end{equation}
It is clear that $d_f(w')=\langle f,w'\rangle $ for all $f\in
(L_j^{\max})^*$, $w'\in W'$. It remains to prove that
\begin{align}
  d_f(bc) = (f\_{-1} \cdot b)\, d_{f\_{0}}(c) + d_{f}(b)\, c
  \label{eq:dfbc}
\end{align}
for all $b,c\in \K\#\toba(V^*)$ and $f = \pi(F)\in
\big(L_j^{1-a^M_{ij}}\big)^*$. Similarly to the proof of Theorem
\ref{theo:skew-der-qdo}, we use Lemma~\ref{le:fambrderTVI}
(iii)$\Rightarrow $(i) and Lemma~\ref{le:KBgenrel} to show that
$d:(L_j^{\max})^*\to \End ( \K\#\toba(V^*))$, given in
\eqref{eqn:df-final}, defines a family of braided derivations of
$\K\#\toba(V^*)$. Again it suffices to check that
\begin{align}
  d'_f(g\otimes x)=d'_f(x\_0\otimes \Ss^{-1}(x\_{-1})\cdot g  +\pl{g}(x))
  \label{eq:d'fgx}
\end{align}
for all $x\in \K$, $g\in V^*$, where $d':(L_j^{\max})^*\to \Hom(T(\K\oplus
V^*),\K\#\toba(V^*))$ denotes the family of braided derivations induced by
$d'_f|_\K=d_f$ and $d'_f|_{V^*}=0$.
The right-hand side of \eqref{eq:d'fgx} is
\begin{align*}
&d'_f\big( x\_{0}\otimes
\Ss^{-1}(x\_{-1})\cdot g + \langle g, x\^1\rangle x\^2\big)\\
&\quad = d_F( x\_{0})\, \Ss^{-1}(x\_{-1})\cdot g
+d_F( \langle g, x\^1\rangle x\^2)\\
&\quad = 
\langle F, (x\^2)\_0\rangle \big(\Ss^{-1}((x\^2)\_{-1})\cdot
(x\^1)\_0\big) \, \big( \Ss^{-1}((x\^1)\_{-1}(x\^2)\_{-2})\cdot g\big)
\\ &\qquad + 
\langle g, x\^1\rangle  \langle F, (x\^3)\_0\rangle
\Ss^{-1}((x\^3)\_{-1})\cdot x\^2,
\end{align*}
and the left-hand side is
\begin{align*}
(f\_{-1} \cdot g)& d'_{f\_{0}}(x) + d'_f(g) x =
(F\_{-1} \cdot g) d_{F\_{0}}(x) \displaybreak[0] \\
\overset{\eqref{eqn:dfLj}}{=}& (F\_{-1} \cdot g) \langle F\_0, (x\^2)\_0\rangle
\Ss^{-1}((x\^2)\_{-1})\cdot x\^1\\
\overset{\eqref{properties:duality4-tensor}}{=}& (\Ss^{-1}( (x\^2)\_{-1})\cdot g) \langle F, (x\^2)\_0\rangle
\Ss^{-1}((x\^2)\_{-2})\cdot x\^1\\
\overset{\phantom{\eqref{eqn:dfLj}}}{=}&
\langle F, (x\^2)\_0\rangle \Ss^{-1}( (x\^2)\_{-1})\cdot (gx\^1)\\
\overset{\eqref{eqn:rel-qdogral}}{=}& \langle F, (x\^2)\_0\rangle \Ss^{-1}( (x\^2)\_{-1})\cdot \big(
(x\^1)\_0 \Ss ^{-1}( (x\^1)\_{-1})\cdot g\big)\\
\displaybreak[0]
&+ \langle F, (x\^3)\_0\rangle \Ss^{-1}( (x\^3)\_{-1})\cdot
\langle g,x\^1\rangle x\^2.
\end{align*}
This proves \eqref{eq:d'fgx} and completes the proof of the lemma.
\epf

Recall the notation from Theorem~\ref{theo:main}.

\begin{prop}\label{pr:defOm}
  There exists a unique surjective algebra map
  $$\Omega :\K\# \toba(V^*)\to \toba(W')$$
  which is the identity on $W'$.
  Define a $\zt $-grading on $W'$ by $\deg x=\sE (\alpha _j)$ for all
  $x\in M'_j$, $j\in \Ib $.
  Then $\Omega $ is a $\zt $-graded map in $\ydh $, and for
  all $v\in V$, $f\in V^*$, $x\in \K$ the following equations
      hold.
      \begin{align}
        \label{eq:Omplf}
        \Omega (\pl{f}(x))=&\ad _{c^{-1}}f(\Omega (x)),\\
        \Omega ( \ad _c v(x))=&-d^{W'}_{\phi_V(v)}(\Omega (x)).
        \label{eq:Omadv}
      \end{align}
\end{prop}

\pf By Lemma~\ref{lema:gen-by-Wprime} there is a unique surjective
algebra map $T(W')\to \K\#\toba(V^*)$ which is the identity on
$W'$. Let $I$ be the kernel of this map. Since $\K\#\toba(V^*)$ is
an $H$-module, $I$ is invariant under the action of $H$. By Main
Lemma~\ref{le:exderiv} there is a unique family $d:W'{}^*\to \End
(\K \# \toba(V^*))$ of braided derivations satisfying
\eqref{eq:dfw'}. Thus Theorem
\ref{theo:quotients-with-derivations} (ii) applies, that is, the
algebra map $\Omega $ exists and is unique. By definition of the
$\zt $-gradings and the Yetter-Drinfeld structures, $\Omega $ is a
$\zt $-graded map in $\ydh $.

\eqref{eq:Omadv} follows from \eqref{factorizacion-omega} by using
the second part of Main Lemma~\ref{le:exderiv}, Equations
\eqref{eq:plbxi}, \eqref{plinverso-prim}, and Definition
\ref{de:canfam}. \eqref{eq:Omplf} follows from the formulas
\begin{align*}
  \Omega (\pl{f}(x))
  \overset{\eqref{eqn:rel-qdogral}}{=}&
  \Omega (fx-x\_0 (\Ss^{-1}(x\_{-1})\cdot f))\\
  \overset{\phantom{\eqref{eqn:rel-qdogral}}}{=}&
  f\Omega (x)-\Omega (x\_0) (\Ss^{-1}(x\_{-1})\cdot f)
  =\ad _{c^{-1}}f(\Omega (x)).
\end{align*}
\epf

%

Now we are prepared to prove Theorem \ref{theo:main}.

\pf[Proof of Theorem \ref{theo:main}] We follow the strategy
explained below Theorem \ref{theo:main}. Recall that $L_j=\ad
\toba(V)(M_{\alpha _j})$ and
\begin{equation}
\begin{aligned}
  L_j=\ku\text{-span of }\{\pl{f_1}\cdots \pl{f_n}(x)\,|\,&
  x\in L_j^{\max}=M'_j,\\
  &n\ge 0, f_1,\ldots ,f_n\in V^*\}
\end{aligned}
  \label{eq:Lj}
\end{equation}
for all $j\in \Ib \setminus \{i\}$ by Theorem
\ref{lmaxirreducible}. Let $M'=\rfl_i(M)$ as in \eqref{eq:rfliM},
\begin{align}
  L'_j=\ad \toba(V^*)(M'_{\beta _j})\subset \toba(W'),
  \label{eq:L'j}
\end{align}
and $\Omega :\K\#\toba(V^*)\to \toba(W')$ the epimorphism in
Prop.~\ref{pr:defOm}.

We first claim that
\begin{align}
  \widetilde{\Omega }|_{L_j}:L_j\to L'_j\quad \text{is bijective,}
  \label{eqn:isodemodulos}
\end{align}
where $\widetilde{\Omega }=\Ss _{\toba(W')}\Omega $.
Indeed, let $x\in L_j^{\max}$, $n\ge 0$, and $f_1,\ldots ,f_n\in V^*$.
Then
\begin{align*}
  &\Ss_{\toba(W')}\Omega(\pl{f_1}\cdots \pl{f_n}(x))
  \overset{\eqref{eq:Omplf}}{=}\Ss_{\toba(W')}\big( \ad _{c^{-1}}(f_1)(\Omega
  (\pl{f_2}\cdots \pl{f_n}(x)))\big)\\
  &\quad \overset{\eqref{S1}}{=}\ad _c (f_1)\big(\Ss _{\toba(W')}\Omega
  (\pl{f_2}\cdots \pl{f_n}(x))\big)\\
  &\quad \overset{\phantom{\eqref{S1}}}{=}
  \ad _c (f_1)\big(\cdots \big(\ad _c (f_n)\big(\Ss _{\toba(W')}\Omega (x)
  \big)\big)\big)\\
  &\quad \overset{\phantom{\eqref{S1}}}{=}
  \ad _c (f_1)\big(\cdots \big(\ad _c (f_n)\big(\Ss _{\toba(W')}(x)\big)
  \big)\big)\\
  &\quad \overset{\phantom{\eqref{S1}}}{=}
  -\ad _c (f_1)(\cdots (\ad _c (f_n)(x)))\in L'_j.
\end{align*}
Since $\widetilde{\Omega }(x)=-x$ for all $x\in
L_j^{\max}=M'_{\beta _j}$, \eqref{eq:Lj} and \eqref{eq:L'j} imply
that $\widetilde{\Omega }(L_j)=L'_j$. We now prove that $\ker
\widetilde{\Omega }\cap L_j=\ker \Omega \cap L_j$ is a
Yetter-Drinfeld module over $\ac(V)$. Together with the
irreducibility of $L_j$, see Prop.~\ref{lirreducible}, this
implies that $\widetilde{\Omega }$ is injective and hence
Claim~\eqref{eqn:isodemodulos} holds.

Since $\Omega $ is a map in $\ydh $, see Prop.~\ref{pr:defOm},
one obtains that $\ker \Omega \cap L_j$ is an object in $\ydh $.
Further, for all $x\in L_j\cap \ker \Omega $ one has
\begin{align*}
  \Omega (\toba(V)\cdot x)=0\quad \text{by \eqref{eq:Omadv}, and}\quad
  \Omega(V^*\cdot x)=0 \quad \text{by \eqref{eq:Omplf}.}
\end{align*}
Thus $L_j\cap \ker \Omega $ is an object in $\ydvh $, and
Claim~\eqref{eqn:isodemodulos} is proven.

Now we prove Theorem \ref{theo:main}(2). Since $\Omega $ and
$\Ss_{\toba(W')}$ are $\zt $-graded maps, \eqref{eqn:isodemodulos}
implies that $\supp L_j=\supp L'_j$ for all $j\in \Ib \setminus
\{i\}$. In particular, $\supp L'_j$ is finite for all $j\in \Ib
\setminus \{i\}$, that is, Condition~$(F_i)$ is fulfilled for
$M'=\rfl_i(M)$, and hence $M'':=\rfl_i(M')$ is well-defined.
For all $j\in \Ib $ let $\gamma _j=s_{i,M'}s_{i,M}(\alpha _j)$. Then
by Eq.~\eqref{eqn:betaj}
one obtains for all $j\in \Ib \setminus \{i\}$ the equations
\begin{align}\label{eqn:defi-cartan-matrix-prime}
  -a^{M'}_{ij} = &\sup \{h\in \N_0\,|\,\sE (\alpha_j) + h \sE (\alpha _i)
  \in \supp L'_j\}
  =-a^M_{i j}, \\
  \gamma_j = &\sE (\alpha_j) - a^{M'}_{ij}\sE (\alpha_i)=\alpha _j,\\
  M''_j = &L'_j\cap \toba(W')_{\gamma _j}\simeq
  L_j\cap \toba(W)_{\gamma _j}=M_j,
  \label{eq:M''}
\end{align}
where the last equation follows from the fact that
$\widetilde{\Omega }|_{L_j}:L_j\to L'_j$ is a $\zt $-graded
isomorphism in $\ydh $, see Prop.~\ref{pr:defOm} and
Claim~\eqref{eqn:isodemodulos}. Since $M''_i=(M'_i)^*=M_i^{**}
\simeq M_i$
by Remark \ref{lema:double-dual}, one obtains that
$\rfl_i(M')\simeq M$, that is, Theorem \ref{theo:main}(2) is
proven.

It remains to prove that $\Omega$ is an isomorphism. Let
$\K'=\toba(W')^{\co \toba(M'_{\beta _i})}$, $W''=\oplus _{j\in \Ib
}M''_{\gamma _j}$, and $\K''=\toba(W'')^{\co \toba(M''_{\gamma
_i})}$. Since $\K$ resp. $\K'$ is generated as an algebra by
$\oplus _{j\in \Ib \setminus \{i\}}L_j$ resp. $\oplus _{j\in \Ib
\setminus \{i\}}L'_j$, see Prop.~\ref{Kgrnlj}, we conclude from
Lemma~\ref{S} (ii) and Claim~\eqref{eqn:isodemodulos} that
$\widetilde{\Omega}(\K) = \K'$. By the same argument we have
$\widetilde{\Omega}'(\K') = \K''$, where
$\widetilde{\Omega}':\K'\#\toba(V^{**})\to \toba(W'')$ is the map
in Prop.~\ref{pr:defOm} obtained by starting with the family $M'$
instead of $M$. Thus $\widetilde{\Omega}$ and
$\widetilde{\Omega'}$ define surjective
$\mathbb{Z}^{\theta}$-homogeneous maps $$\K
\xrightarrow{\widetilde{\Omega}|_\K} \K'
\xrightarrow{\widetilde{\Omega'}|_{\K'}} \K''$$ of Yetter-Drinfeld
modules over $H$. But $\K \simeq \K''$ as
$\mathbb{Z}^{\theta}$-graded Yetter-Drinfeld modules since $W
\simeq W''$ by Theorem \ref{theo:main} (2). The
$\mathbb{Z}^{\theta}$-homogeneous components of $\K$ are all \fd{}
since $W$ is \fd. Hence the map $\K
\xrightarrow{\widetilde{\Omega}|_\K} \K'
\xrightarrow{\widetilde{\Omega'}|_{\K'}} \K''$ is bijective, and
$\widetilde{\Omega}|_\K : \K \rightarrow \K'$ is an isomorphism.
Next, let
$$\mu : \mathfrak{B}(V^*) \otimes \K \to \K \# \mathfrak{B}(V^*)\;
\text{ and } \mu' : \K' \otimes \mathfrak{B}(V^*) \to
\mathfrak{B}(W')$$ be the multiplication maps. By Remark
\ref{re:multinv} resp. Lemma~\ref{basicK} (ii), both maps are
bijective. Let $f \in \mathfrak{B}(V^*)$ and $x \in \K$. Then
\begin{align*}
\widetilde{\Omega}(fx) &= \Ss_{\mathfrak{B}(W')}(f \Omega(x)) =
(f\_{-1} \cdot \widetilde{\Omega}(x)) \Ss_{\mathfrak{B}(W')}(f\_0)
\\& =(f\_{-1} \cdot \widetilde{\Omega}(x))
\Ss_{\mathfrak{B}(V^*)}(f\_0).
\end{align*}
Thus $\widetilde{\Omega} \mu   = \mu' c (\Ss_{\mathfrak{B}(V^*)}
\otimes (\Ss_{\toba(W')} \Omega|_\K))$. Hence
$\widetilde{\Omega}$, and \emph{a fortiori} $\Omega$, are
bijective. This completes the proof of Theorem \ref{theo:main}.
\epf

\begin{obs}
The proof of Theorem \ref{theo:main} does not use the fact that
  $V=M_i$ is irreducible in $\ydh $. However,
  $M_i$ has to be irreducible if one wants to apply the theorem
  for an index $j\in \Ib $, $j\not=i$, which satisfies Condition $(F_j)$.
\end{obs}

\bigbreak The algebras $\toba(W)$ and $\toba(W')$ are not
necessarily isomorphic. However, we have the following
consequences of Theorem \ref{theo:main}.

\begin{cor}\label{co:WW'supp}
  Let $M$, $i$, $W$, $W'$ and the $\zt $-gradings of $W$ and $W'$
  be as in Thm. ~\ref{theo:main}. Then
  $\toba(W)\#\toba(W^*)$ and $\toba(W')\#\toba(W'{}^*)$ are isomorphic
  as $\zt $-graded objects in $\ydh $. In particular,
  $\supp \toba(W)\#\toba(W^*)=\supp \toba(W')\#\toba(W'{}^*)$.
\end{cor}

\pf Since the homogeneous components of $\K$ are \fd, the graded
dual $\K^{\grd}$ of $\K$ is a $\zt $-graded object in $\ydh$. By
definition of $\K $ and the isomorphism $\toba(W^*)\simeq
\toba(W)^{\grd}$, see \eqref{eq:BVgrdual}, one has
$$\toba(W)\#\toba(W^*)\simeq \K \otimes \toba(V) \otimes \K ^{\grd}\otimes
\toba(V^*)$$ as $\zt $-graded objects in $\ydh $. Further, Theorem
\ref{theo:main} implies that
$$\toba(W')\#\toba(W'{}^*)\simeq \K \otimes \toba(V^*) \otimes \K ^{\grd}\otimes
\toba(V)$$
as $\zt $-graded objects in $\ydh $. Since
$A\ot B\simeq B\ot A$ for all $\zt$-graded objects $A,B$ in $\ydh $, the
above equations prove the corollary.
\epf

In many applications it will be more convenient to use the following
reformulation of Corollary~\ref{co:WW'supp}.

\begin{cor}\label{co:WW'supp2}
  Let $M$, $i$, $W$, $W'$ be as in Thm. ~\ref{theo:main}. Define
  $\zt $-gradings on $W$ and $W'$ by $\deg x=\alpha _j$ for all $x\in M_j$
  and all $x\in M'_j$, $j\in \Ib $.
  Then for all $\alpha \in \zt $ the homogeneous components
  $(\toba(W)\#\toba(W^*))_{\alpha }$ and $(\toba(W')\#\toba(W'{}^*))_{\sE
  (\alpha )}$
  are isomorphic in $\ydh $. In particular,
  $$\supp \toba(W')\#\toba(W'{}^*)=\sE (\supp \toba(W)\#\toba(W^*)).$$
\end{cor}

\begin{cor}\label{cor:dimfin} If $\dim \toba (W) < \infty$, then
$\dim \toba (W) = \dim\toba(W')$.
\end{cor}

\pf We compute $\dim \toba (W) \overset{} = \dim \K \dim \toba (V)
= \dim \K \dim \toba (V^*) = \dim\toba(W')$. Here the first
equality holds by Lemma~\ref{K0} (ii), the second by
Prop.~\ref{prop:duality}, and the third by
Theorem \ref{theo:main}. \epf

\subsection{Weyl groupoid and real roots}
\label{ss:Weylgroupoid}

In this subsection we define and study invariants of finite
families of finite dimensional irreducible Yetter-Drinfeld
modules. The definitions are based on Theorem \ref{theo:main}.


\medbreak

Recall the definition of $\C_\theta$ from Subsection
\ref{subsection:reflected}. If $M$, $M'\in \C_{\theta }$, then we
say that
$$M\sim M'$$ if there exists an index $i$ such that
Condition~$(F_i)$ holds for $M$, see Subsection
\ref{subsection:reflected}, and if $\rfl_i(M)\simeq M '$. By
Theorem \ref{theo:main}(2), the relation $\sim $ is symmetric. The
equivalence relation $\approx$ generated by $\sim$ is called
\textit{Weyl equivalence}.

Recall that $(\alpha _1,\ldots ,\alpha _\theta )$ is the standard basis of
$\zt $ and $\Ib =\{1,2,\dots ,\theta \}$.

\begin{definition}\label{def:rsys}
  Let $M \in \C_\theta $.
  Define
  \begin{align*}
    \W (M)=&\{M'\in \C_\theta \,|\,M'\approx M\}.
  \end{align*}
  Let $\Wg (M)$ denote the following category with
  $\mathrm{Ob}(\Wg (M))=\W (M)$.
  For each $M'\in \W (M)$ such that $M'$ satisfies $(F_i)$
  consider the reflection $s_{i,M'}\in \Aut (\zt )$,
  $s_{i,M'}(\alpha _j)=\alpha_j-a^{M'}_{ij}\alpha_i$ for all $j\in \Ib $,
  as a morphism
  $M'\to \rfl _i(M')$.
  Let $\Wg (M)$ be the category in which all morphisms are compositions of the
  morphisms $s_{i,M'}$, where $i\in \Ib $ and $M'\in \W (M)$ satisfies
  $(F_i)$. The category $\Wg (M)$ is called the \textit{Weyl groupoid of} $M$.

  Let
  \begin{align*}
    \rsys (M)=\{ w(\alpha _j)\,|\,w\in \Hom (M',M),M'\in \W (M)\}
  \subset \zt.
  \end{align*}
Following the notation in \cite[\S 5.1]{K},
  $\rsys (M)$ is called the
  \textit{set of real roots of} $M$.
\end{definition}

\begin{obs}
  Let $M\in \C_\theta $. Then the category $\Wg (M)$ is a connected
  groupoid.
  Indeed, if $i\in \Ib $ and $M'\in \W (M)$ satisfies $(F_i)$,
  then $\rfl _i(M')$ satisfies $(F_i)$, $\rfl _i(\rfl _i(M'))
  \simeq M'$ and $s_{i,\rfl _i(M')}s_{i,M'}=\id _{\zt }$
  by Theorem~\ref{theo:main}.
  Therefore the generating morphisms (and hence all morphisms)
  of $\Wg (M)$ are invertible. Further, for any two $M',M''\in \W (M)$
  there is a morphism in $\Hom (M',M'')$ by the definition of $\W (M)$.
\end{obs}

\begin{obs}
If
$M'\in\W(M)$, then equation
$$\dim\toba(\oplus _{n\in \Ib }M_n)
=\dim\toba(\oplus _{n\in \Ib }M'_n)$$
holds by Corollary~\ref{cor:dimfin}.
\end{obs}

\begin{obs} Assume that $M$ is of diagonal type,
  that is, $\dim M_n=1$ for all $n\in \Ib $.
  Let $M'\in \W(M)$, $j\in \Ib $, $k\in \mathbb{N}$, $i_1,\dots ,i_k\in \Ib $,
  $M^0=M',M^1,M^2,\dots ,M^k=M\in \W (M)$, $M^l\simeq \rfl
  _{i_{l+1}}(M^{l+1})$ for all $l\in \mathbb{N}_0$ with $l<k$,
  and $\beta =s_{i_k,M^{k-1}}\cdots s_{i_2,M^1}s_{i_1,M^0}(\alpha _j)\in \rsys
  (M)$. Let $m_1,\dots
  ,m_\theta \in \mathbb{Z}$ such that $\beta =\sum _{n\in \Ib }m_n\alpha _n$.
  Then $M'_j\simeq \otimes _{i\in \Ib }M_n^{\otimes m_n}$ in $\ydh $,
  where $M_n^{\otimes m_n}=(M_n^*)^{\otimes -m_n}$ if $m_n<0$.
  This means, that $M'$ depends only on $M$ and $\beta $, but not on the
  particular choice of $i_1, \dots ,i_k$.
  However, for more general Yetter-Drinfeld modules $M$ it is in general not
  clear, if $M'_j$ can be recovered from $j$, $\beta $ and $M$.
\end{obs}

\begin{prop}\label{pr:finrsys}
  Let $M\in \C_\theta $
  and $W=M_1\oplus \cdots \oplus M_\theta $. Then
  $\rsys (M )\subset \supp \toba(W)\#\toba (W^*)$. In particular, if
  $\toba (W) $ is \fd, then $\rsys (M )$ is a finite subset of $\zt $.
\end{prop}

\pf
Clearly, $\supp W =\{\alpha _1,\dots ,\alpha _\theta \}
\subset \supp \toba(W)\#\toba (W^*)$.
Let
$i\in \Ib $, $M '=\rfl_i(M )$, and
$W'=\oplus _{n\in \Ib }M'_n$.
Then
$$\supp W'\subset \supp \toba(W')\#\toba(W'{}^*)
=s_{i,M}(\supp \toba(W)\#\toba(W^*))$$
by Corollary~\ref{co:WW'supp2}.
Thus equation $s_{i,M'}\sE =\id $ gives that
$s_{i,M'}(\alpha _j)\in \supp \toba(W)\#\toba(W^*)$ for all $j\in \Ib $.
By iteration one obtains that
$\rsys (M )\subset \supp \toba(W)\#\toba (W^*)$.
If $\dim \toba(W)<\infty $ then
the finiteness of $\rsys (M )$ follows from the equations
\begin{align*}
\supp \toba (W)\otimes \toba(W^*)=
&\supp \toba (W)+\supp \toba(W)^{\grd}\\
=& \supp \toba (W)-\supp \toba(W),
\end{align*}
see \eqref{eq:BVgrdual}, and the fact that $\supp \toba (W)$ is
finite.
\epf

\begin{lema}\label{exa:fatqls} Let $M\in \C_{\theta }$
  and let $i,j\in \Ib $ such that $a^M_{ij} = 0$.
Then $a^M_{ji} = 0$, and $\toba(M_i\oplus M_j)
\simeq \toba(M_i)\ot \toba(M_j)$ as graded
vector spaces. \end{lema}

\pf Let $x\in M_i$, $y\in M_j$. Then
\eqref{eq:adcxy} gives that
\begin{equation}\label{eqn:fatqls}
\Delta(\ad_c x(y)) = \ad_c x(y)\ot 1 + x\ot y - c^2(x\ot y) +
1\ot \ad_c x (y).
\end{equation}
Thus, $a_{i j}^M=0$ implies that $\ad_c x (y) = 0$. Hence $x\ot y
- c^2(x\ot y) = 0$, that is, $(\id - c^2)(M_i\ot M_j) = 0$.
Then $(\id - c^2)c(M_j\ot M_i)=0$, but $c$ is invertible, so that $(\id -
c^2)(M_j\ot M_i)=0$. Eq.~\eqref{eqn:fatqls} gives
that $\ad _c x(y)$ is primitive in $\toba (M_i\oplus M_j)$
for all $x\in M_j$, $y\in M_i$, hence zero. This yields $a^M_{ji}=0$.
The last claim of the
lemma is \cite[Thm. \,2.2]{G-jalg}. \epf

\begin{lema}\label{lema:gen-cartan-matrix}
  Let $M = ( M_j)_ {j\in \Ib }$ be an object in
$\C_{\theta }$ which satisfies $(F_i)$ for all $i\in \Ib $.
Then $A = (a^M_{ij})$ is a generalized Cartan
matrix. In particular, the subgroup
$$\wo := \langle \sE\,|\,i\in \Ib\rangle$$ of $GL(\theta ,\Z )$ is
isomorphic to the Weyl group of the Kac-Moody algebra $\g(A)$.
\end{lema}

\pf The first claim follows from Lemma~\ref{exa:fatqls}. Let
$(\h, \Pi, \Pi^{\vee})$ be a realization of $A$ \cite[\S 1.1]{K}
and let $W$ be the Weyl group of $\g(A)$ \cite[\S 3.7]{K}. Then
$W$ preserves the subspace $V$ of $\h^*$ generated by $\Pi^{\vee}$
and the morphism $W\to GL(V)$ is injective \cite[Ex. 3.6]{K}. Now
$V \simeq \zt\ot_{\Cc} \Cc$ by \cite[(1.1.1)]{K} and the
image of $W$ in $GL(V)$ coincides with $\wo$ by \cite[(1.1.2)]{K}.
\epf

It follows that $\wo$ is a Coxeter group \cite[Prop. 3.13]{K} but
we do not need this fact in the sequel. The group $\wo$ is
important in the study of Nichols algebras in the following
special case.

\begin{definition}\label{defi:standard} We say that $M\in \C_{\theta }$
  is \emph{standard} if $M'$ satisfies Condition~$(F_i)$
   and $a^{M'}_{i j}=a^M_{i j}$ for all $M'\in \W (M)$ and
   $i,j\in \Ib $.
\end{definition}

\begin{obs}\label{rem:standard}
In the following two special cases the family $M\in \C_\theta $ is standard.

  1. Let $H$ be the group algebra of an abelian group $\Gamma $ and $M$ a
  family of 1-dimensional Yetter-Drinfeld modules $\ku v_i=M_i$
  over $H$, where $i\in \Ib $. Let $\delta (v_i)=g_i\ot v_i$ and
  $g\cdot v_i=\chi _i(g)v_i$ denote the coaction and action of $H$,
  respectively, where $g_i\in \Gamma $, $\chi _i\in \widehat{\Gamma }$,
  $i\in \Ib $. Define $q_{i j}=\chi _j(g_i)\in \ku $ for $i,j\in \Ib $.
  If $M$ is of Cartan type, that is, for all $i\not=j$ there exist
  $a_{i j}\in \Z$ such that $0\le -a_{i j}<\mathrm{ord}\,q_{i i}$ and
  $q_{i j}q_{j i}=q_{i i}^{a_{i j}}$, then $M$ is standard. This can be seen
  from \cite[Lemma\,1(ii), Eq.\,(24)]{H3}.

  2. Assume that $M\in \C_\theta $ satisfies Condition $(F_i)$
  and $\rfl _i(M)_{s_{i,M}(\alpha _j)}\simeq M_j$ in $\ydh $
  for all $i,j\in \Ib $. Then $M$ is standard by
  definition of $a^M_{i j}$.
\end{obs}


%

\begin{lema}\label{le:wrsys}
  Assume that $M\in \C_\theta $ is standard. Then
  $$\rsys (M)=\{w(\alpha _j)\,|\,w\in \wo ,j\in \Ib \}.$$
  In particular, $w(\rsys (M))=\rsys (M)$ for all $w\in \wo $.
\end{lema}

\pf This follows immediately from the definitions of $\rsys (M)$ and $\wo $
and the relations $a^{M'}_{ij}=a^M_{ij}$ for all $i,j\in \Ib $
and $M'\in \W (M)$.
\epf

\begin{theorem}\label{theorem:gpd-nichols-finite}
  Let $M=(M_i)_{i\in \Ib }\in \C_\theta $ and
  $W=\oplus _{i\in \Ib }M_i$.
  If $M$ is standard and $\dim \toba(W)<\infty $,
then the generalized Cartan matrix $(a^M_{ij})_{i,j\in \Ib }$
is of finite type.
\end{theorem}

\pf Since $\dim \toba(W)<\infty $, the set $\rsys (M)$ is finite
by Prop.~\ref{pr:finrsys}. Since $M$ is standard, $\rsys (M)$ is stable under
the action of $\wo $ by Lemma~\ref{le:wrsys}.
The corresponding permutation representation
$\wo \to S(\rsys (M))$ is
injective, since $\wo \subset GL(\theta ,\Z )$ and $\rsys (M)$ contains the
standard basis of $\Z ^\theta $. Therefore $\wo $ is finite. Thus the claim
follows from Lemma~\ref{lema:gen-cartan-matrix}
and \cite[Prop.\,4.9]{K}.
\epf

%
%
%
%

\section{Applications}\label{section:appl}

\subsection{Hopf algebras with few \fd{} Nichols algebras}


\begin{lema}\label{lema:one-irred-finite}
Let $H$ be a Hopf algebra. Assume that, up to isomorphism, there
is exactly one  \fd{} irreducible Yetter-Drinfeld module
$L\in \ydh $ such that $\dim \toba(L)<
\infty$.
Let $M = (M_1,M_2)\in \C_2$ such that $M_1\simeq M_2\simeq L$.

\begin{enumerate}
    \item[(i)] 
      If $M$ satisfies $(F_1)$ then $M$ satisfies $(F_2)$
      and $a^M_{12} = a^M_{21}$.
      If additionally $\dim \toba(L^2)<\infty $ then $M$ is standard.
    \item[(ii)] If 
      $M$ does not fulfill $(F_1)$ or if $a^M_{12} \leq -2$,
      then $\dim(\toba(L^n)) = \infty$ for $n\geq 2$.
    \item[(iii)] If $a^M_{12} =0$,
      then $\dim \toba(L^n) = (\dim \toba(L))^n$ for all $n\in \N$.
    \item[(iv)] When $a^M_{12} = -1$, then $\dim \toba(L^n) = \infty$
      for $n\geq 3$.
\end{enumerate}
\end{lema}

Note that if $a^M_{12}=-1$ then
Lemma~\ref{lema:one-irred-finite} gives no
information about $\dim \toba(L^2)$.

\pf
If 
$M$ does not fulfill Condition~$(F_1)$ then $\dim \toba(L^2) =
\infty$. Otherwise $a^M_{12}\in \Z_{\le 0}$, and $a^M_{12} =
a^M_{21}$ by symmetry. Moreover, if $\dim \toba (L^2)<\infty $,
then for $i\in \{1,2\}$ the Nichols algebra of
$\rfl_i(M)_1\oplus \rfl_i(M)_2$ is also \fd{} by
Corollary~\ref{cor:dimfin}, and hence $\rfl_i(M)_j\simeq L$ for
$j\in \{1,2\}$. Therefore $M$ is standard, and (i) is proven.

The generalized Cartan matrix $ \begin{pmatrix} 2 & a^M_{12} \\
a^M_{12} & 2 \end{pmatrix}$ is of finite type iff $a^M_{12} = 0$
or $a^M_{1 2}=-1$. Then (ii) follows from Theorem
\ref{theorem:gpd-nichols-finite}. Now (iii) follows from
\cite{G-jalg}, see Lemma~\ref{exa:fatqls}. If $a^M_{12} = -1$,
then the generalized Cartan matrix of $L^3$ has Dynkin diagram
$A_2^{(1)}$; hence $\dim \toba(L^3) = \infty$, and \emph{a
fortiori} the same holds for $L^n$ for $n\ge 3$. This shows (iv).
\epf

\begin{theorem}\label{theo:one-irred-finite}
Let $H$ be a Hopf algebra such that the category of \fd{}
Yetter-Drinfeld modules is semisimple. Assume that up to
isomorphism there is exactly one irreducible $L\in \ydh$ such that
$\dim \toba(L)< \infty$. Let $M = (M_1,M_2)\in \C_2$,
where $M_1=M_2=L$. If $M$ satisfies $(F_1)$ then $M$
satisfies $(F_2)$ and $ a^M_{12} = a^M_{21}\in \Z_{\le 0}$.
%

\begin{enumerate}
    \item[(i)] If $a^M_{12} = - \infty$ or $a^M_{12} \leq -2$,
      then $L$ is the only Yetter-Drinfeld module over $H$ with
      \fd{} Nichols algebra.
    \item[(ii)] If $a^M_{12} =0$, then
      a Yetter-Drinfeld module $W$ over $H$ has \fd{} Nichols algebra if
      and only if $W\simeq L^n$ for some $n\in \N$.
      Furthermore,
    $\dim \toba(L^n) = (\dim \toba(L))^n$.
    \item[(iii)]  If $a^M_{12} = -1$, then the only possible
      Yetter-Drinfeld modules over $H$ with \fd{}
    Nichols algebra are $L$ and (perhaps) $L^2$.
\end{enumerate}
\end{theorem}

\pf By hypothesis, the only Yetter-Drinfeld module candidates to
have \fd{} Nichols algebras are those of the form
$L^n$, $n\in \N$. The theorem follows then from
Lemma~\ref{lema:one-irred-finite}. \epf

Now we state another general result that can be obtained from
Theorem \ref{theorem:gpd-nichols-finite}. We shall use it when
considering Nichols algebras over $\sk$.

\begin{lema}\label{lema:one-twocopies-finite}
  Let $M_1, \dots, M_s\in \ydh $,
  where $s\in \N $, be a maximal set of pairwise
  nonisomorphic irreducible Yetter-Drinfeld modules, such that
  $\dim \toba(M_i)< \infty$ for $1\le i \le s$.
  Assume that there exist $i, j \in \{1, \dots, s\}$
  (the possibility $i=j$ is not excluded) such that
\begin{enumerate}
    \item[(i)] $\dim \toba(M_i\oplus M_j) < \infty$.
    \item[(ii)] If $\{\ell,m\}\neq \{i,j\}$, then
      $\dim \toba(M_\ell\oplus M_m) = \infty$.
    \item[(iii)] $M_i \not\simeq M_j^*$.
\end{enumerate}
Let $M = (M_i, M_j)\in \C_2$.
Then $M$ is standard.
\end{lema}

\pf By (i) the Nichols algebra of $(M_i\oplus M_j)^*\simeq
M_i^*\oplus M_j^*$ is \fd. By (ii) one has $M_i^*\oplus
M_j^*\simeq M_i\oplus M_j$, and (iii) implies that $M_i^*\simeq
M_i$ and $M_j^*\simeq M_j$. Thus it suffices to consider the
reflection $\rfl _i$. By (i) and Corollary~\ref{cor:dimfin}, $M
'=(M'_1,M'_2):=\rfl _i(M )$ is well-defined and
$\dim \toba(M'_1\oplus M'_2)<\infty $. By (ii)
one has $M'_1\oplus M'_2\simeq M_i\oplus M_j$.
Since $M'_1\simeq M_i^*\simeq M_i$ by the
beginning of the proof, one has $M'_2\simeq M_j$. Hence
$M$ is standard by Remark~\ref{rem:standard}. \epf
%

\subsection{Pointed Hopf algebras with group $\st$}
In the rest of this section, it is
assumed that the base field is $\ku =\Cc $.
Let $G$ be a finite non-abelian group. We shall use the rack
notation $x\trid y := xyx^{-1}$, $x,y\in G$. Since the group algebra
$\Cc G$ is semisimple, the corresponding category $\ydg$ of Yetter-Drinfeld modules
is semisimple. It is well-known that the irreducible objects in $\ydg$ are
parametrized by pairs $(\Oc, \rho)$, where $\Oc$ is a conjugacy class of $G$ and
$\rho$ is an
irreducible representation of the centralizer $G^s$ of a fixed
$s\in \Oc$. Let $M(\Oc, \rho)$ denote the irreducible Yetter-Drinfeld
module corresponding to $(\Oc, \rho)$ and let $\toba(\Oc, \rho)$
be its Nichols algebra. Then $M(\Oc, \rho)$ is the induced module
$\Ind_{G^s}^G \rho$, and the comodule structure is given by the following rule.
Let $g_1 = g$, \dots, $g_{t}$ be a numeration of $\Oc$ and let
$x_i\in G$ such that $x_i \trid g = g_i$ for all $1\le i \le t$.
Then $M(\Oc, \rho) = \oplus_{1\le i \le t} x_i\otimes V$. If
$x_iv := x_i\otimes v \in M(\Oc,\rho)$,
then $\delta(x_iv) = g_i\ot x_iv$, for  $1\le i \le t$, $v\in V$.
The braiding in $M(\Oc, \rho)$ is given by
$c(x_iv\otimes x_jw) = g_i\cdot(x_jw)\otimes x_iv =
x_h\,\rho(\gamma)(w) \otimes x_iv$ for any $1\le
i,j\le t$, $v,w\in V$, where $g_ix_j = x_h\gamma$ for unique $h$,
$1\le h \le t$ and $\gamma \in G^s$.

\medbreak If $G = \sn$, then $\Oc^n_2$ is the conjugacy class of
the involutions and $\sgn$ is the restriction of the sign
representation to the isotropy group.

\medbreak  Before stating our first classification result, we need
to recall the construction of some Hopf algebras from
\cite{AG-ama}.

\begin{definition}\label{defi:abofet} Let $\lambda\in \ku$.
Let $\abofet$ be the algebra presented by generators $e_t$, $t\in
T:= \{(12), (23)\}$, and $a_\sigma$, $\sigma \in \Oc_2^3$; with
relations
\begin{align}
\label{relgpo} e_t e_se_{t} &= e_se_{t}e_{s}, \quad e_t^2 = 1,
\quad s\neq t\in T; \\
\label{relsym4} e_t a_\sigma &= -a_{t\sigma t} e_t \hspace{55pt} t
\in T,\, \sigma\in \Oc_2^3;
\\\label{relsym41}
a_\sigma^2 &= 0, \hspace{80pt} \sigma\in \Oc_2^3;
\\\label{relsym43}
a_{(12)} a_{(23)} &+ a_{(23)} a_{(13)} +a_{(13)} a_{(12)} =
\lambda (1- e_{(12)} e_{(23)});
\\\label{relsym43bis}
a_{(12)} a_{(13)} &+ a_{(13)} a_{(23)} +a_{(23)} a_{(12)} =
\lambda (1- e_{(23)} e_{(12)}).
\end{align}
Set $e_{(13)} = e_{(12)}e_{(23)}e_{(12)}$. Then $\abofet$ is a
Hopf algebra of dimension $72$ with comultiplication determined by
\begin{equation}\label{eqn:abofet-comult}
\Delta(a_{\sigma}) = a_{\sigma}\ot 1 + e_{\sigma}\ot a_{\sigma},
\quad \Delta(e_{t}) = e_{t}\ot e_{t}, \qquad \sigma\in \Oc_2^3,
t\in T.
\end{equation}

Observe that the Hopf algebra $\abofet$ is isomorphic to $\ac(\st,
\Oc_2^3, \lambda c^2)$ (via $a_{\sigma} \mapsto
c^{-1}a'_{\sigma}$, where $a'_{\sigma}$ are the generators of the
latter). Also $\ac(\st, \Oc_2^3, 0)$ $\simeq\toba(\Oc_2^3, \sgn)
\#\ku \st$. But $\ac(\st, \Oc_2^3, 0) \not \simeq \ac(\st,
\Oc_2^3, 1)$ since the former is self-dual but the latter is not.
\end{definition}

\begin{theorem}\label{theo:s3}
  Let $H$ be a \fd{} pointed Hopf algebra with
$G(H)\simeq \st$. Then either $H\simeq \ku \st$, or $H\simeq
\toba(\Oc_2^3, \sgn) \#\ku \st$ or $H\simeq \ac(\st, \Oc_2^3, 1)$.
\end{theorem}

\pf It is known that $\dim \toba(\Oc_2^3, \sgn) = 12$ \cite{MS};
it is also known that this is the only \fd{} Nichols algebra with
irreducible Yetter-Drinfeld module of primitives \cite{AZ}. We can
then apply Theorem \ref{theo:one-irred-finite}. Let $M =
(M(\Oc_2^3, \sgn), M(\Oc_2^3, \sgn))$. Assume that $a^M_{12}\in
\Z_{\le 0}$, notation as above. We claim that $-a^M_{12}\ge 2$.

Let $\sigma_1 = (12)$, $\sigma_2 = (23)$, $\sigma_3 = (13)\in
\st$. The Yetter-Drinfeld module $M(\Oc_2^3, \sgn) \oplus
M(\Oc_2^3, \sgn)$ has a basis $x_1$, $x_2$, $x_3$ (from the first
copy), $y_1$, $y_2$, $y_3$ (from the second copy) with
\begin{equation}\label{eqn:oc32}
\delta(x_i) = \sigma_i \otimes x_i, \, \delta(y_i) = \sigma_i
\otimes y_i, \,  t\cdot x_i = \sgn(t)x_{t \trid i}, \, t\cdot y_i
= \sgn(t)y_{t \trid i}, \end{equation} for $1\le i \le 3$, $ t \in
\st$. Here $\sigma_{t \trid i} := t\trid \sigma_i =
t\sigma_it^{-1}$. Also, $j\trid i$ means $\sigma_{j \trid i} :=
\sigma_j \trid \sigma_i$. The braiding in the vectors of the basis
gives
\begin{align*}
c(x_j\otimes x_i) &= - x_{j\trid i}\otimes  x_j,&\,
c(y_j\otimes y_i) &= - y_{j\trid i} \otimes y_{j},\\
c(x_j\otimes y_i) &= - y_{j\trid i} \otimes x_j, &\,
c(y_j\otimes x_i) &= - x_{j\trid i} \otimes y_j.
\end{align*}
To prove our claim, we need to find
$i,j,k$ such that $\ad_c(x_i)(\ad_c(x_j) (y_k))\neq 0$. Let
$\partial_{x_i}$, $\partial_{y_i}$ be the skew-derivations as in
\cite{MS}. Now
\begin{align*}
\ad_c(x_2)(\ad_c(x_1) (y_2)) &= \ad_c(x_2)(x_1y_2 + y_3x_1) \\
&= x_2x_1y_2 + x_2y_3x_1 - x_3y_2x_2 - y_1x_3x_2,
\end{align*}
hence $\partial_{x_3} \partial_{y_1}\left(\ad_c(x_2)(\ad_c(x_1)
(y_2))\right) = \partial_{x_3}\left(-x_2x_3\right) = -x_2 \neq 0$,
and the claim is proved. Thus $\dim \toba(M(\Oc_2^3, \sgn) \oplus
M(\Oc_2^3, \sgn)) = \infty$ by Theorem
\ref{theo:one-irred-finite}, and $\toba(\Oc_2^3, \sgn)$ is the
only \fd{} Nichols algebra over $\st$.

\medbreak Let $H\not\simeq \ku \st$ be a \fd{} pointed Hopf
algebra with $G(H)\simeq \st$. Then the infinitesimal braiding of
$H$, see \cite{AS-cambr}, is isomorphic to $M(\Oc_2^3, \sgn)$.
Hence $H$ is generated as algebra by group-like and skew-primitive
elements \cite[Theorem \,2.1]{AG-ama} and the theorem follows from
\cite[Thm.\,3.8]{AG-ama}. \epf

\subsection{Nichols algebras over the group $\sk$}

Let us recall the general terminology for $\mathbb S_n$.
If $\pi=(12)\in \Oc_2^n$, then the isotropy subgroup is
$\sn^{\pi}\simeq \Z_2\times \mathbb S_{n-2}$. Any irreducible
representation of $\sn^{\pi}$ is of the form $\chi\ot \rho$, where
$\chi\in \widehat{\Z_2}$, $\rho\in \widehat{\mathbb S_{n-2}}$. If
$\chi = \cou$, then $\chi\ot \rho(\pi) = 1$ and $\dim
\toba(\Oc_2^n, \cou\ot \rho) = \infty$. Thus, we are
actually interested in the Nichols algebras $\toba(\Oc_2^n,
\sgn\ot \rho)$. If $\rho = \sgn$, then $\sgn\ot \rho$ is just the
restriction to $\sn^\pi$ of the sign representation of $\sn$; we
denote in this case $\toba(\Oc_2^n, \sgn) = \toba(\Oc_2^n, \sgn\ot
\sgn)$.

\medbreak The proof of Theorem \ref{theo:s3} gives the following
result.

\begin{lema}\label{lema:o2n}
The Nichols algebras $\toba\big(M(\Oc_2^n, \sgn\ot \rho)\oplus
M(\Oc_2^n, \sgn\ot \rho')\big)$, $n\ge 4$, $\rho, \rho'\in
\widehat{\mathbb S_{n-2}}$, have infinite dimension.
\end{lema}

\pf The braided vector space $M(\Oc_2^3, \sgn) \oplus M(\Oc_2^3,
\sgn)$ is a braided subspace of any of these braided vector
spaces. \epf

The isotropy group of the 4-cycle $(1234)$ in $\sk$ is the cyclic
group $\langle(1234)\rangle$. Let $\chi_-$ be its character
defined by $\chi_-(1234) = -1$. Let $\Oc_4^4$ be the conjugacy
class of 4-cycles in $\sk$.

\begin{theorem}\label{thm:s4}
The only Nichols algebras of Yetter-Drinfeld modules over $\sk$
with finite dimension, up to isomorphism, are those in the
following list. All of them have dimension 576.

\begin{enumerate}
    \item $\toba(\Oc_2^4, \sgn)$.
    \item $\toba(\Oc_2^4, \sgn\ot \cou)$.
    \item $\toba(\Oc_4^4, \chi_{-})$.
\end{enumerate}
\end{theorem}

\pf The Nichols algebras in the list have the claimed dimension by
\cite{FK, MS, AG-adv}, respectively. These are the only Nichols
algebras of irreducible Yetter-Drinfeld modules over $\sk$ with
finite dimension by \cite{AZ}.

It remains to show: If $M$, $M'$ are two of $M(\Oc_2^4, \sgn)$,
$M(\Oc_2^4, \sgn\ot \cou)$, $M(\Oc_4^4, \chi_{-})$, then
$\dim \toba(M\oplus M') = \infty$. Some possibilities are covered
by Lemma~\ref{lema:o2n}. The rest are:
\begin{enumerate}
    \item[(i)] $\toba(M(\Oc_4^4, \chi_{-}) \oplus M(\Oc_4^4, \chi_{-}))$.
    \item[(ii)] $\toba(M(\Oc_2^4, \sgn)\oplus M(\Oc_4^4, \chi_{-}))$.
    \item[(iii)] $\toba(M(\Oc_2^4, \sgn\ot \cou)\oplus M(\Oc_4^4, \chi_{-}))$.
\end{enumerate}

(i). We claim that there is a surjective rack homomorphism
$\Oc_4^4 \to \Oc_2^3$ that induces a surjective morphism of
braided vector spaces $M(\Oc_4^4, \chi_{-}) \oplus M(\Oc_4^4,
\chi_{-}) \to M(\Oc_2^3, \sgn) \oplus M(\Oc_2^3, \sgn)$; since the
Nichols algebra of the latter is \infd{} by the proof of Theorem
\ref{theo:s3}, $\dim\toba(M(\Oc_4^4, \chi_{-}) \oplus M(\Oc_4^4,
\chi_{-})) = \infty$ too. Let us now verify the claim. We numerate
the elements in the orbit $\Oc_4^4$ as follows:
\begin{align*}
\tau_1 &= (1234), &\tau_3 &= (1243), &\tau_5 &= (1324),
\\ \tau_2 &= (1432) = \tau_1^{-1}, &\tau_4 &= (1342) = \tau_3^{-1},
&\tau_6 &= (1423) = \tau_5^{-1};
\end{align*}
set accordingly
\begin{align*}
h_1 &= \tau_1,  &h_2 &= (24), &h_3 &= \tau_6, &h_4 &= \tau_5, &h_5
&= \tau_3, & \quad h_6 &= \tau_4;
\end{align*}
so that $h_i\trid \tau_1 = \tau_i$, $1\le i \le 6$. The
Yetter-Drinfeld module $M(\Oc_4^4, \chi_{-}) \oplus M(\Oc_4^4,
\chi_{-})$ has a basis $u_1, \dots, u_6$ (from the first copy),
$w_1, \dots, w_6$ (from the second copy) with
\begin{equation}\begin{aligned}\label{eqn:deltaoc44} \delta(u_i) &=
\tau_i \otimes u_i, & t\cdot u_i &= \chi_{-}(\widetilde t)u_{t
\trid i},
\\ \delta(w_i) &= \tau_i \otimes w_i,  & t\cdot w_i &=
\chi_{-}(\widetilde t)w_{t \trid i},
\end{aligned}\end{equation}  for $1\le i \le 6$, $ t \in \sk$. Here $t\trid i$
and $\widetilde t \in \sk^{\tau_1} = \langle\tau_1\rangle$ have
the meaning that $t h_i = h_{t \trid i}\widetilde t$. Let now
\begin{align}\label{eqn:defIa}
I_3 &= \{1,2\}, & I_2 &= \{3,4\}, & I_1 &= \{5,6\}.
\end{align}
Let $a,b\in \{1,2,3\}$, $i\in I_a$, $j \in I_b$. If $a=b$, then
the braiding in the corresponding vectors of the basis is
\begin{align*}
c(u_i\otimes u_{j}) &= - u_{j}\otimes u_{i}, & c(u_i\otimes
w_{j}) &= - w_{j}\otimes u_{i}, \\
c(w_i\otimes w_{j}) &= - w_{j}\otimes w_{i}, & c(w_i\otimes u_{j})
&= - u_{j}\otimes w_{i};
\end{align*}
and if $a\neq b$, then for some $\ell\in I_c$, where $c\neq a,b$,
one has
\begin{align*}c(u_i\otimes u_{j}) &= - u_{\ell}\otimes u_{i}, & c(u_i\otimes
w_{j}) &= - w_{\ell}\otimes u_{i}, \\
c(w_i\otimes w_{j}) &= - w_{\ell}\otimes w_{i}, & c(w_i\otimes
u_{j}) &= - u_{\ell}\otimes w_{i}.
\end{align*}

Thus, the map $\pi: M(\Oc_4^4, \chi_{-}) \oplus M(\Oc_4^4,
\chi_{-}) \to M(\Oc_2^3, \sgn) \oplus M(\Oc_2^3, \sgn)$ given by
$\pi(u_i) = x_a$, $\pi(x_i) = y_a$, for $i\in I_a$, $a=1,2,3$,
preserves the braiding. This proves the claim.

(ii). We claim that there is a surjective morphism of
braided vector spaces $ M(\Oc_2^4, \sgn)\oplus M(\Oc_4^4,
\chi_{-}) \to M(\Oc_2^3, \sgn) \oplus M(\Oc_2^3, \sgn)$. Again,
this implies that the Nichols algebra in (ii) is \infd.
Let us check the claim. We numerate the elements in
the orbit $\Oc_2^4$ as follows:
\begin{align*}
\sigma_1 &= (12), &\sigma_2 &= (23), &\sigma_3 &= (13), & \sigma_4
&= (14), &\sigma_5 &= (24), &\sigma_6 &= (34);
\end{align*}
set accordingly
\begin{align*}
g_1 &= \sigma_1,  &g_2 &= \sigma_3, &g_3 &= \sigma_2, &g_4 &=
\sigma_5, &g_5 &= \sigma_4, & \quad g_6 &= (1324);
\end{align*}
so that $g_i\trid \sigma_1 = \sigma_i$, $1\le i \le 6$. Let
$\tau_i$ and $h_i$, $1\le i \le 6$, be as in the previous part of the
proof. The Yetter-Drinfeld module $M(\Oc_4^4, \sgn) \oplus
M(\Oc_4^4, \chi_{-})$ has a basis $z_1, \dots, z_6$ (from the
first summand), $w_1, \dots, w_6$ (from the second summand) with
\begin{equation}\begin{aligned}\label{eqn:deltaoc42} \delta(z_i) &=
\sigma_i \otimes z_i, & t\cdot z_i &= \sgn (t')z_{t \trid i},
\\ \delta(w_i) &= \tau_i \otimes w_i,  & t\cdot w_i &=
\chi_{-}(\widetilde t)w_{t \trid i},
\end{aligned}\end{equation}  for $1\le i \le 6$, $ t \in \sk$. Here,
in the first line $t\trid i$ and $t' \in \sk^{\sigma_1}$ have the
meaning that $t  g_i = g_{t \trid i}t'$; and in the second line,
$t\trid i$ and $\widetilde t \in \sk^{\tau_1} =
\langle\tau_1\rangle$ have the meaning that $t  h_i = h_{t \trid
i}\widetilde t$. Set $ t_1 :=\sigma_1$, $t_2 :=\sigma_6$, so that
$\sk^{\sigma_1} = \langle t_1, t_2 \rangle$.

\medbreak Let now $I_a$ be as in \eqref{eqn:defIa} and let $J_1 =
\{1, 6\}$, $J_2 = \{2,4\}$, $J_3 = \{3,5\}$. Let $a,b, c\in
\{1,2,3\}$ such that $\sigma_a\trid \sigma_b = \sigma_c$. Let
$i\in I_a$, $j \in I_b$. Then there exist $k\in I_c$, $\ell,m\in
J_c$, $\epsilon \in \{\pm 1\}$, $p,q\in \{1,2\}$ such that
\begin{align*}
\sigma_ih_j &= h_k\tau_1^{\epsilon},&\sigma_i g_j &= g_\ell t_p,
&\tau_i g_j &= g_m t_q;
\end{align*}
see the Table \ref{tab:1}.

\begin{table}[ht]
\caption{Multiplication in $\sk$.}\label{tab:1}
\begin{center}
\begin{tabular}{r|cccccc}
  $\cdot $ & $h_1$ & $h_2$ & $h_3$ & $h_4$ & $h_5$ & $h_6$\\
  \hline
  $\sigma _1$ & $h_4\tau_1^{-1}$ & $h_3\tau_1^{-1}$ & $h_2\tau _1$ &
  $h_1\tau_1$ & $h_6\tau_1$ & $h_5\tau _1^{-1}$ \\
  $\sigma_2$ & $h_5\tau_1^{-1}$ & $h_6\tau_1$ & $h_4\tau _1$ &
  $h_3\tau_1^{-1}$ & $h_1\tau_1$ & $h_2\tau _1^{-1}$ \\
  $\sigma_3$ & $h_2\tau_1$ & $h_1\tau _1$ & $h_5\tau _1^{-1}$ &
  $h_6\tau_1$ & $h_3\tau_1^{-1}$ & $h_4\tau _1^{-1}$ \\
  $\sigma _4$ & $h_6\tau_1^{-1}$ & $h_5\tau_1$ & $h_4\tau _1^{-1}$ &
  $h_3\tau_1$ & $h_2\tau_1^{-1}$ & $h_1\tau _1$ \\
  $\sigma _5$ & $h_2\tau_1$ & $h_1\tau_1^{-1}$ & $h_6\tau _1^{-1}$ &
  $h_5\tau_1^{-1}$ & $h_4\tau_1$ & $h_3\tau _1$ \\
  $\sigma _6$ & $h_3\tau_1^{-1}$ & $h_4\tau_1^{-1}$ & $h_1\tau _1$ &
  $h_2\tau_1$ & $h_6\tau_1^{-1}$ & $h_5\tau _1$
\end{tabular}
\end{center}

\begin{center}
\begin{tabular}{r|cccccc}
  $\cdot $ & $g_1$ & $g_2$ & $g_3$ & $g_4$ & $g_5$ & $g_6$\\
  \hline
  $\sigma _1$ & $g_1 t_1$ & $g_3 t_1$ & $g_2 t_1$ &
  $g_5 t_1$ & $g_4 t_1$ & $g_6 t_2$ \\
  $\sigma_2$ & $g_3 t_1$ & $g_2 t_1$ & $g_1 t_1$ &
  $g_4 t_2$ & $g_6 t_1$ & $g_5 t_1$ \\
  $\sigma_3$ & $g_2 t_1$ & $g_1 t_1$ & $g_3 t_1$ &
  $g_6 t_2$ & $g_5 t_2$ & $g_4 t_2$ \\
  $\sigma _4$ & $g_5 t_1$ & $g_2 t_2$ & $g_6 t_1$ &
  $g_4 t_1$ & $g_1 t_1$ & $g_3 t_1$ \\
  $\sigma _5$ & $g_4 t_1$ & $g_6 t_2$ & $g_3 t_2$ &
  $g_1 t_1$ & $g_5 t_1$ & $g_2 t_2$ \\
  $\sigma _6$ & $g_1 t_2$ & $g_5 t_2$ & $g_4 t_2$ &
  $g_3 t_2$ & $g_2 t_2$ & $g_6 t_1$
\end{tabular}
\end{center}

\begin{center}
\begin{tabular}{r|cccccc}
  $\cdot $ & $g_1$ & $g_2$ & $g_3$ & $g_4$ & $g_5$ & $g_6$\\
  \hline
  $\tau _1$ & $g_2 t_2$ & $g_6 t_1$ & $g_5 t_1$ &
  $g_1 t_2$ & $g_3 t_2$ & $g_4 t_1$ \\
  $\tau _2$ & $g_4 t_2$ & $g_1 t_2$ & $g_5 t_2$ &
  $g_6 t_1$ & $g_3 t_1$ & $g_2 t_1$ \\
  $\tau _3$ & $g_5 t_2$ & $g_4 t_2$ & $g_1 t_2$ &
  $g_2 t_1$ & $g_6 t_2$ & $g_3 t_2$ \\
  $\tau _4$ & $g_3 t_2$ & $g_4 t_1$ & $g_6 t_2$ &
  $g_2 t_2$ & $g_1 t_2$ & $g_5 t_2$ \\
  $\tau _5$ & $g_6 t_1$ & $g_5 t_1$ & $g_2 t_2$ &
  $g_3 t_1$ & $g_4 t_2$ & $g_1 t_2$ \\
  $\tau _6$ & $g_6 t_2$ & $g_3 t_2$ & $g_4 t_1$ &
  $g_5 t_2$ & $g_2 t_1$ & $g_1 t_1$
\end{tabular}
\end{center}
\end{table}

Hence, the braiding in the  vectors of the basis is
\begin{align*}
c(z_i\otimes w_{j}) &= - w_{k}\otimes z_{i}, & c(z_i\otimes z_{j})
&= - z_{\ell}\otimes z_{i}, & c(w_i\otimes z_{j}) &= -
z_{m}\otimes w_{i}.
\end{align*}
Thus, the map $\pi: M(\Oc_2^4, \sgn) \oplus M(\Oc_4^4, \chi_{-})
\to M(\Oc_2^3, \sgn) \oplus M(\Oc_2^3, \sgn)$ given by $\pi(z_i) =
x_a$, $\pi(w_j) = y_a$, for $i\in I_a$, $j\in J_a$, $a=1,2,3$,
preserves the braiding. This proves the claim.

\medbreak  (iii). The argument in the preceding part can not be
adapted to this one. However, assume that $\dim \toba(M(\Oc_2^4,
\sgn\ot \cou)\oplus M(\Oc_4^4, \chi_{-}))< \infty$. Then
$M(\Oc_2^4, \sgn\ot \cou)\oplus M(\Oc_4^4, \chi_{-})$ is standard
with finite Cartan matrix $(a_{ij})$, by
Lemma~\ref{lema:one-twocopies-finite}. Let $\sigma_i$ and $g_i$,
$\tau_i$ and $h_i$, $1\le i \le 6$, be as in previous part of the
proof. The Yetter-Drinfeld module $M(\Oc_4^4, \sgn\ot \cou) \oplus
M(\Oc_4^4, \chi_{-})$ has a basis $\widetilde z_1, \dots,
\widetilde z_6$ (from the first summand), $w_1, \dots, w_6$ (from
the second summand) with action and coaction given by
$\delta(\widetilde z_i) = \sigma_i \otimes \widetilde z_i$,
$t\cdot \widetilde z_i = (\sgn \ot \cou)(t')\widetilde z_{t \trid
i}$ for $1\le i \le 6$, $ t \in \sk$, and the second line of
\eqref{eqn:deltaoc42}. Here, $t\trid i$ and $t' \in
\sk^{\sigma_1}$ have the meaning that $t  g_i = g_{t \trid i}t'$.
Then
\begin{align*}
\ad (\widetilde z_2)(\ad (\widetilde z_1) (w_1)) &= \widetilde z_2
\widetilde z_1 w_1 + \widetilde z_2 w_4\widetilde z_1 - \widetilde
z_3 w_5\widetilde z_2 - w_3\widetilde z_3 \widetilde z_2 \neq 0\\
\text{since }  \partial_{\widetilde z_1} &\partial_{w_1}
\big(\ad (\widetilde z_2)(\ad (\widetilde z_1) (w_1)) \big) =
\partial_{\widetilde z_1}(\widetilde z_2 \widetilde z_1) = \widetilde
z_2 \neq 0;\\
\ad (w_2)(\ad (w_1) (\widetilde z_1)) &= w_2 w_1\widetilde  z_1 - w_2
\widetilde z_2w_1 + w_1 \widetilde z_4w_2 - \widetilde z_1w_1 w_2 \neq 0\\
\text{since }  \partial_{w_5} \partial_{\widetilde z_2} &\big(
\ad (w_2)(\ad (w_1) (\widetilde z_1)) \big) = \partial_{w_5}(w_2 w_5) =
w_2 \neq 0.
\end{align*}
Hence $a_{12} \leq -2$, $a_{21} \leq -2$, a contradiction. Thus,
$\dim \toba(M(\Oc_2^4, \sgn\ot \cou)\oplus M(\Oc_4^4,
\chi_{-}))= \infty$.  \epf

\subsection{Nichols algebras over the group $\dn$, $n$ odd}

Let $n>1$ be an odd integer and let $\dn$ be the dihedral group of
order $2n$, generated by $x$ and $y$ with defining relations $x^2
= e = y^n$ and $xyx = y^{-1}$. Let $\Oc$ be a conjugacy class of
$\dn$ and let $\rho$ be an irreducible representation of the
centralizer $G^s$ of a fixed $s\in \Oc$.

By \cite[Th. 3.1]{AF2}, we know that there is at most one pair
$(\Oc, \rho)$ such that the Nichols algebra $\toba(\Oc, \rho)$ is
finite-dimensional, namely $(\Oc, \rho) = (\Oc_{x}, \sgn)$, where
$\sgn \in \widehat{\dn^x}$, $\dn^x=\langle x \rangle \simeq
\mathbb Z_2$. However, it is not known if the dimension of
$\toba(\Oc_{x}, \sgn)$ is finite, except when $n= 3$-- since $\dt
\simeq \st$.

\smallbreak The next result generalizes the first part of the
proof of Theorem \ref{theo:s3}.

\begin{theorem}\label{theorem:dn}
The only possible Nichols algebra over $\dn$ with finite
dimension, up to isomorphism, is $\toba(\Oc_{x}, \sgn)$.
\end{theorem}

\pf If $\dim\toba(\Oc_{x}, \sgn) = \infty$, then there is no
finite-dimensional Nichols algebra over $\dn$. Otherwise, we can
apply Theorem \ref{theo:one-irred-finite}. Let $\Mb = M(\Oc_{x},
\sgn) \oplus M(\Oc_{x}, \sgn)$. Assume that $a_{12}\in \Z_{\le
0}$, notation as above. We claim that $-a_{12}\ge 2$. Let
$\sigma_i = xy^i\in \dn$; $\Oc_{x} = \{\sigma_i\,|\,i\in \Z_n\}$. The
Yetter-Drinfeld module $\Mb$ has a basis $v_i$, $i\in \Z_n$ (from
the first copy), $w_i$, $i\in \Z_n$ (from the second copy) with
action, coaction and braiding

\begin{align*}
t\cdot v_i &= \sgn(t)v_{t \trid i}, &\, t\cdot w_i&= \sgn(t)w_{t \trid i}, \\
\delta(v_i) &= \sigma_i \otimes v_i, &\, \delta(w_i) &= \sigma_i
\otimes w_i, \\
c(v_j\otimes v_i) &= - v_{j\trid i}\otimes  v_j,&\, c(w_j\otimes
w_i) &= - w_{j\trid i} \otimes w_{j},\\ c(v_j\otimes w_i) &= -
w_{j\trid i} \otimes v_j, &\, c(w_j\otimes v_i) &= - v_{j\trid i}
\otimes w_j.\end{align*} for $i,j\in \Z_n$, $ t \in \dn$. Here, as
above, $\sigma_{t \trid i} := t\trid \sigma_i = t\sigma_it^{-1}$.
To prove our claim, we need to find $i,j,k$ such that
$\ad_c(v_i)\ad_c(v_j) (w_k)\neq 0$. Let $\partial_{v_i}$,
$\partial_{w_i}$ be the skew-derivations as in \cite{MS}. Now
\begin{align*}
\ad_c(v_2)\ad_c(v_1) (w_2) &= \ad_c(v_2)(v_1w_2 + w_0v_1) \\
&= v_2v_1w_2 + v_2w_0v_1 - v_3w_2v_2 - w_4v_3v_2,
\end{align*}
hence $\partial_{v_6} \partial_{w_4}\left(\ad_c(v_2)\ad_c(v_1)
(w_2)\right) = \partial_{v_6}\left(-v_5v_6\right) = -v_5 \neq 0$.
The claim and the theorem are proved. \epf

\bigbreak
\subsection*{Acknowledgements}
Results of this paper were obtained during a visit of I. H. and
H.-J. S. at the University of C\'ordoba in March 2006, partially
supported through grants of CONICET and Foncyt. I. H. and H.-J. S.
thank the University of C\'ordoba, CONICET and Foncyt for the kind
hospitality and the financial support. N. A. is grateful to Jorge
Vargas for conversations on Coxeter groups; and to Fernando
Fantino for discussions on $M(\Oc_4^4, \chi_{-})$.

\end{document}